\documentclass[a4paper]{article}

\usepackage{amsmath,amssymb} 
\usepackage{amsfonts}
\usepackage{graphicx}
\usepackage[subrefformat=parens]{subcaption}
\usepackage{amsopn}
\usepackage{color} 
\usepackage{amsthm} 
\newtheorem{thm}{Theorem}[section]
\newtheorem{pro}[thm]{Proposition}
\newtheorem{lem}[thm]{Lemma}
\newtheorem{rmk}[thm]{Remark}

\numberwithin{equation}{section}

\def\eps{\varepsilon}
\def\dis{\displaystyle}
\def\pl{\parallel}
\def\pd{\partial}

\begin{document}

\begin{center}
{\large Matched asymptotic expansion approach to pulse dynamics for a three-component reaction diffusion systems}

\vspace{3mm}

Yasumasa Nishiura\footnote{Research Institute for Electronic Science, Hokkaido University, Sapporo, 060-0812, Japan}
 and Hiromasa Suzuki\footnote{Faculty of Education, Shiga University, Hiratsu, Otsu, 520-0862, Japan}

\vspace{4mm}

\begin{abstract}
We study the existence and stability of standing pulse solutions to a singularly perturbed three-component reaction 
diffusion system with one-activator and two-inhibitor type. We apply the MAE (matched asymptotic expansion) method to 
the construction of solutions and the SLEP (Singular Limit Eigenvalue Problem) method to their stability properties. 
This approach is not just an alternative approach to geometric singular perturbation and the associated Evans function, 
but gives us two advantages: one is the extendability to higher dimensional case, and the other is to allow us to 
obtain more precise information on the behaviors of critical eigenvalues. This implies the existence of codimension 
two singularity of drift and Hopf bifurcations for the standing pulse solution and it is numerically confirmed that 
stable standing and traveling breathers emerge around the singularity in a physically-acceptable regime.
\end{abstract}


\end{center}

\section{Introduction}

Spatially localized patterns such as pulses and spots are fundamental components to understand the complex dynamics in dissipative systems. Since they are localized in space, they form a variety of patterns from a molecular shape to a crystal structure through the aggregation. Loosely speaking, each pattern behaves like an atom and a new structure emerges through the interaction among them. They interact each other either attractive or repulsive way depending on the parameters, and the interaction is either weak or strong depending on the distance: weak interaction means that they are well-separated and communicate through their tails; strong one means that they may lose their identities, typically as in collisions, and the final output is either annihilation, fusion or splitting (\cite{L}, \cite{N2}). 
Surprisingly, even though the transient dynamics looks very complicated, they eventually settle down to a coherent structure in most cases. At present it is difficult to give a rigor scenario to the whole process in strong interaction regime partly because large deformation is accompanied and the associated model system consists of three components, however an extensive numerics and semi-rigorous arguments reveal hidden structures that drive complex transient dynamics (\cite{N3}, \cite{NU2}, \cite{NU}). 

One of the lessons we learnt from those studies is to focus on a role of saddle solutions with high codimension. For instance, when two traveling spots collide strongly, they merge into one body or repel each other depending on the parameter values. There is a saddle solution of peanut shape that controls the dynamics as a separator. It is therefore quite important to search for a network of those saddles and study the connections among them (\cite{NTU2}, \cite{NTU}, \cite{NTU3}). 
Meanwhile it is known that there is a long history about the existence and stability of pulses for two-component reaction diffusion systems typically like the FitzHugh-Nagumo equations.
Despite its great success, the two-component systems and their variants are not adequate to study the above dynamics unfortunately. For instance, it is known that such a class does not support the coexistence of multiple stable traveling spots in 2D.
Mathematical model suitable for those dynamics is a class of three-component reaction diffusion systems. In fact, one of the roles of the third component $w$ is to prevent elongation or shrinkage of the spot orthogonal to the traveling direction so that it can sustain the localized shape and the coexistence of multiple stable traveling spots in higher dimensional space (\cite{VE}). Such a class gives a nice framework to consider collision dynamics and various ordered patterns consisting of many moving spots. 

One of the representative three-component systems was proposed as a qualitative model for the gas-discharged phenomenon with one activator and two inhibitors (Bode et al.\cite{BLSP}, Or-Guil et al.\cite{OBSP}, Purwins et al.\cite{PS}, Schenk et al.\cite{SOBP}). Another example is an activator-substrate-depleted model designed for shell patterns (Meinhardt \cite{ME}). Those models show a variety of dynamics including the annihilation, repulsion, fusion, and splitting even for 1D pulses upon collision (Nishiura, Teramoto and Ueda \cite{NTU2}, \cite{NTU}, \cite{NTU3}). Here we adopt the following $\eps$-scaled gas-discharged system for our model system:
  \begin{equation}
    \left \{
    \begin{array}{l}
       u_{t} = \eps^{2} u_{xx} + f(u) 
            - \eps (\alpha v + \beta w + \gamma),   \\
               \\
       \tau v_{t} = v_{xx} + u - v, \hspace{30mm} 
               (x, t) \in {\bf R} \times (0, \infty),  \\
               \\
       \theta w_{t} = D^2 w_{xx} + u - w, 
    \end{array}
    \right.
  \end{equation}
where $0<\eps <<1$, $D>0$, $\alpha>0$, $\beta>0$, $\gamma>0$, $\tau>0$, $\theta>0$ and $f(u)=u-u^{3}$. In this article, we maintain the assumptions $\alpha>0$, $\beta>0$ and $\gamma>0$ in which $u$ plays as an activator and $(v,w)$ as inhibitors. 
The above scaling was introduced by Doelman, van Heijster and Kaper \cite{DVK}, van Heijster, Doelman and Kaper \cite{vHDK} and fits for the singularly perturbed method. In fact, they succeeded to show the existence and the stability of pulse solutions to (1.1) by using the geometric singular perturbation method (GSP) and the Evans function (Alexander, Gardner and Jones \cite{AGJ}). One of the challenges for the above system is to study how the stability depends on the relaxation parameters $\tau$ and $\theta$. In particular, for the case $\tau=O(1/\eps^{2})$ and $\theta=O(1/\eps^{2})$, it was shown in \cite{vHDK} that the standing pulses lose their stabilities in various ways under several technical conditions. The highlight of this paper is to push forward their analysis and clarify the global behavior of critical eigenvalues in the above parameter regime, which allows us to find a higher codimension point in a wider parameter space. For this purpose, instead of GSP, we employ the MAE (Matched Asymptotic Expansion) approach and the SLEP (Singular Limit Eigenvalue Problem) method (\cite{NF}, \cite{NMIF}) for the control of critical eigenvalues. This is not just an alternative approach, in fact, we are able to remove several technical conditions assumed in \cite{vHDK} (see Remark 2.7) and give a nice perspective of all critical eigenvalues relevant to the stability of pulses. Our goal and the advantage are the following: firstly the SLEP allows us to control the precise parametric dependency of critical eigenvalues, in particular, for complex ones, namely we can show how they emerge from a degenerated real eigenvalue, cross the imaginal axis and come back to the real axis so that we can identify both the drift and Hopf bifurcations. This leads to finding a location of codimension two point more clearly in the full parameter space $(\tau, \theta)$ (see also Remark 2.7). Secondly MAE and SLEP method employed here remains valid for higher dimensional problems as was shown in \cite{NS1} and \cite{NS2}, which is very promising to accomplish the final goal in the future.

A couple of remarks are in order.
A relation between the Evans function and the SLEP equation was discussed in Ikeda, Nishiura and Suzuki \cite{INS} for the case of the front solutions. Basically these two methods are equivalent for 1D case. As for the front solutions, there appears a paper discussed from the view point of Bogdanov-Takens bifurcation and the reduction to center manifold by M.Chirilus-Bruckner, P.van Heijster, H.Ikeda and J.D.M.Rademacher \cite{BHIR}, which shows a mathematically nice structure of (1.1), although they violate our sign conditions of $\alpha>0$ and $\beta>0$ in some parameter regime (see also \cite{BDHR}). 
Finally, for the $\eps$-scaled system with heterogeneities, a new approach by using action functional was presented to prove the stability properties of the pulse solution (van Heijster, Chen, Nishiura and Teramoto \cite{vHCNT1}, \cite{vHCNT2}).

The article is organized as follows. In section 2, we mention the main results that give s perspective to the reader.  We outline the construction of stationary one-pulse solution of (1.1) (proof of Theorem 2.3) in section 3 by using the matched asymptotic expansion method. In section 4, the linearized eigenvalue problem (4.1) is reduced to the SLEP equation, and the stability properties of the one-pulse solution is studied for the case where $\tau$ and $\theta$ are of order $O(1)$ as $\eps \downarrow 0$ (that is, the proof of Theorem 2.5). In section 5, the stability properties for the case $\tau$ and $\theta$ are of order $O(1/\eps^{2})$ are presented (the proofs of Theorems 2.8 and 2.9, and Proposition 2.11). Finally, in section 6, we conclude the paper and discuss about the future problems. 


\section{Main results}

In this section, we summarize the main results and give a perspective of the article.
First of all, in view of the symmetry of the pulse solution, we can restrict the analysis on the interval $(0, \infty)$ without loss of generality. That is, we obtain the pulse solution by constructing the half of it on $(0,\infty)$ with Neumann boundary condition and flipping it. Concerning the stability analysis, it suffices to consider the eigenfunctions on the same interval with the Neumann (even) and Dirichlet (odd) boundary conditions at $x=0$.
The existence of stationary pulse solution has been already proved in \cite{DVK} by using the GSP theory, but we present another proof by the MAE method, which gives us necessary asymptotic forms for stability analysis later. The stationary problem on $(0,\infty)$ with Neumann boundary condition at $x=0$ reads
  \begin{equation}
    \left \{
    \begin{array}{l}
       0=\eps^{2} u_{xx} + f(u) - \eps (\alpha v + \beta w + \gamma),   \\
               \\
       0 = v_{xx} + u - v, \hspace{32mm} x \in (0, \infty),  \\
               \\
       0 = D^2 w_{xx} + u - w, 
    \end{array}
    \right.
  \end{equation}
  \begin{equation}
    \begin{array}{l}
       u_x(0) = 0, \ \ \ v_x (0) = 0, \ \ w_x (0)=0,  \\
                \\
       \dis \lim_{x \rightarrow \infty} (u(x),v(x),w(x)) = 
               (\overline{u}^{\eps}, \overline{v}^{\eps}, \overline{w}^{\eps}),
    \end{array}
  \end{equation}
where, $(\overline{u}^{\eps}, \overline{v}^{\eps}, \overline{w}^{\eps})$ is a negative constant solution of the following system:
  \begin{equation}
    \left \{
    \begin{array}{l}
       0 = f(u) - \eps (\alpha v + \beta w + \gamma),   \\
       0 = u - v,  \\
       0 = u - w.
    \end{array}
    \right.
  \end{equation}
The asymptotic form of it as $\eps \downarrow 0$ is given by
  \begin{displaymath}
    (\overline{u}^{\eps}, \overline{v}^{\eps}, \overline{w}^{\eps})
      = (-1, -1, -1) 
      + \eps \left( \dis\frac{\alpha+\beta-\gamma}{2}, 
          \dis\frac{\alpha+\beta-\gamma}{2}, \dis\frac{\alpha+\beta-\gamma}{2}
      \right) + o(\eps).
    \nonumber
  \end{displaymath}
It is clear that this is asymptotically stable in PDE sense. Let us denote the true solutions of (2.1) by $(u^{\eps}(x), v^{\eps}(x), w^{\eps}(x))$. In a singularly perturbed setting, there appears an internal layer at which $u$ becomes discontinuous as $\eps \downarrow 0$. The layer position $x=x^{*}$ is determined by solving the following reduced problem, which also gives the asymptotic forms of the inhibitors $(v^{\eps}(x), w^{\eps}(x))$ as $\eps \downarrow 0$. The second derivatives of the inhibitors do not degenerate so that they are expected to be $C^{1}$ even at the layer position. We construct smooth solutions of them except the layer position and match them at $x=x^{*}$. Recall that the reduced problem for the activator becomes an algebraic-like equation so that it can be represented as a function of two inhibitors (see section 3 for details).

  \begin{equation}
    \begin{array}{l}
    \left \{
    \begin{array}{l}
       0 = V_{xx}^{+} + 1 - V^{+}, \ \  x \in (0, x^{*}), \\
               \\
       V_{x}^{+}(0) = 0, \ \  V^{+}(x^{*}) = v^{*},
    \end{array}
    \right.  \\
                      \\
    \left \{
    \begin{array}{l}
       0 = V_{xx}^{-} - 1 - V^{-}, \ \  x \in (x^{*},\infty), \\
               \\
       V_{x}^{-}(0) = v^{*}, \ \  V^{+}(\infty) = -1,
    \end{array}
    \right.
    \end{array}
  \end{equation}
  \begin{equation}
    V_{x}^{-}(x^{*}) = V_{x}^{+}(x^{*})
  \end{equation}
  \begin{equation}
    \begin{array}{l}
    \left \{
    \begin{array}{l}
       0 = D^2 W_{xx}^{+} + 1 - W^{+}, \ \  x \in (0, x^{*}), \\
               \\
       W_{x}^{+}(0) = 0, \ \  W^{+}(x^{*}) = w^{*},
    \end{array}
    \right. \\
                      \\
    \left \{
    \begin{array}{l}
       0 = D^2 W_{xx}^{-} - 1 - W^{-}, \ \  x \in (x^{*},\infty), \\
               \\
       W_{x}^{-}(0) = w^{*}, \ \  W^{+}(\infty) = -1,
    \end{array}
    \right.
    \end{array}
  \end{equation}
  \begin{equation}
    W_{x}^{-}(x^{*}) = W_{x}^{+}(x^{*}),
  \end{equation}
where $x^{*}$, $v^{*}$, $w^{*}$ are determined by the following relations:
  \begin{equation}
    \alpha e^{-2 x^{*}} + \beta e^{-2 x^{*}/D} = \gamma, \ \ 
    v^{*} = -e^{-2 x^{*}}, \ \ w^{*} = -e^{-2 x^{*}/D}.
  \end{equation}
That is, $v^{*}$ and $w^{*}$ are the limiting values of $v^{\eps}(x)$ and $w^{\eps}(x)$ as $\eps \downarrow 0$, respectively. 

\begin{pro}[Reduced problem] 
Let $\alpha>0$, $\beta>0$, $\gamma>0$, $D>0$. If $\gamma < \alpha + \beta$, there exists a unique $x^{*}>0$ satisfying the first equation of $(2.8)$. The values at layer position $v^{*}$ and $w^{*}$ are defined by the second and third equations of $(2.8)$. Moreover there exist solutions $V^{+}(x)$ $($resp. $W^{+}(x)$$)$ $\in C^{2}(0,x^{*})$ and $V^{-}(x)$ $($resp. $W^{-}(x)$$)$ $\in C^{2}(x^{*},\infty)$ of $(2.4)$ $($resp. $(2.6)$$)$ satisfying $(2.5)$ $($resp. $(2.7)$$)$. 
\end{pro}

\begin{rmk}
The relations (2.8) are, what is called, the $C^{1}$-{\it matching conditions}. They are equivalent to (2.17) of Doelman et.al.\cite{DVK}, which determines the principal parts of the jump point between two slow manifolds ${\mathcal M}_{\eps}^{-}$ and ${\mathcal M}_{\eps}^{+}$, and the values of $v$ and $w$ at the point.
\end{rmk}

\vspace{2mm}

 In the spirit of MAE method,  we can smooth out the discontinuity by using the inner solution and obtain a family of smooth pulse solutions for $\eps>0$. 

\begin{thm}[Existence theorem of one-pulse solutions]
Under the assumptions of Proposition 2.1, there exists $\eps_{0}>0$ such that {\rm (2.1)}-{\rm (2.2)} has an $\eps$-family of one-pulse solutions $(u^{\eps}(x), v^{\eps}(x), w^{\eps}(x)) \in {\bf C}_{\rm unif}({\bf R})$ for $0 < \eps < \eps_{0}$. They satisfy
  \begin{displaymath}
    \dis\lim_{\eps \downarrow 0} u^{\eps}(x) = U^{*}(x) \hspace{5mm}
    \mbox{\it compact uniformly on } {\bf R} 
      \backslash I_{\delta}^{-} \cup I_{\delta}^{+} 
      \mbox{\it \ for any } \delta>0
  \end{displaymath}
and
  \begin{displaymath}
    \dis\lim_{\eps \downarrow 0} v^{\eps}(x) = V^{*}(x) \hspace{5mm}
    \mbox{\it compact uniformly on } {\bf R}, \hspace{21mm}
  \end{displaymath}
  \begin{displaymath}
    \dis\lim_{\eps \downarrow 0} w^{\eps}(x) = W^{*}(x) \hspace{5mm}
    \mbox{\it compact uniformly on } {\bf R}, \hspace{20mm}
  \end{displaymath}
where 
  \begin{displaymath}
    \begin{array}{l}
      {\bf C}_{\rm unif}({\bf R}) := C_{\rm unif}({\bf R}) \times 
                     C_{\rm unif}({\bf R}) \times C_{\rm unif}({\bf R}),  \\
             \\
      C_{\rm unif}({\bf R}) := \{ u \ | \ u \ 
          \mbox{is bounded and uniformly continuous on} \ {\bf R} \},
    \end{array}
  \end{displaymath}
$I_{\delta}^{\pm} = (\pm x^{*}-\delta, \pm x^{*}+\delta)$ and $x^{*}$ denotes the layer position of the reduced solution $(U^{*}, V^{*}, W^{*})$. Here the reduced solution $(U^{*}(x),V^{*}(x),W^{*}(x))$ is given by
  \begin{displaymath}
    U^{*}(x) = \left \{
    \begin{array}{rl}
        1, & \hspace{7mm} x \in (-x^{*}, x^{*})  \\
                   \\
       -1, & \hspace{7mm} x \in (-\infty, -x^{*}) \cup (x^{*}, \infty)
    \end{array}
    \right.
  \end{displaymath}
  \begin{displaymath}
    V^{*}(x) = \left \{
    \begin{array}{rl}
       V^{-}(-x), & x \in (-\infty, -x^{*})  \\
                   \\
       V^{+}(-x), & x \in [-x^{*}, 0)  \\
                   \\
       V^{+}(x), & x \in [0, x^{*})  \\
                   \\
       V^{-}(x), & x \in [x^{*}, \infty)  \\
    \end{array}
    \right. 
    W^{*}(x) = \left \{
    \begin{array}{rl}
       W^{-}(-x), & x \in (-\infty, -x^{*})  \\
                   \\
       W^{+}(-x), & x \in [-x^{*}, 0)  \\
                   \\
       W^{+}(x), & x \in [0, x^{*})  \\
                   \\
       W^{-}(x), & x \in [x^{*}, \infty)  \\
    \end{array}
    \right. 
  \end{displaymath}
\end{thm}

\begin{figure}[htbp]
  \begin{center}
     \includegraphics[width=7.2cm,height=4.2cm]{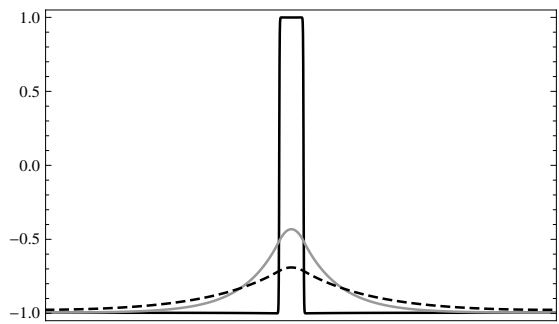} 
     \caption{A profile of stationary one-pulse solution: solid line for $u$, gray line for $v$, and dotted line for $w$. The parameters are given by $\eps=0.012$, $D=2.0$, $\alpha=1.0$, $\beta=2.0$, $\gamma=2.0$.}
    \end{center}
  \end{figure}
\noindent
Note that the pulse solutions constructed above are independent of the parameters $\tau$ and $\theta$. 

Stability properties of the pulse solutions are determined by the spectra of the associated linearized problem. As for the case $\tau = O(1)$ and $\theta = O(1)$, it is given by the following:
  \begin{equation}
    \left \{
    \begin{array}{l}
       \lambda^{\eps} p = \eps^{2} p_{xx} + f_{u}^{\eps} p 
            - \eps \alpha q - \eps \beta r,  \\
               \\
       \tau \lambda^{\eps} q = q_{xx} + p - q, 
                    \hspace{23mm} x \in {\bf R} \\
               \\
       \theta \lambda^{\eps} r = D^2 r_{xx} + p - r,  \\
               \\
       \dis\lim_{x \rightarrow \pm \infty} p(x) = 0, \ \ 
             \dis\lim_{x \rightarrow \pm \infty} q(x) = 0, \ \ 
             \dis\lim_{x \rightarrow \pm \infty} r(x) = 0,
    \end{array}
    \right.
  \end{equation}
where  $f_{u}^{\eps} := f'(u^{\eps}(x))$ and $u^{\eps}(x)$ denotes the $u$-component. We consider the spectral distribution of (2.9) in ${\bf L}_{2}({\bf R}) = L_{2}({\bf R}) \times L_{2}({\bf R}) \times L_{2}({\bf R})$ (for details, see Nishiura, Mimura, Ikeda and Fujii \cite{NMIF}). 

Similarly, the problem (2.9) can also be decomposed into an equivalent pair of eigenvalue problems on $(0, \infty)$, that is ${\rm (EP)}_{ev}$ and ${\rm (EP)}_{od}$ defined below.   \begin{equation}
    \left \{
    \begin{array}{rcl}
       \lambda p &=& \eps^{2} p_{xx} + f_{u}^{\eps} p 
            - \eps \alpha q - \eps \beta r  \\
               \\
       \tau \lambda q &=& q_{xx} + p - q, 
                    \hspace{20mm} x \in I := (0, \infty) \\
               \\
       \theta \lambda r &=& D^2 r_{xx} + p - r
    \end{array}
    \right.
  \end{equation}
with the boundary conditions
  \begin{equation}
    \lim_{x \rightarrow \infty} p(x) = 0, \ \ 
    \lim_{x \rightarrow \infty} q(x) = 0, \ \ 
    \lim_{x \rightarrow \infty} r(x) = 0.
  \end{equation}
Let ${\rm (EP)}_{ev}$ be the eigenvalue problem (2.10), (2.11) and 
  \begin{equation}
    p_{x}(0) = 0, \ \ q_{x}(0) = 0, \ \ r_{x}(0) = 0,
  \end{equation}
and ${\rm (EP)}_{od}$ be the eigenvalue problem (2.10), (2.11) and 
  \begin{equation}
    p(0) = 0, \ \ q(0) = 0, \ \ r(0) = 0.
  \end{equation}
Then, we have the following proposition. The proof is given in section 4 (see also Ikeda and Nishiura \cite{IN}).

\begin{pro} The stability properties of the stationary one-pulse solution are determined by the real parts of the spectra of ${\rm (EP)}_{ev}$ and ${\rm (EP)}_{od}$. 
\end{pro}

For later use, we define that an eigenvalue belongs to the class of {\it critical eigenvalues} when the real part of it tends to zero as $\eps \downarrow 0$. In view of the degeneracy of the second derivatives and noting that $u^{\eps}(x)$ becomes discontinuous as $\eps \downarrow 0$, the eigenvalues and associated eigenfunctions behave in a singular way, namely the eigenvalues may converge to zero (i.e., {\it critical eigenvalues}) and eigenfunctions are no more usual functions as $\eps \downarrow 0$. To overcome the difficulty, we resort to the {\it singular limit eigenvalue problem} (SLEP) method originated in Nishiura and Fujii \cite{NF}. The key idea of the SLEP method is to introduce suitable scalings and eliminate the singularity as $\eps \downarrow 0$. Using this method, finding critical eigenvalue of (2.10) is reduced to solving an algebraic-like equation called the SLEP equation. There is a huge amount of applications of this method including the higher dimensional problems (see for instance, Nishiura , Ikeda, Mimura and Fujii \cite{NMIF}, Nishiura and Mimura \cite{NM}, Nishiura and Suzuki \cite{NS1}, \cite{NS2}).

Now we are ready to state the stability property of the one-pulse solution of (1.1) constructed in Theorem 2.3 for the case $\tau = O(1)$ and $\theta = O(1)$. Note that there always exists a trivial zero eigenvalue coming from translation invariance, however this does not affect the stability properties.

\begin{thm}
Let $\tau$ and $\theta$ be $O(1)$ as $\eps \downarrow 0$. Then, the stability of one-pulse solution $(u^{\eps}(x),v^{\eps}(x),w^{\eps}(x))$ is determined by the following critical eigenvalue $\lambda^{\eps}$ of {\rm (2.9)}, which is real and simple, and has the principal part as $\eps \downarrow 0$ given by 
  \begin{displaymath}
    \lambda^{\eps} \approx - 3\sqrt{2} 
       \left( \alpha e^{-2 x^{*}} + \dis\frac{\beta}{D} e^{-2 x^{*}/D} 
       \right) \eps^{2}.
  \end{displaymath}
This implies from the assumptions of Theorem 2.3 that one-pulse solution $(u^{\eps}(x)$, $v^{\eps}(x)$, $w^{\eps}(x))$ is asymptotically stable when $\eps \downarrow 0$. The coefficient in front of $\eps^{2}$ is a solution of the SLEP equation {\rm (4.16)} in section 4. The essential spectrum of {\rm (2.9)} is uniformly bounded away from the imaginary axis for small $\eps$. Also all the other eigenvalues, if they exist, do not influence the stability of the one-pulse solution.
\end{thm} 

\begin{rmk}
The principal part of $\lambda^{\eps}$ is exactly the same as the leading order of (4.20) of \cite{vHDK} (see also (4.1) therein). Since we are interested in the physically-acceptable parameter regime, i.e., $\alpha$ and $\beta$ are positive, $\lambda^{\eps}$ is negative and hence the one-pulse solution is asymptotically stable. On the other hand, van Heijster et.al. \cite{vHDK} also study the case where $\alpha$ and $\beta$ are not necessarily positive.
\end{rmk}

\vspace{2mm}

The one-pulse solutions constructed above are independent of the parameters $\tau$ and $\theta$, however, their stability properties may depend on $\tau$ and $\theta$. In fact, a new regime was introduced such as $\tau = \hat{\tau}/\eps^{2}$ and $\theta = \hat{\theta}/\eps^{2}$ by \cite{vHDK}, and various instabilities were indicated there. We also move on to this case. The linearized eigenvalue problem for $(\hat{\tau}, \hat{\theta})$ is recast as
  \begin{equation}
    \left \{
    \begin{array}{rcl}
       \lambda p &=& \eps^{2} p_{xx} + f_{u}^{\eps} p 
            - \eps \alpha q - \eps \beta r,  \\
               \\
       \dis\frac{1}{\eps^2} \hat{\tau} \lambda q &=& q_{xx} + p - q,
                    \hspace{23mm} x \in {\bf R} \\
               \\
       \dis\frac{1}{\eps^2} \hat{\theta} \lambda r &=& D^2 r_{xx} + p - r,
    \end{array}
    \right.
  \end{equation}
with boundary conditions
  \begin{equation}
    \lim_{x \rightarrow \pm \infty} p(x) = 0, \ \ 
    \lim_{x \rightarrow \pm \infty} q(x) = 0, \ \ 
    \lim_{x \rightarrow \pm \infty} r(x) = 0.
  \end{equation}
In this case, the response of the inhibitors $v$ and $w$ are very slow so that it causes the instability of the one-pulse solution in various ways. In fact, we will see that the stability is determined by the three critical eigenvalues: one is real and the other is a pair of complex conjugate ones. In this regime, the essential spectrum dos not affect the stability property of the pulse (see Proposition 5.1) so that we can concentrate on the behavior of these three critical eigenvalues.
%

\begin{rmk}
The stability analysis of \cite{vHDK} for $\tau = O(1/\eps^{2})$ and $\theta = O(1/\eps^{2})$ is restricted to the case $\theta=1$, namely $\hat{\theta}=0$, and the bifurcation analysis is also considered there. As for the existence of codimension two (drift $+$ Hopf) point is discussed under the assumption $D \rightarrow \infty$. We don't need these assumptions and study the problem under a general setting, which allows us to search for the codimension two points in a wider parameter space.
\end{rmk}

For later convenience, we rewrite the linearized eigenvalue problem (2.14) as 
  \begin{displaymath}
    {\mathcal L} \varPhi = \lambda \widehat{T} \varPhi, 
  \end{displaymath}
where
  \begin{displaymath}
    {\mathcal L} = 
    \left(
      \begin{array}{ccc}
         L^{\eps} & -\eps \alpha & -\eps \beta  \\
                \\
         1 & M & 0  \\
                \\
         1 &  0 & N
      \end{array}
    \right), \ \ 
    \widehat{T} = 
    \left(
      \begin{array}{ccc}
         1 & 0 & 0  \\
                \\
         0 & \hat{\tau}/\eps^{2} & 0  \\
                \\
         0 & 0 & \hat{\theta}/\eps^{2}
      \end{array}
    \right), \ \ 
  \end{displaymath}
  \begin{displaymath}
     L^{\eps} = \eps^{2} \dis\frac{d^{2}}{dx^{2}} + f_{u}^{\eps}, \ \ 
     M := \dis\frac{d^{2}}{dx^{2}} - 1 , \ \ 
     N:= D^{2} \dis\frac{d^{2}}{dx^{2}} - 1.
  \end{displaymath}
In the same manner as before, the problem can be decomposed into an equivalent pair of problems ${\rm (EP)}_{od}$ and ${\rm (EP)}_{ev}$.
First let us consider the case ${\rm (EP)}_{od}$ satisfying (2.13).

\begin{thm}
Let $\tau$ and $\theta$ be $O(1/\eps^{2})$ $(\eps \downarrow 0)$. Then the critical eigenvalue $\lambda^{\eps}$ of {\rm (2.14)} satisfying the odd symmetry {\rm (2.13)} is real and simple. The principal part of $\lambda^{\eps}$ as $\eps \downarrow 0$ is given by
  \begin{displaymath}
    \lambda^{\eps} \approx 
          \hat{\lambda}_{d}^{*}(\hat{\tau}, \hat{\theta}) \ \eps^{2},
  \end{displaymath}
where $\hat{\lambda}_{d}^{*}(\hat{\tau}, \hat{\theta})$ is a solution of the {\rm SLEP} equation $G_{od}(\hat{\lambda}^{*}; \hat{\tau}, \hat{\theta}) = 0$ {\rm (}see {\rm (5.8))} and satisfies

\noindent
{\rm (i)}
  \begin{displaymath}
    \hat{\lambda}_{d}^{*}(\hat{\tau},\hat{\theta}) \left \{
    \begin{array}{l}
      < 0 \ \ \ \mbox{\it for} \ \ (\hat{\tau},\hat{\theta}) \in \Omega_{-}  \\
                     \\
      = 0 \ \ \ \mbox{\it for} \ \ (\hat{\tau},\hat{\theta}) \in \Gamma_{d}  \\
                     \\
      > 0 \ \ \ \mbox{\it for} \ \ (\hat{\tau},\hat{\theta}) \in \Omega_{+},
    \end{array}
    \right. \hspace{30mm}
  \end{displaymath}
{\rm (ii)}
  \begin{equation}
    \dis\frac{\pd}{\pd \hat{\tau}} 
          \hat{\lambda}_{d}^{*} (\hat{\tau}_{0}, \hat{\theta}_{0}) > 0, \ \ \ 
    \dis\frac{\pd}{\pd \hat{\theta}} 
          \hat{\lambda}_{d}^{*} (\hat{\tau}_{0}, \hat{\theta}_{0}) > 0 \ \ \ \ 
    \mbox{for} \ (\hat{\tau}_{0},\hat{\theta}_{0}) \in \Gamma_{d},
    \nonumber
  \end{equation}
where
  \begin{displaymath}
    \begin{array}{l}
      \Omega_{\pm} = \left \{
        (\hat{\tau},\hat{\theta}) \in {\bf R}_{+} \times {\bf R}_{+} ; 
      \dis\frac{\pd}{\pd \hat{\lambda}^{*}} 
         G_{od}(0; \hat{\tau}, \hat{\theta}) \lessgtr 0,
         \hat{\tau}>0, \hat{\theta}>0 \right \},  \\
                  \\
      \Gamma_{d} = \left \{
        (\hat{\tau},\hat{\theta}) \in {\bf R}_{+} \times {\bf R}_{+} ; 
      \dis\frac{\pd}{\pd \hat{\lambda}^{*}} 
          G_{od}(0; \hat{\tau}, \hat{\theta}) = 0,
          \hat{\tau}>0, \hat{\theta}>0 \right \}
    \end{array}
  \end{displaymath}
and
  \begin{displaymath}
    \dis\frac{\pd}{\pd \hat{\lambda}^{*}} G_{od}(0; \hat{\tau}, \hat{\theta})
         = 1 - C_{1} \, \hat{\tau} - C_{2} \, \hat{\theta}.
  \end{displaymath}
Here ${\bf R}_{+} = \{x \in {\bf R} \ | \ x>0  \}$, $C_{1}$ and $C_{2}$ are positive constants {\rm (}for more details, see section {\rm 5.1)}. That is, $\Gamma_{d}$ is a drift bifurcation set and given by the straight line in $\hat{\tau}$-$\hat{\theta}$ plane as depicted in Fig.3. The one-pulse solution can be destabilized when the parameters $(\hat{\tau},\hat{\theta})$ cross $\Gamma_{d}$ from $\Omega_{-}$ to $\Omega_{+}$ transversally. In fact, the algebraic multiplicity of the eigenvalue $\hat{\lambda}_{d}^{*}(\hat{\tau}_{0},\hat{\theta}_{0})=0$ $((\hat{\tau}_{0},\hat{\theta}_{0}) \in \Gamma_{d})$ to the operator $\widehat{T}_{0}^{-1} {\mathcal L}$ is two, where $\widehat{T}_{0}$ is $\widehat{T}$ at $(\hat{\tau},\hat{\theta})=(\hat{\tau}_{0},\hat{\theta}_{0})$. Finally, all the other real eigenvalues, if they exist, do not influence the stability of the one-pulse solution.
\end{thm}

On the other hand, for the symmetric case ${\rm (EP)}_{ev}$, we can derive the SLEP equation $G_{ev}(\hat{\lambda}^{*}; \hat{\tau}, \hat{\theta})=0$ (see (5.9)) and it is found that there occurs an instability of Hopf type when the parameters $(\hat{\tau},\hat{\theta})$ cross the curve $\Gamma_{H}$ defined below in $\hat{\tau}$-$\hat{\theta}$ plane. In order to present our results, we introduce the polar coordinate $(s, \psi)$ in $(\hat{\tau}, \hat{\theta})$, and set
  \begin{displaymath}
    \hat{\tau} = s \cos \psi, \ \ \ \hat{\theta} = s \sin \psi, \ \ \ 
    s>0, \ \ \ 0 < \psi < \dis\frac{\pi}{2}.
  \end{displaymath}
Then the SLEP equation for ${\rm (EP)}_{ev}$ takes the following form
  \begin{equation}
    \widehat{G}_{ev}(\hat{\lambda}^{*}; s, \psi)
        := G_{ev}(\hat{\lambda}^{*}; s \cos \psi, s \sin \psi) = 0
  \end{equation}
(for details, see section 5.2). Then we have

\begin{thm}
 Let $\tau$ and $\theta$ be $O(1/\eps^{2})$ as $\eps \downarrow 0$. Then, the critical eigenvalue $\lambda^{\eps}$ satisfying the even symmetry {\rm (2.12)} consists of a unique pair of complex numbers of $O(\eps^{2})$. More precisely, for any $\psi \in (0, \pi/2)$, there exist a unique $s^{*}(\psi)>0$, $\xi^{*}(\psi)>0$ and a constant $\delta_{0}>0$ such that a unique isolated pair of complex solutions $\hat{\lambda}_{h}^{*}(s,\psi)$ and $\overline{\hat{\lambda}_{h}^{*}(s,\psi)}$ of {\rm (2.16)} for $s \in (s^{*}(\psi)-\delta_{0}, s^{*}(\psi)+\delta_{0})$ becomes the principal part of critical eigenvalue of $O(\eps^{2})$. They cross the imaginary axis transversally from left to right at $s=s^{*}(\psi)$ when $s$ increases. Namely, we have
  \begin{displaymath}
    \lambda^{\eps} \approx 
          \hat{\lambda}_{h}^{*}(s,\psi) \ \eps^{2}, \ \ 
          \overline{\hat{\lambda}_{h}^{*}(s,\psi)} \eps^{2}
  \end{displaymath}
for $s \in (s^{*}(\psi)-\delta_{0}, s^{*}(\psi)+\delta_{0})$, 
  \begin{displaymath}
    \hat{\lambda}_{h}^{*}(s^{*}(\psi),\psi) = \xi(\psi) i, \ \ \ 
    \dis\frac{\pd}{\pd s} 
        {\mathrm Re} \ \hat{\lambda}_{h}^{*}(s^{*}(\psi), \psi) > 0.
  \end{displaymath}
Here $s^{*}(\psi)$ and $\xi^{*}(\psi)$ are defined by
  \begin{displaymath}
    \begin{array}{l}
    \hat{\zeta}_{0}^{*} = 4 (\kappa^{*})^{2} \left [
      \alpha x^{*} R \left ( 
      \dis\frac{1}{2} \tan^{-1} (s^{*}(\psi) \xi^{*}(\psi) \cos \psi), 2 x^{*}
      \right) \right.  \\
                   \\ \hspace{25mm}
      \left. + \dis\frac{\beta x^{*}}{D^{2}} R 
      \left  ( 
      \dis\frac{1}{2} \tan^{-1} (s^{*}(\psi) \xi^{*}(\psi) \sin \psi), 
     \dis\frac{2 x^{*}}{D} 
      \right )
      \right ],
    \end{array}
  \end{displaymath}
  \begin{displaymath}
    \begin{array}{l}
      \xi^{*}(\psi) = - 4 (\kappa^{*})^{2} \left [
      \alpha x^{*} I 
        \left ( 
        \dis\frac{1}{2} \tan^{-1} (s^{*}(\psi) \xi^{*}(\psi) \cos \psi),2 x^{*}
        \right) \right.  \\
                  \\ \hspace{25mm}
      \left. + \dis\frac{\beta x^{*}}{D^{2}} I 
        \left  ( 
           \dis\frac{1}{2} \tan^{-1} (s^{*}(\psi) \xi^{*}(\psi) \sin \psi), 
           \dis\frac{2 x^{*}}{D}
       \right )
      \right ],
    \end{array}
  \end{displaymath}
where
  \begin{displaymath}
    \renewcommand{\arraystretch}{1.7}
    \begin{array}{l}
      R(z;d):= \dis\frac{1}{d} 
             \sqrt{\cos 2z} \ [\cos z + e^{-x(z;d)} \cos(y(z;d) + z) ], \\
                \\
      I(z;d):= - \dis\frac{1}{d} 
               \sqrt{\cos 2z} \ [\sin z + e^{-x(z;d)} \sin(y(z;d) + z) ], \\
                \\
      x(z;d):=\dis\frac{d \cos z}{\sqrt{\cos 2z}}, \ \ 
      y(z;d):=\dis\frac{d \sin z}{\sqrt{\cos 2z}}.
    \end{array}
  \end{displaymath}
Then Hopf bifurcation curve $\Gamma_{H}$ is given by
  \begin{displaymath}
    \Gamma_{H} = 
    \{ (\hat{\tau}, \hat{\theta}) \ | \ \hat{\tau} = s^{*}(\psi) \cos \psi, 
     \ \hat{\theta} = s^{*}(\psi) \sin \psi, \ 0 < \psi < \pi/2 \}.
  \end{displaymath}
\end{thm} 

\begin{rmk}
Note that we have no restrictions on $\theta$ and $D$ in the above theorem, which is not the case for \cite{vHDK}.
\end{rmk}

Theorem 2.9 claims that a pair of complex eigenvalues crosses the imaginary axis, and Hopf bifurcation occurs there. A natural question is how these eigenvalues behave before and after the Hopf bifurcation. In fact we are able to trace the behavior globally and find the transitions from real eigenvalues to a complex pair and vice versa as in the next proposition. This partly stems from the fact that for any solution $(\hat{\lambda}_{0}^{*}, \hat{\tau}_{0}, \hat{\theta}_{0}) \in {\bf C} \times {\bf R}_{+} \times {\bf R}_{+}$ of $G_{ev}(\hat{\lambda}^{*}; \hat{\tau}, \hat{\theta}) =0$, we have
  \begin{displaymath}
    \dis\frac{d}{d \hat{\lambda}^{*}} 
      G_{ev}(\hat{\lambda}_{0}^{*}; \hat{\tau}_{0}, \hat{\theta}_{0}) \neq 0
  \end{displaymath}
if ${\rm Im} \ \hat{\lambda}_{0}^{*} \neq 0$ (for more details, see Proposition 5.8 in section 5). 

\begin{pro}
For any $\psi \in (0, \pi/2)$, there exist positive constants $\underline{s}(\psi)$ and $\overline{s}(\psi)$ with $\underline{s}(\psi) < \overline{s}(\psi)$ such that the critical complex eigenvalues behave as follows in ${\bf C}$ (see Fig.2):

\vspace{2mm}

{\rm (i)} There exist exactly two real negative eigenvalues of {\rm (2.16)} for $0 < s < \underline{s}(\psi)$.

\vspace{2mm}

{\rm (ii)} There exists a unique negative eigenvalue $\hat{\lambda}^{*,-}$ at 
$s = \underline{s}(\psi)$  with double multiplicity. Near $s = \underline{s}(\psi)$, $\hat{\lambda}^{*,-}$ splits into two eigenvalues in the following way:
  \begin{displaymath}
    \begin{array}{l}
      \hat{\lambda}^{*} \approx \hat{\lambda}^{*,-} 
            \pm \sqrt{c_{-} (s - \underline{s}(\psi))} \ \ \ 
            \mbox{\it for} \ \ s < \underline{s}(\psi),   \\
             \\
      \hat{\lambda}^{*} \approx \hat{\lambda}^{*,-} 
            \pm i \sqrt{c_{-} (\underline{s}(\psi) - s)} \ \ \ 
            \mbox{\it for} \ \ s > \underline{s}(\psi),
    \end{array}
  \end{displaymath}
where $c_{-}$ is a negative constant defined by
  \begin{displaymath}
    c_{-} = - \dis\frac{2 \widehat{G}_{e s}(\hat{\lambda}^{*,-}; \underline{s}(\psi), \psi)}
         {\widehat{G}_{e \lambda \lambda}(\hat{\lambda}^{*,-}; \underline{s}(\psi), \psi)} < 0
  \end{displaymath}
and the subscripts $s$ and $\lambda \lambda$ mean the partial derivatives.

\vspace{2mm}

{\rm (iii)} There are no real eigenvalues for $\underline{s}(\psi) < s < \overline{s}(\psi)$.

\vspace{2mm}

{\rm (iv)} There exists a unique positive eigenvalue $\hat{\lambda}^{*,+}$ at 
$s = \overline{s}(\psi)$  with double multiplicity. Near $s = \overline{s}(\psi)$, $\hat{\lambda}^{*,+}$ splits into two eigenvalues in the following way:
  \begin{displaymath}
    \begin{array}{l}
      \hat{\lambda}^{*} \approx \hat{\lambda}^{*,+} 
            \pm \sqrt{c_{+} (s - \overline{s}(\psi))} \ \ \ 
            \mbox{\it for} \ \ s > \overline{s}(\psi),   \\
             \\
      \hat{\lambda}^{*} \approx \hat{\lambda}^{*,+} 
            \pm i \sqrt{c_{+} (\overline{s}(\psi) - s)} \ \ \ 
            \mbox{\it for} \ \ s < \overline{s}(\psi),
    \end{array}
  \end{displaymath}
where $c_{+}$ is a positive constant defined by
  \begin{displaymath}
    c_{+} = - \dis\frac{2 \widehat{G}_{e s}(\hat{\lambda}^{*,+}; \overline{s}(\psi), \psi)}
         {\widehat{G}_{e \lambda \lambda}(\hat{\lambda}^{*,+}; \overline{s}(\psi), \psi)} > 0.
  \end{displaymath}

{\rm (v)} There exist exactly two real positive eigenvalues for 
$\overline{s}(\psi) < s$.
\end{pro}

Fig.2 shows the global behavior of eigenvalues as $s>0$ varies. Although the dotted part remains to be proved rigorously, it shows how a pair of complex eigenvalues behave qualitatively.

  \begin{figure}[htbp]
    \begin{center}
       \includegraphics[width=7cm]{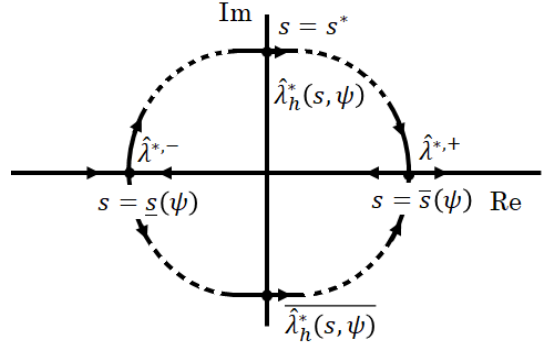}
       \caption{Global behavior of complex eigenvalues: the two real eigenvalues merge at $s=\underline{s}(\psi)$, cross the imaginary axis at $s=s^{*}$, and fall into the real axis at $s= \overline{s}(\psi)$.}
    \end{center}
  \end{figure}
The two degenerated points $\underline{s}(\psi)$ and $\overline{s}(\psi)$ and the Hopf bifurcation in between are not yet obtained explicitly, however it is possible to locate them numerically. 

\begin{pro}
Under the assumptions of Theorem 2.3, it is confirmed numerically that there are two codimension two bifurcation points consisting of drift and Hopf bifurcations in $(\hat{\tau}, \hat{\theta})$-space when $\alpha$, $\beta$, $\gamma$, and $D$ are appropriately chosen (see two large black dots in Fig.3). There appear three different pulse dynamics: traveling pulse, standing and traveling breathers in addition to the standing pulse by unfolding the codimension two points. In fact, Fig.4 (a)-(d) shows those dynamics around the bottom singularity.
\end{pro}

In fact, the drift bifurcation line (the solid line in Fig.3) is explicitly written by the set $\Gamma_{d}$ of Theorem 2.8, and Hopf bifurcation line (dotted curve in Fig.3) can be obtained by solving the SLEP equation (5.9) by the Newton method. It is confirmed numerically that the two bifurcation curves cross transversally at exactly two points. Since these lines are smooth functions of parameters $\alpha$, $\beta$, $\gamma$ and $D$, and there are no restrictions for the parameters $\hat{\tau}$, $\hat{\theta}$ and $D$ unlike \cite{vHDK}, it is expected that codimension two bifurcation points exist for generic parameters' range.

Fig.4 displays four different pulse dynamics around the lower codimension two point in Fig.3, in which the time evolutions of contour line $u(t,x)=0$ are presented. They are computed by using implicit scheme with $\Delta x = 7.0 \times 2^{-10}$ and $\Delta t = 0.012$ and the system size is $14.0$ subject to Neumann boundary conditions. 
In this paper we do not go into the details about the unfolding the codimension two points , which is left for the future work.
%
\begin{figure}[htbp]
  \begin{center}
    \includegraphics[width=5.8cm, height=5.8cm]{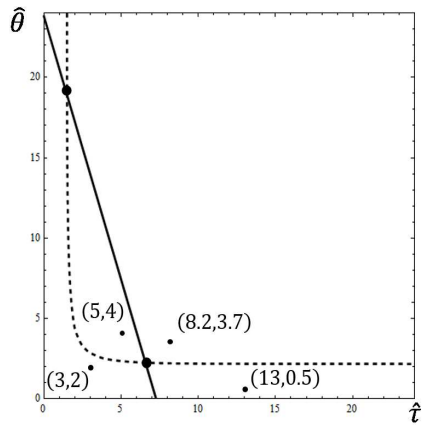}
    \caption{Bifurcation diagram in $(\hat{\tau}$, $\hat{\theta})$-space. The solid and dotted lines indicate the set of drift and Hopf bifurcations, respectively. Other parameters are set as $\alpha=1.0$, $\beta=2.0$, $\gamma=2.0$, $D=2.0$, $x^{*}=0.311905$.}
  \end{center}
\end{figure}

\begin{figure}[htbp]
 \begin{tabular}{c}
  \begin{minipage}{0.4\textwidth}
    \centering
    \includegraphics[width=5.2cm]{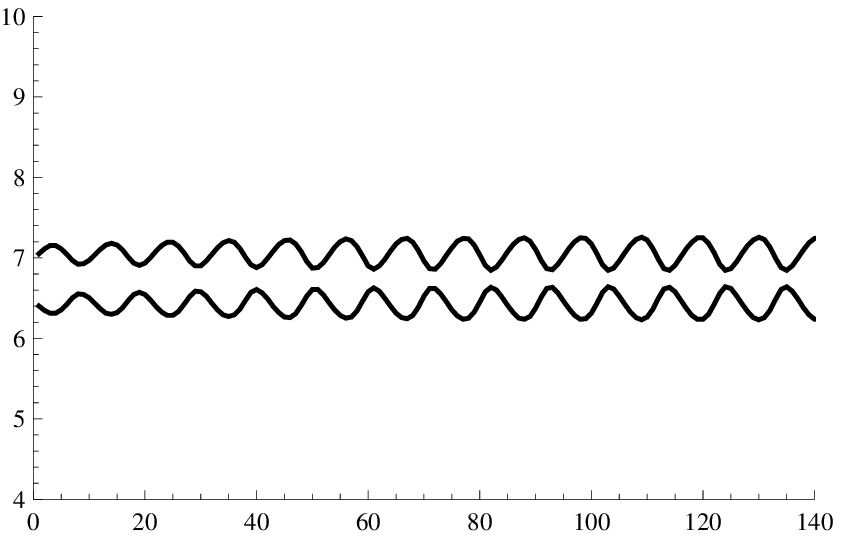}
    \subcaption{$(\hat{\tau}, \hat{\theta}) = (5,4)$}
    \includegraphics[width=5.2cm]{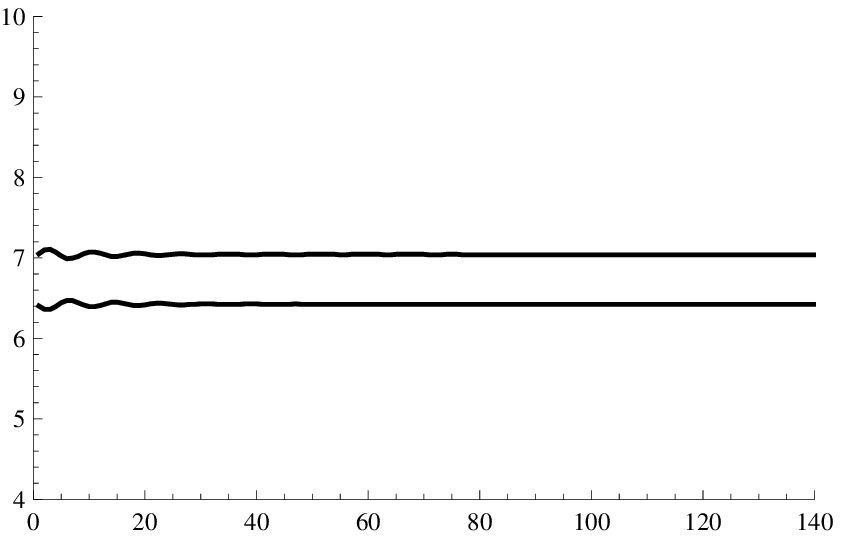}
    \subcaption{$(\hat{\tau}, \hat{\theta}) = (3,2)$}
  \end{minipage}

\hspace{10mm}

  \begin{minipage}{0.4\textwidth}
    \centering
    \includegraphics[width=5.2cm]{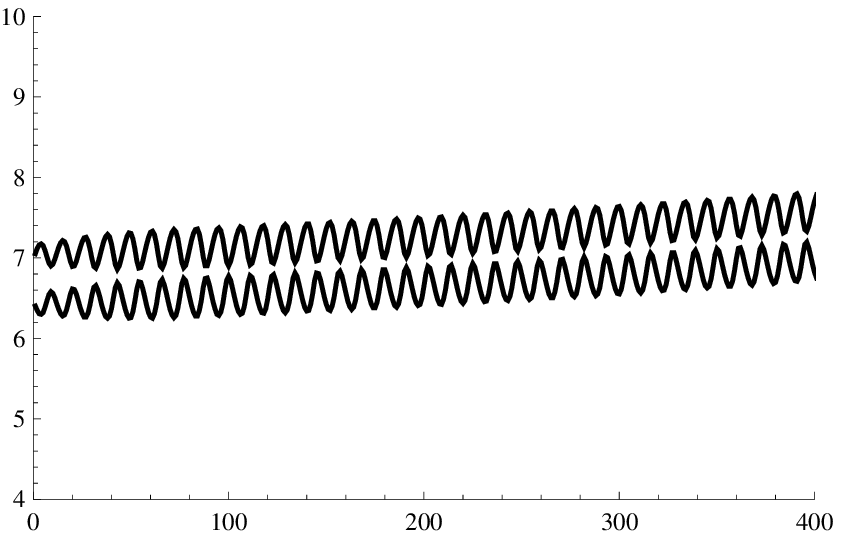}
    \subcaption{$(\hat{\tau}, \hat{\theta}) = (8.2,3.7)$}
    \includegraphics[width=5.2cm]{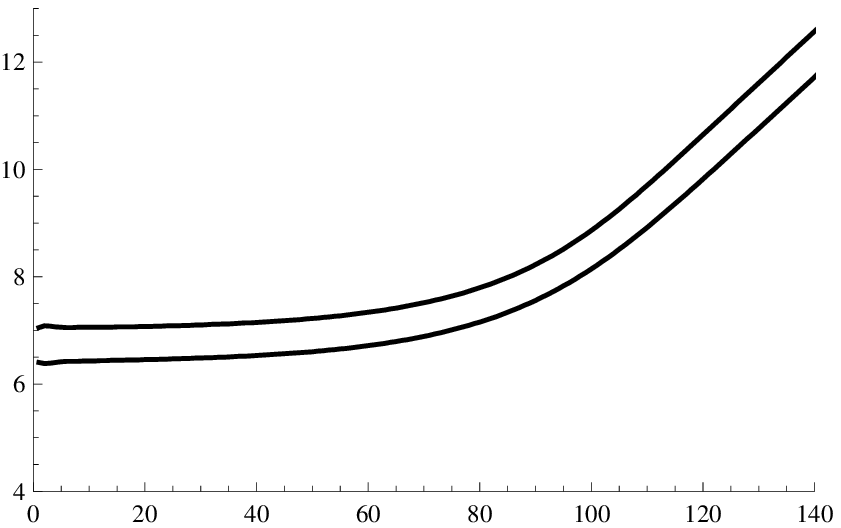}
    \subcaption{$(\hat{\tau}, \hat{\theta}) = (13,0.5)$}
  \end{minipage}
 \end{tabular}
 \caption{Pulse dynamics around the lower codimension-two point in Fig.3. The contour line $u(t,x)=0$ for each case is shown as time evolves (horizontal axis). Here we used an implicit scheme with $\Delta x = 7.0 \times 2^{-10}$ and $\Delta t = 0.012$ and the system size is $14.0$ subject to Neumann boundary conditions. }
\end{figure}

\section{Construction of standing pulse solution}

In the following, we only outline the proofs of Proposition 2.1 and Theorem 2.3. For the details, see \cite{N}, \cite{NF} and \cite{IMN}. As was mentioned in section 2, it suffices to construct the pulse on the half interval $(0, \infty)$ with the Neumann boundary condition at $x=0$.

If we set $\eps=0$ in (2.1) formally, the first equation is no more a differential equation so that $u$-component becomes discontinuous at some point $x=x^{*}$ at which two solutions satisfying the boundary conditions at $x=0$ and $x=\infty$ respectively should be matched there. First we divide the problem (2.1) into an equivalent pair of problems on the interval $(0, x^{*})$ and $(x^{*}, \infty)$ with $x=x^{*}$ being a free parameter. We determine the layer position $x=x^{*}$ later by using $C^{1}$-matching conditions.
  \begin{equation}
    \left \{
    \begin{array}{l}
       0=\eps^{2} u_{xx}^{+} + f(u^{+}) 
           - \eps (\alpha v^{+} + \beta w^{+} + \gamma),   \\
               \\
       0 = v_{xx}^{+} + u^{+} - v^{+}, \hspace{32mm} x \in (0,x^{*}),  \\
               \\
       0 = D^2 w_{xx}^{+} + u^{+} - w^{+}, \\
               \\
       u_{x}^{+}(0) = 0, \ \ \ v_{x}^{+} (0) = 0, \ \ w_{x}^{+} (0)=0,  \\
               \\
       u^{+}(x^{*}) = a_{0}, \ \ v^{+}(x^{*}) = b_{0} + \eps b_{1}, \ \ 
            w^{+}(x^{*}) = c_{0} + \eps c_{1}
    \end{array}
    \right.
  \end{equation}
  \begin{equation}
    \left \{
    \begin{array}{l}
       0=\eps^{2} u_{xx}^{-} + f(u^{-}) 
            - \eps (\alpha v^{-} + \beta w^{-} + \gamma),   \\
               \\
       0 = v_{xx} + u^{-} - v^{-}, \hspace{32mm} x \in (x^{*}, \infty),  \\
               \\
       0 = D^2 w_{xx}^{-} + u^{-} - w^{-}, \\
               \\
       u^{-}(x^{*}) = a_{0}, \ \ v^{-}(x^{*}) = b_{0} + \eps b_{1}, \ \ 
            w^{-}(x^{*}) = c_{0} + \eps c_{1}  \\
               \\
       \dis \lim_{x \rightarrow \infty} (u^{-}(x),v^{-}(x),w^{-}(x)) = 
            (\overline{u}^{\eps}, \overline{v}^{\eps}, \overline{w}^{\eps}),
    \end{array}
    \right.
  \end{equation}
where $-1<a_{0}<1$, $b_{j}$ and $c_{j}$ ($j=0,1$) are the parameters determined by the $C^{1}$-matching conditions later.

\vspace{2mm}

\noindent
{\bf 3.1. Outer expansion}

\vspace{2mm}

Let us substitute
  \begin{displaymath}
    u = U^{0,\pm} + \eps U^{1,\pm} + \eps^{2} U^{2,\pm} + o(\eps^{2}), \
    v = V^{0,\pm} + \eps V^{1,\pm} + \eps^{2} V^{2,\pm} + o(\eps^{2}),
  \end{displaymath}
  \begin{displaymath}
    w = W^{0,\pm} + \eps W^{1,\pm} + \eps^{2} W^{2,\pm} + o(\eps^{2})
  \end{displaymath}
into (3.1) and (3.2), and equating like powers of $\eps^{k}$, we have the following problems for $(U^{k,\pm}$, $V^{k,\pm}$, $W^{k,\pm})$ ($k = 0, 1, 2$):
  \begin{displaymath}
    \left \{
    \begin{array}{l}
       0 = f(U^{0,\pm}),   \\
               \\
       0 = V_{xx}^{0,\pm} + U^{0,\pm} - V^{0,\pm},   \\
               \\
       0 = D^{2} W_{xx}^{0,\pm} + U^{0,\pm} - W^{0,\pm}
    \end{array}
    \right.
  \end{displaymath}
  \begin{equation}
    \left \{
    \begin{array}{l}
       0 = f'(U^{0,\pm}) U^{1,\pm} 
           - (\alpha V^{0,\pm} + \beta W^{0,\pm} + \gamma), \\
               \\
       0 = V_{xx}^{1,\pm} + U^{1,\pm} - V^{1,\pm},  \\
               \\
       0 = D^{2} W_{xx}^{1,\pm} + U^{1,\pm} - W^{1,\pm}.
    \end{array}
    \right.
  \end{equation}
  \begin{equation}
    \left \{
    \begin{array}{l}
       0 = f'(U^{0,\pm}) \ U^{2,\pm}(x) 
           + F^{2,\pm}(U^{0,\pm},V^{0,\pm},W^{0,\pm}, 
                       U^{1,\pm},V^{1,\pm},W^{1,\pm}), \\
               \\
       0 = V_{xx}^{2,\pm} - V^{2,\pm} 
             + G^{2,\pm} (U^{0,\pm},V^{0,\pm},W^{0,\pm}, 
                          U^{1,\pm},V^{1,\pm},W^{1,\pm}), \\
               \\
       0 = D^2 W_{xx}^{2,\pm} - W^{2,\pm} 
             + H^{2,\pm} (U^{0,\pm},V^{0,\pm},W^{0,\pm}, 
                          U^{1,\pm},V^{1,\pm},W^{1,\pm}),
    \end{array}
    \right.
  \end{equation}
  \begin{displaymath}
    V_{x}^{2,+}(0) = 0, \ \ \ W_{x}^{2,+}(0) = 0, \ \ \ 
    V^{2,\pm}(x^{*}) = 0, \ \ \ W^{2,\pm}(x^{*}) = 0
  \end{displaymath}
  \begin{displaymath}
    \begin{array}{l}
      \dis \lim_{x \rightarrow \infty} V^{2,-}(x)
      = \dis \lim_{x \rightarrow \infty} 
        G^{2,-} (U^{0},V^{0},W^{0},U^{1},V^{1},W^{1})   \\
               \\
      \dis \lim_{x \rightarrow \infty} W^{2,-}(x)
      = \dis \lim_{x \rightarrow \infty} 
        H^{2,-} (U^{0},V^{0},W^{0},U^{1},V^{1},W^{1})
    \end{array}
  \end{displaymath}
where $f'(u) = 1 - 3 u^{2}$, $F^{2,\pm}$, $G^{2,\pm}$ and $H^{2,\pm}$ are functions with respect to $x$ depending on $(U^{i,\pm}, V^{i,\pm}$, $W^{i,\pm})$ ($i = 0, 1$).

In view of bistability, we choose $U^{0,\pm}(x)$ as $U^{0,\pm}(x) \equiv \pm 1$. Then the equations for $V^{0,\pm}$ and $W^{0,\pm}$ are recast as
  \begin{equation}
    \left \{
    \begin{array}{l}
       0 = V_{xx}^{0, \pm} \pm 1 - V^{0, \pm},   \\
               \\
       V_{x}^{0,+}(0) = 0, \ V^{0,\pm}(x^{*}) = b_{0}, \\
               \\
       \dis \lim_{x \rightarrow \infty} V^{0,-}(x) = -1,
    \end{array}
    \right. \hspace{5mm}
    \left \{
    \begin{array}{l}
       0 = D^{2} W_{xx}^{0, \pm} \pm 1 - W^{0, \pm},  \\
               \\
       W_{x}^{0,+}(0) = 0, \ W^{0,\pm}(x^{*}) = c_{0}, \\
               \\
       \dis \lim_{x \rightarrow \infty} W^{0,-}(x) = -1
    \end{array}
    \right.
  \end{equation}
Since $f'(\pm 1)=-2<0$, we can solve the first equation of (3.3) with respect to $U^{1,\pm}$ as
  \begin{equation}
    U^{1,\pm}(x) 
     = - \dis\frac{1}{2} (\alpha V^{0,\pm}(x) + \beta W^{0,\pm}(x) + \gamma).
  \end{equation}
Substituting (3.6) into the second and the third equation of (3.3), we have the equations for $V^{1}$ and $W^{1}$:
  \begin{equation}
    \left \{
    \begin{array}{l}
       0 = V_{xx}^{1,\pm} - V^{1,\pm} - \dis\frac{1}{2} 
             (\alpha V^{0,\pm} + \beta W^{0,\pm} + \gamma),  \\
               \\
       V_{x}^{1,+}(0) = 0, \ V^{1,\pm}(x^{*}) = b_{1}, \\
               \\
       \dis \lim_{x \rightarrow \infty} V^{1,-}(x) 
       =\dis \lim_{x \rightarrow \infty} \dis\frac{1}{f'(U^{0,-})}
         (\alpha V^{0,-}(x) + \beta W^{0,-}(x) + \gamma)  \\
               \\ \hspace{21mm}
       = - \dis\frac{1}{2} (\alpha + \beta + \gamma)
    \end{array}
    \right.
  \end{equation}

\vspace{-5mm}

  \begin{equation}
    \left \{
    \begin{array}{l}
       0 = D^2 W_{xx}^{1,\pm} - W^{1,\pm} - \dis\frac{1}{2} 
             (\alpha V^{0,\pm} + \beta W^{0,\pm} + \gamma),  \\
               \\
       W_{x}^{1,+}(0) = 0, \ W^{1,\pm}(x^{*}) = c_{1}, \\
               \\
      \dis \lim_{x \rightarrow \infty} W^{1,-}(x) 
      = \dis \lim_{x \rightarrow \infty} \dis\frac{1}{f'(U^{0,-})} 
        (\alpha V^{0,-}(x) + \beta W^{0,-}(x) + \gamma)  \\
               \\ \hspace{21mm}
      = - \dis\frac{1}{2} (\alpha + \beta + \gamma)
    \end{array}
    \right.
  \end{equation}
For the outer solutions for the inhibitors $(V,W)$, we can solve the above equations successively. In fact, noting that the fundamental solutions of homogeneous equation $D W_{xx} - W = 0$ are $e^{-x/D}$ and $e^{x/D}$, we are able to solve (3.8) with respect to $W^{1,\pm}$. Similarly for $V^{1,\pm}$. Once $(V^{i,\pm}, W^{i,\pm})$ (i=0,1) are determined, we can solve (3.4) with respect to $V^{2,\pm}$ and $W^{2,\pm}$. On the other hand, $U^{0,+}$ and $U^{0,-}$ are not continuous at $x=x^{*}$, we need the inner expansion in order to correct the approximation in the neighborhood of $x=x^{*}$.

\vspace{2mm}

\noindent
{\bf 3.2. Inner expansion}

\vspace{2mm}

We introduce a stretched variable $y = (x-x^{*})/\eps$ and make inner expansion in the neighborhood of $x=x^{*}$ to smooth out the gap of $U^{+,\eps}$ and $U^{-,\eps}$ at $x=x^{*}$.
In the following, we only consider the problem on $(0,x^{*})$. 

We determine the functions $u^{i,\pm}$, $v^{i,\pm}$ and $w^{i,\pm}$ ($i=0,1,2$) in the following expressions:
  \begin{equation}
    \left \{
    \begin{array}{l}
      u = U^{0,\pm}(x) + \eps U^{1,\pm}(x) + \eps^{2} U^{2,\pm}(x) \\
                   \\ \hspace{7mm}
      + u^{0,\pm}((x-x^{*})/\eps) 
      + \eps u^{1,\pm}((x-x^{*})/\eps) + \eps^{2} u^{2,\pm}((x-x^{*})/\eps), \\
                   \\
      v = V^{0,\pm}(x) + \eps V^{1,\pm}(x) + \eps^{2} V^{2,\pm}(x) \\
                   \\ \hspace{7mm}
        + \eps^{2} [v^{0,\pm}((x-x^{*})/\eps) + \eps v^{1,\pm}((x-x^{*})/\eps)
        + \eps^{2} v^{2,\pm}((x-x^{*})/\eps)  ], \\
                   \\
      w = W^{0,\pm}(x) + \eps W^{1,\pm}(x) + \eps^{2} W^{2,\pm}(x) \\
                   \\ \hspace{7mm}
        + \eps^2 [w^{0,\pm}((x-x^{*})/\eps) + \eps w^{1,\pm}((x-x^{*})/\eps)
        + \eps^{2} w^{2,\pm}((x-x^{*})/\eps) ],
    \end{array}
    \right.
  \end{equation}
Substituting (3.9) into (2.1) and equating like powers of $\eps$, we have 
  \begin{equation}
    \left \{
    \begin{array}{l}
      u_{y y}^{0,\pm} + f(U^{0,\pm}(x^{*})+u^{0,\pm}) = 0,  \\
            \\
      u^{0,\pm}(\mp \infty) = 0, \ \ u^{0,\pm}(0) = a_{0} - U^{0,\pm}(x^{*}),
    \end{array}
    \right.
  \end{equation}
  \begin{equation}
    \left \{
    \begin{array}{l}
      v_{y y}^{0,\pm} + u^{0,\pm} = 0,  \\
          \\
      v^{0,\pm}(\mp \infty) = 0, \ \ v^{0,\pm}(0) = 0,
    \end{array}
    \right. \ \ 
    \left \{
    \begin{array}{l}
      D^{2} w_{y y}^{0,\pm} + u^{0,\pm} = 0,  \\
          \\
      w^{0,\pm}(\mp \infty) = 0, \ \ w^{0,\pm}(0) = 0,
    \end{array}
    \right.
  \end{equation}
  \begin{equation}
    \left \{
    \begin{array}{l}
      u_{y y}^{1,\pm} + \tilde{f}_{u}^{\pm} u^{1,\pm} 
         = P_{1}^{\pm}(y), \\
            \\
      u^{1,\pm}(\mp \infty) = 0, \ \ u^{1,\pm}(0) = - U^{1,\pm}(x^{*}),
    \end{array}
    \right.
  \end{equation}
where
  \begin{equation}
    P_{1}^{\pm}(y) = - \tilde{f}_{u}^{\pm}  [U_{x}^{0,\pm}(x^*) y + U^{1,\pm}(x^*)]
        + \alpha V^{0,\pm}(x^{*}) + \beta W^{0,\pm}(x^{*}) + \gamma.
        \nonumber
  \end{equation}
  \begin{equation}
    \left \{
    \begin{array}{l}
      v_{y y}^{1,\pm} + u^{1,\pm} = 0,  \\
          \\
      v^{1,\pm}(\mp \infty) = 0, \ \ v^{1,\pm}(0) = 0,
    \end{array}
    \right.
    \left \{
    \begin{array}{l}
      D^{2} w_{y y}^{1,\pm} + u^{1,\pm} = 0,  \\
          \\
      w^{1,\pm}(\mp \infty) = 0, \ \ w^{1,\pm}(0) = 0,
    \end{array}
    \right.
  \end{equation}
  \begin{equation}
    \left \{
    \begin{array}{l}
      u_{y y}^{2,\pm} + \tilde{f}_{u}^{\pm} u^{2,\pm} = P_{2}^{\pm}(y), \\ 
                    \\
      u^{2,\pm}(\mp \infty) = 0, \ \ u^{2,\pm}(0) = - U^{2,\pm}(x^{*}),
    \end{array}
    \right.
  \end{equation}
where
  \begin{displaymath}
    \begin{array}{l}
      P_{2}^{\pm}(y) 
       = - \dis\frac{1}{2} \tilde{f}_{uu}^{\pm} [U^{1,\pm}(x^{*}) 
       + u^{1,\pm}]^2 - \tilde{f}_{u}^{\pm} 
          [U_{x}^{1,\pm}(x^*) y + U^{2,\pm}(x^*)] \\
                       \\ \hspace{15mm}
        + \alpha (V_{x}^{0,\pm}(x^{*}) y + V^{1,\pm}(x^{*}))
        + \beta  (W_{x}^{0,\pm}(x^{*}) y + W^{1,\pm}(x^{*})),
    \end{array}
  \end{displaymath}
  \begin{equation}
    \left \{
    \begin{array}{l}
      v_{y y}^{2,\pm} + u^{2,\pm} - v^{0,\pm} = 0,  \\
          \\
      v^{2,\pm}(\mp \infty) = 0, \ \ v^{2,\pm}(0) = 0,
    \end{array}
    \right.
    \left \{
    \begin{array}{l}
      D^{2} w_{y y}^{2,\pm} + u^{2,\pm} - w^{0,\pm} = 0,  \\
          \\
      w^{2,\pm}(\mp \infty) = 0, \ \ w^{2,\pm}(0) = 0,
    \end{array}
    \right.
  \end{equation}
where
  \begin{displaymath}
    \tilde{f}_{u}^{\pm} := f'(U^{0,\pm}(x^{*})+u^{0,\pm}(y)), \ \ \ 
    \tilde{f}_{uu}^{\pm} := f''(U^{0,\pm}(x^{*})+u^{0,\pm}(y))
    \nonumber
  \end{displaymath}

The next lemma is basic for the solvability of (3.10). See, for instance, \cite{N}.

\begin{lem} 
There exists a unique solution $\tilde{u}^{*}(y)$ of
  \begin{displaymath}
    \left \{
    \begin{array}{l}
      u_{yy} + f(u) = 0, \hspace{10mm} y \in {\bf R} \\
             \\
      \dis\lim_{y \rightarrow -\infty} u(y) = U^{0,+}(x^{*}) = 1, \ \ \ 
      \dis\lim_{y \rightarrow \infty} u(y) = U^{0,-}(x^{*}) = -1, \ \ \ 
      u(0) = 0.
    \end{array}
    \right.
  \end{displaymath}
In fact, $\tilde{u}^{*}(y)$ is explicitly represented as
  \begin{equation}
    \tilde{u}^{*}(y) = - \tanh \dis\frac{y}{\sqrt{2}}
  \end{equation}
\end{lem}

By using $\tilde{u}^{*}(y)$, $u^{0,+}(y)$ is given by
  \begin{equation}
    u^{0,\pm}(y) = \tilde{u}^{*}(y+s) - U^{0,\pm}(x^{*})
  \end{equation}
where, $s$ is a real number satisfying $\tilde{u}^{*}(s)=a_{0}$, which is uniquely determined since $-1<a_{0}<1$ and the monotonicity of $\tilde{u}^{*}(y)$.
By using the expression of $u^{0,+}(y)$, $v^{0,+}(y)$ and $w^{0,+}(y)$ are determined by (3.11) respectively.

Taking account of the fact that $u_{y}^{0,\pm}(y) = \tilde{u}_{y}^{*}(y+s)$ is a monotone solution of the homogeneous equation
  \begin{displaymath}
    \phi_{y y} + \tilde{f}' \phi = 0,
  \end{displaymath}
the solution $u^{1,\pm}(y)$ can be expressed as
  \begin{equation}
    \begin{array}{l}
    u^{1,\pm}(y) = - \dis\frac{U^{1,\pm}(x^{*})}
      {\tilde{u}_{y}^{*}(s)} \tilde{u}_{y}^{*}(y+s)   \\
                      \\ \hspace{20mm}
      + \tilde{u}_{y}^{*}(y+s) \dis\int_{0}^{y} [\tilde{u}_{y}^{*}(z+s)]^{-2}
      \left( \int_{\mp \infty}^{z} P_{1}^{\pm}(\tau) \tilde{u}_{y}^{*}(\tau+s) 
          d\tau \right) dz.
    \end{array}
  \end{equation}
Once $u^{0,\pm}(y)$ and $u^{1,\pm}(y)$ are determined, $u^{2,\pm}(y)$ is given by
  \begin{equation}
    \begin{array}{l}
    u^{2,\pm}(y) = - \dis\frac{U^{2,\pm}(x^{*})}
      {\tilde{u}_{y}^{*}(s)} \tilde{u}_{y}^{*}(y+s)   \\
                     \\  \hspace{20mm}
      + \tilde{u}_{y}^{*}(y+s) \dis\int_{0}^{y} [\tilde{u}_{y}^{*}(z+s)]^{-2}
      \left( \int_{\mp \infty}^{z} P_{2}^{\pm}(\tau) \tilde{u}_{y}^{*}(\tau+s) 
          d\tau \right) dz.
      \end{array}
  \end{equation}
Then $v^{i,\pm}(y)$ and $w^{i,\pm}(y)$ ($i=1,2$) are determined by (3.13) and (3.15) respectively.

\vspace{2mm} 

\noindent
{\bf 3.3. Justification and $C^{1}$-matching conditions}

\vspace{2mm} 

Now we are ready to construct the singularly perturbed solutions of (3.1) and (3.2) with the following asymptotic forms:
  \begin{equation}
    \left \{
    \begin{array}{l}
      u^{\eps,\pm}(x) = U^{0,\pm}(x) + \eps U^{1,\pm}(x)
        + \eps^{2} U^{2,\pm}(x)   \\
                 \\ \hspace{2mm}
        + [ u^{0,\pm}((x-x^{*})/\eps) + \eps u^{1,\pm}((x-x^{*})/\eps)
        + \eps^{2} u^{2,\pm}((x-x^{*})/\eps) ] \omega(x),  \\
                   \\
      v^{\eps,\pm}(x) = V^{0,\pm}(x) + \eps V^{1,\pm}(x) 
        + \eps^{2} V^{2,\pm}(x) \\
                 \\ \hspace{2mm}
        + \eps^2 [v^{0,\pm}((x-x^{*})/\eps) + \eps v^{1,\pm}((x-x^{*})/\eps)
        + \eps^{2} v^{2,\pm}((x-x^{*})/\eps)] \omega(x), \\
                   \\
      w^{\eps,\pm}(x) = W^{0,\pm}(x) + \eps W^{1,\pm}(x) 
        + \eps^{2} W^{2,\pm}(x)    \\
                 \\ \hspace{2mm}
        + \eps^2 [w^{0,\pm}((x-x^{*})/\eps) + \eps w^{1,\pm}((x-x^{*})/\eps)
        + \eps^{2} w^{2,\pm}((x-x^{*})/\eps)] \omega(x),
    \end{array}
    \right.
  \end{equation}
where $\omega(x)$ is a smooth cutoff function such that
  \begin{displaymath}
    \omega(x) = 1, \ |x-x^{*}| \leq \dis\frac{x^{*}}{4}, \ \ 
        \omega(x) = 1, \ |x-x^{*}| \geq \dis\frac{x^{*}}{2}
  \end{displaymath}

In order to obtain a smooth pulse solution, the derivatives of 
$(u^{\eps,-}(x), $ $v^{\eps,-}(x),$ $w^{\eps,-}(x))$ and $(u^{\eps,+}(x), v^{\eps,+}(x), w^{\eps,+}(x))$ should be matched at $x=x^{*}$. The difference of the derivatives at $x=x^{*}$ can be computed as
  \begin{equation}
    \left \{
    \begin{array}{l}
      \eps u_{x}^{\eps,-}(x^{*}) - \eps u_{x}^{\eps,+}(x^{*})
        = u_{y}^{0,-}(0) - u_{y}^{0,+}(0) 
        + \eps [u_{y}^{1,-}(0) - u_{y}^{1,+}(0)]   \\
                   \\ \hspace{10mm}
        + \eps^{2} [ \{U_{x}^{1,-}(x^{*}) + u_{y}^{2,-}(0) \} - 
           \{U_{x}^{1,+}(x^{*}) + u_{y}^{2,+}(0) \}] + o(\eps^{2})  \\
                   \\
      v_{x}^{\eps,-}(x^{*}) - v_{x}^{\eps,+}(x^{*})
      = V_{x}^{0,-}(x^{*}) - V_{x}^{0,+}(x^{*})  \\
                   \\ \hspace{10mm}
      + \eps [ \{V_{x}^{1,-}(x^{*}) + v_{y}^{0,-}(0) \} - 
           \{V_{x}^{1,+}(x^{*}) + v_{y}^{0,+}(0) \}] + o(\eps) \\
                   \\
      w_{x}^{\eps,-}(x^{*}) - w_{x}^{\eps,+}(x^{*})
      = W_{x}^{0,-}(x^{*}) - W_{x}^{0,+}(x^{*})  \\
                   \\ \hspace{10mm}
      + \eps [ \{W_{x}^{1,-}(x^{*}) + w_{y}^{0,-}(0) \} - 
           \{W_{x}^{1,+}(x^{*}) + w_{y}^{0,+}(0) \}] + o(\eps)
    \end{array}
    \right.
  \end{equation}
Noting the expression (3.17) of $u^{0,\pm}(y)$, we see that $u_{y}^{0,-}(0) - u_{y}^{0,+}(0)=0$ is already satisfied. 
On the other hand,
  \begin{equation}
    V_{x}^{0,-}(x^{*}) - V_{x}^{0,+}(x^{*}) = 0 \ \  \mbox{and} \ \ 
    W_{x}^{0,-}(x^{*}) - W_{x}^{0,+}(x^{*}) = 0
  \end{equation}
are the conditions for the reduced problems (see (2.5) and (2.7)). They are satisfied if the solutions of the reduced problems exist.

For the inner expansion, we see from (3.18)
  \begin{equation}
    u_{y}^{1,\pm}(0) = \dis\frac{a_{1} - U^{1,\pm}(x^{*})}
      {\tilde{u}_{y}^{*}(s)} \tilde{u}_{yy}^{*}(s)
      + \dis \frac{1}{\tilde{u}_{y}^{*}(s)} 
      \dis\int_{-\infty}^{0} P_{1}^{\pm}(\tau) \tilde{u}_{y}^{*}(\tau+s) d\tau,
  \end{equation}
where
  \begin{displaymath}
    P_{1}^{\pm}(y) = - \tilde{f}' U^{1,\pm}(x^*)
        + \alpha b_{0} + \beta c_{0} + \gamma.
  \end{displaymath}
The integral part of (3.23) can be computed as
  \begin{displaymath}
    \dis\int_{-\infty}^{0} P_{1}^{\pm}(\tau) \tilde{u}_{y}^{*}(\tau+s) d\tau
      = U^{1,\pm}(x^*) \tilde{u}_{yy}^{*}(s)
        + (\alpha b_{0} + \beta c_{0} + \gamma) (\tilde{u}^{*}(s) - 1),
  \end{displaymath}
so we have
  \begin{displaymath}
    \tilde{u}_{y}^{*}(s) u_{y}^{1,\pm}(0) = a_{1} \tilde{u}_{yy}^{*}(s)
     + (\alpha b_{0} + \beta c_{0} + \gamma) (\tilde{u}^{*}(s) \mp 1).
  \end{displaymath}
This yields the fact that $u_{y}^{1,-}(0) - u_{y}^{1,+}(0)=0$ is equivalent to
  \begin{equation}
    \alpha b_{0} + \beta c_{0} + \gamma = 0
  \end{equation}
since $\tilde{u}_{y}^{*}(s) \neq 0$, which gives a relation between $b_{0}$ and $c_{0}$.

On the other hand, we have the following expression for (3.22).

\begin{lem}
The conditions of $(3.22)$ are equivalent to
  \begin{equation}
    b_{0} = -e^{-2x^{*}} \ \ \mbox{and} \ \ c_{0} = -e^{-2x^{*}/D},
  \end{equation}
respectively.
\end{lem}

\begin{proof}
Noting that the fundamental solutions of homogeneous equation 
$D^{2} W_{xx} - W = 0$ are $e^{-x/D}$ and $e^{x/D}$, we obtain
  \begin{displaymath}
    W^{0,+}(x) = \dis\frac{(c_{0}-1) \cosh (x/D)}{\cosh (x^{*}/D)} + 1, \ \ \ 
    W^{0,-}(x) = (c_{0}+1) e^{(x^{*} - x)/D} - 1.
  \end{displaymath}
Then their derivatives are given by
  \begin{displaymath}
    W_{x}^{0,+}(x) = \dis\frac{(c_{0}-1) \sinh (x/D)}{D \cosh (x^{*}/D)}, \ \ \ 
    W_{x}^{0,-}(x) = - \dis\frac{1}{D}(c_{0}+1) e^{(x^{*} - x)/D},
  \end{displaymath}
so we have
  \begin{displaymath}
    W_{x}^{0,+}(x^{*}) 
         = \dis\frac{(c_{0}-1) \sinh (x^{*}/D)}{D \cosh (x^{*}/D)}, \ \ \ 
    W_{x}^{0,-}(x^{*}) = - \dis\frac{1}{D} (c_{0}+1).
  \end{displaymath}
Solving the equation $W_{x}^{0,+}(x^{*}) = W_{x}^{0,-}(x^{*})$ with respect to $c_{0}$, we have
  \begin{displaymath}
    c_{0} = -e^{-2 x^{*}/D}.
  \end{displaymath}
Similarly,
  \begin{displaymath}
    b_{0} = -e^{-2 x^{*}}
  \end{displaymath}
(set $D=1$ and replace $W$ by $V$).
\end{proof}

Substituting (3.25) into (3.24), we have
  \begin{equation}
     \alpha e^{-2 x^{*}} + \beta e^{-2 x^{*}/D} = \gamma.
  \end{equation}

\begin{lem}
Let $\alpha>0$, $\beta>0$, $\gamma>0$, $D>0$. If $\gamma < \alpha + \beta$, there exists a unique $x^{*}>0$ satisfying $(3.26)$.
\end{lem}
\begin{proof}
Let $g(z):=\alpha e^{-z} + \beta e^{-z/D}$. It is easy to verify that
  \begin{displaymath}
    g'(z)<0, \ \ g(0)=\alpha + \beta, \ \ \lim_{z \rightarrow \infty} g(z) = 0
  \end{displaymath}
by direct computation. The intermediate value theorem guarantees the conclusion.
\end{proof}

\begin{rmk} 
(i) (3.26) is the principal part of the $C^{1}$-matching condition. It is equivalent to (3.17) of \cite{DVK}, which determines the principal parts of the jump point between two slow manifolds ${\mathcal M}_{\eps}^{-}$ and ${\mathcal M}_{\eps}^{+}$.

(ii) We concentrate on the case where $\alpha>0$, $\beta>0$ and $\gamma>0$. If these are not satisfied, there exist more than one solution $x^{*}$ satisfying (3.26), which yields the existence of the two-pulse solution as was discussed in section 5 of \cite{DVK}.
\end{rmk}

\vspace{2mm}

Once $x^{*}$ is determined, $b_{0}$ and $c_{0}$ can be obtained by (3.25). 
Note that $a_{0}$ and $s$ satisfy $\tilde{u}^{*}(s)=a_{0}$ and $\tilde{u}^{*}(s)$ is a monotonically decreasing function.

Now we have the asymptotic forms of the stationary solution $(v^{\eps}(x), w^{\eps}(x))$ as $\eps \downarrow 0$, which is the claim of Proposition 2.1. 
The principal term $a_{0}$ of the boundary value at $x=x^{*}$ of $u^{\eps}$ remains unknown yet. 
In order to obtain a true one-pulse solution, we have to match also the higher order terms of (3.21).
In the following, we show that $(s, b_{1}, c_{1})$ can be determined by $C^{1}$-matching condition for the next order $O(\eps)$-term:
  \begin{equation}
    \left \{
    \begin{array}{l}
       \{U_{x}^{1,-}(x^{*}) + u_{y}^{2,-}(0) \} - 
           \{U_{x}^{1,+}(x^{*}) + u_{y}^{2,+}(0) \} = 0,  \\
             \\
       \{V_{x}^{1,-}(x^{*}) + v_{y}^{0,-}(0) \} - 
           \{V_{x}^{1,+}(x^{*}) + v_{y}^{0,+}(0) \} = 0,  \\
             \\
       \{W_{x}^{1,-}(x^{*}) + w_{y}^{0,-}(0) \} - 
           \{W_{x}^{1,+}(x^{*}) + w_{y}^{0,+}(0) \} = 0.
    \end{array}
    \right.
  \end{equation}
Before computing $u_{y}^{2,\pm}(0)$, we see from (3.6), (3.26) and Lemma 3.3 that
  \begin{displaymath}
      U^{1,\pm}(x^{*}) = \dis\frac{1}{2} (\alpha e^{-2 x^{*}} + \beta e^{-2 x^{*}/D} - \gamma) = 0.
  \end{displaymath}
Moreover combining this and the fact that $U^{0,\pm}(x)$ are constant functions, we have
$u^{1,\pm}(y) \equiv 0$ since the inhomogeneous terms $P_{1}^{\pm}(y)$ of the equation (3.12) 
for $u^{1,\pm}(y)$ is equal to zero. In view of (3.19), we have
  \begin{equation}
    u_{y}^{2,\pm}(0) = - \dis\frac{U^{2,\pm}(x^{*})}
      {\tilde{u}_{y}^{*}(s)} \tilde{u}_{yy}^{*}(s)
      + \dis \frac{1}{\tilde{u}_{y}^{*}(s)} 
      \dis\int_{\mp \infty}^{0} \{ P_{21}^{\pm}(\tau) + P_{22}^{\pm}(\tau) \} \tilde{u}_{y}^{*}(\tau+s)
          d\tau
  \end{equation}
where
  \begin{displaymath}
    \begin{array}{l}
    P_{21}^{\pm}(y) = - \tilde{f}' \{ U_{x}^{1,\pm}(x^*) y + U^{2,\pm}(x^*) \}  \\
                \\
    P_{22}^{\pm}(y) =     \alpha \{ V_{x}^{0,\pm}(x^{*}) y + V^{1,\pm}(x^{*}) \}
        + \beta \{ W_{x}^{0,\pm}(x^{*}) y + W^{1,\pm}(x^{*}) \}.
    \end{array}
  \end{displaymath}
Next we compute each term of the integral part of (3.28):
  \begin{displaymath}
    \begin{array}{l}
      - \dis\int_{\mp \infty}^{0} P_{21}^{\pm}(y) \ \tilde{u}_{y}^{*}(y+s) dy
         = U_{x}^{1,\pm}(x^*) \Big[f(\tilde{u}^{*}(y+s)) y \Big]_{\mp \infty}^{0}  \\
                   \\ \hspace{15mm}
         - U_{x}^{1,\pm}(x^*) \dis\int_{\mp \infty}^{0} f(\tilde{u}^{*}(y+s))dy
                   + U^{2,\pm}(x^*) \dis\int_{\mp \infty}^{0} 
               \tilde{f}' \tilde{u}_{y}^{*}(y+s) dy  \\
                    \\ \hspace{5mm}
        = U_{x}^{1,\pm}(x^*) \tilde{u}_{y}^{*}(s) - 
           U^{2,\pm}(x^*) \tilde{u}_{yy}^{*}(s)
    \end{array}
  \end{displaymath}
  \begin{displaymath}
    \begin{array}{l}
      \dis\int_{\mp \infty}^{0} P_{22}^{\pm}(y) \ \tilde{u}_{y}^{*}(y+s) dy \\
                 \\
       = [\alpha V_{x}^{0,\pm}(x^{*}) + \beta W_{x}^{0,\pm}(x^{*})] 
       \dis\int_{\mp \infty}^{0} y \tilde{u}_{y}^{*}(y+s) dy
       + (\alpha b_{1} + \beta c_{1})
         \dis\int_{\mp \infty}^{0} \tilde{u}_{y}^{*}(y+s) dy \\
                 \\
       = [\alpha V_{x}^{0,\pm}(x^{*}) + \beta W_{x}^{0,\pm}(x^{*})] 
       \dis\int_{\mp \infty}^{s} z \tilde{u}_{y}^{*}(z) dz  \\
                 \\ \hspace{10mm}
       + [\alpha b_{1} + \beta c_{1} 
          - \{ \alpha V_{x}^{0,\pm}(x^{*}) + \beta W_{x}^{0,\pm}(x^{*}) \} s]
          \dis\int_{\mp \infty}^{s} \tilde{u}_{y}^{*}(z) dz \\
                 \\
       = [\alpha V_{x}^{0,\pm}(x^{*}) + \beta W_{x}^{0,\pm}(x^{*})] 
       \dis\int_{\mp \infty}^{s} z \tilde{u}_{y}^{*}(z) dz  \\
                 \\ \hspace{10mm}
       + \Big[\alpha b_{1} + \beta c_{1} 
     - \{ \alpha V_{x}^{0,\pm}(x^{*}) + \beta W_{x}^{0,\pm}(x^{*}) \} s \Big]
          (\tilde{u}^{*}(s) \mp 1)
    \end{array}
  \end{displaymath}
Using these results, we have
  \begin{equation}
    \begin{array}{l}
      \tilde{u}_{y}^{*}(s) [\{ u_{y}^{2,-}(0) + U_{x}^{1,-}(x^{*}) \}
         - \{ u_{y}^{2,+}(0) + U_{x}^{1,+}(x^{*}) \} ]  \\
                     \\ \hspace{7mm}
      = 2 \Big[\alpha b_{1} + \beta c_{1} 
        - \{ \alpha V_{x}^{0}(x^{*}) + \beta W_{x}^{0}(x^{*}) \} s \Big] 
    \end{array}
  \end{equation}
Here we used the fact that
  \[
    \dis\int_{-\infty}^{\infty} z \tilde{u}_{y}^{*}(z) dz = 0
  \]
thanks to the oddness of  $y \tilde{u}_{y}^{*}(y)$. Also we omitted the superscript $\pm$ of $V^{0}(x)$ and $W^{0}(x)$ because they are already matched at $x=x^{*}$ in $C^{1}$-sense. 

\begin{lem} 
The second and the third equations of $(3.27)$ give the following relations:
  \begin{equation}
    b_{1} = - (1 + e^{-2x^{*}}) s + C_{1}, \ \ \ 
    c_{1} = - \dis\frac{1}{D} (1 + e^{-2x^{*}/D}) s + C_{2},
  \end{equation}
where $C_{1}$ and $C_{2}$ are constants independent of $s$
  \end{lem}
\begin{proof}
By using the variation of constants formula, we can compute $W^{1,+}(x)$ and $W^{1,-}(x)$ as
  \begin{displaymath}
    W^{1,+}(x) = \dis\frac{(c_{1}-1) \cosh (x/D)}{\cosh (x^{*}/D)} + 1, \ \ \ 
    W^{1,-}(x) = (c_{1}+1) e^{(x^{*} - x)/D} - 1.
  \end{displaymath}
Here we used the fact
  \begin{displaymath}
    \alpha b_{0} + \beta c_{0} + \gamma = 0.
  \end{displaymath}
Therefore we have
  \begin{displaymath}
    \begin{array}{l}
      W_{x}^{1,-}(x^{*}) - W_{x}^{1,+}(x^{*}) = 
       - \left( \dis\frac{1}{D} + \dis\frac{\sinh (x^{*}/D)}{D \cosh (x^{*}/D)}
       \right ) c_{1} + \dis\frac{\sinh (x^{*}/D)}{D \cosh (x^{*}/D)}
       - \dis\frac{1}{D}  \\
                  \\ \hspace{30mm}
      = - \dis\frac{e^{x^{*}/D}}{D \cosh (x^{*}/D)} c_{1}
        - \dis\frac{e^{-x^{*}/D}}{D \cosh (x^{*}/D)}
    \end{array}
  \end{displaymath}

On the other hand, integrating the equations $w^{0,\pm}(y)$ from $\pm \infty$ to $0$, we have 
  \begin{displaymath}
    \begin{array}{l}
      D^{2} [w_{y}^{0,-}(0) - w_{y}^{0,+}(0)] = 
      - \dis\int_{\infty}^{0} (\tilde{u}^{*}(y+s) + 1) dy 
      + \dis\int_{-\infty}^{0} (\tilde{u}^{*}(y+s) - 1) dy  \\
                 \\ \hspace{33mm}
      = \dis\int_{s}^{\infty} (\tilde{u}^{*}(z) + 1) dz 
        + \dis\int_{-\infty}^{s} (\tilde{u}^{*}(z) - 1) dz  \\
                 \\ \hspace{33mm}
      = -2 s
    \end{array}
  \end{displaymath}
(see also (3.11) and (3.17)). This yields
  \begin{displaymath}
    w_{y}^{0,-}(0) - w_{y}^{0,+}(0) = - \dis\frac{2}{D^{2}} s.
  \end{displaymath}
Then the third equation of (3.27) becomes
  \begin{equation}
    \begin{array}{l}
      0 = \{W_{x}^{1,-}(x^{*}) + w_{y}^{0,-}(0) \} - 
           \{W_{x}^{1,+}(x^{*}) + w_{y}^{0,+}(0) \}  \\
                 \\ \ \ 
      = - \dis\frac{e^{x^{*}/D}}{D \cosh (x^{*}/D)} c_{1} 
        - \dis\frac{2}{D^{2}} s
        - \dis\frac{e^{-x^{*}/D}}{D \cosh (x^{*}/D)}.
    \end{array}
  \end{equation}
Similarly we have the following relation:
  \begin{equation}
    \begin{array}{l}
      0 = \{V_{x}^{1,-}(x^{*}) + v_{y}^{0,-}(0) \} - 
           \{V_{x}^{1,+}(x^{*}) + v_{y}^{0,+}(0) \}  \\
                 \\ \ \ 
      = - \dis\frac{e^{x^{*}}}{\cosh (x^{*})} b_{1} - 2 s
        - \dis\frac{e^{-x^{*}}}{\cosh (x^{*})}.
    \end{array}
  \end{equation}
The above two relations lead to the conclusion (3.30).
\end{proof}

Now we are ready to determine $(s, b_{1}, c_{1})$. Substituting (3.30) and
  \begin{displaymath}
    V_{x}^{0}(x^{*}) = e^{-2 x^{*}} - 1, \ \ \ 
    W_{x}^{0}(x^{*}) = \dis\frac{1}{D} (e^{-2 x^{*}/D} - 1)
  \end{displaymath}
(see the proof of Lemma 3.2) into (3.29), the coefficient of $s$ is given by
  \begin{displaymath}
    \alpha (-2 e^{-2 x^{*}}) 
        + \beta \left ( - \dis\frac{2}{D} e^{-2 x^{*}/D} \right )
    = -2 \left( \alpha e^{-2 x^{*}} 
        + \dis\frac{\beta}{D} e^{-2 x^{*}/D} \right )<0.
  \end{displaymath}
This means that the equation $(3.29) = 0$ can be solved with respect to $s$. Therefore $b_{1}$ and $c_{1}$ are determined by (3.31) and (3.32).

Now we have $C^{1}$-matched solution by using the above parameter values $a_{0}$, $b_{0}$, $b_{1}$, $c_{0}$, $c_{1}$. Applying the implicit function theorem, we obtain the singularly perturbed solution $(u^{\eps}(x), v^{\eps}(x), w^{\eps}(x))$ of (2.1), which leads to the conclusion of Theorem 2.3. 

\section{Stability of standing pulse for $\tau=O(1)$ and $\theta = O(1)$ : Proof of Theorem 2.5}

Stability properties of the pulse solutions are determined by the spectrum of the following linearized problem:
  \begin{equation}
    \left \{
    \begin{array}{l}
       \lambda p = \eps^{2} p_{xx} + f_{u}^{\eps} p 
            - \eps \alpha q - \eps \beta r,  \\
               \\
       \tau \lambda q = q_{xx} + p - q, 
                    \hspace{23mm} x \in {\bf R} \\
               \\
       \theta \lambda r = D^2 r_{xx} + p - r,  \\
               \\
       \dis\lim_{x \rightarrow \pm \infty} p(x) = 0, \ \ 
             \dis\lim_{x \rightarrow \pm \infty} q(x) = 0, \ \ 
             \dis\lim_{x \rightarrow \pm \infty} r(x) = 0,
    \end{array}
    \right.
  \end{equation}
where $f_{u}^{\eps} := f'(u^{\eps}(x))$ and $u^{\eps}(x)$ is $u$-component of the pulse solution.

First we consider the location of the essential spectrum $\sigma_{ess} \{ {\rm (4.1)} \}$ of (4.1). 

\begin{pro}
There exists a negative constant $C_{1}=C_{1}(\tau,\theta)$ independent of $\eps$ such that
  \begin{displaymath}
    \sigma_{ess} \{ {\rm (4.1)}  \} \leq C_{1} < 0
  \end{displaymath}
holds for small $\eps>0$.
\end{pro}

\begin{proof}
The essential spectrum is given by
  \begin{displaymath}
    S = \{ \lambda \in {\bf C} \ | \ \det M = 0, \ -\infty < \xi < \infty \},
  \end{displaymath}
where
  \begin{displaymath}
    M = \left (
    \begin{array}{ccc}
       \lambda + \eps^2 \xi^2 - m & \eps \alpha & \eps \beta  \\
       -1 & \tau \lambda + \xi^{2} + 1 & 0   \\
       -1 & 0 & \theta \lambda + D \xi^{2}+1
    \end{array}
    \right ), 
    \end{displaymath}
    \begin{displaymath}
        m := \lim_{x \rightarrow \pm \infty} f'(u^\eps(x)) = -2< 0.
    \nonumber
  \end{displaymath}
(see Henry \cite{H}). 

Introduce a new variable $\Lambda$ and parameter $\delta$ defined by
  \begin{displaymath}
    \Lambda := \dis\frac{\lambda}{\xi^{2}}, \ \ \ 
        \delta := \dis\frac{1}{\xi^{2}}.
  \end{displaymath}
Then $\det M$ is expressed as 
  \begin{displaymath}
    \begin{array}{l}
      \dis\frac{1}{\xi^{6}} \det M = 
        (\Lambda + \eps^{2} - m \delta)(\tau \Lambda + 1 + \delta)
          (\theta \Lambda + D^{2} + \delta)  \\
                   \\ \hspace{20mm} 
          + \eps \delta ^{2} (\alpha \theta + \beta \tau) \Lambda 
          + \eps \delta^{2} \{ \alpha (D^{2} + \delta) 
             + \beta (1 + \delta) \}.
   \end{array}
  \end{displaymath}
For a fixed constant $\eps_{0}>0$, there exists a small $\delta_{0} = \delta_{0}(\eps_{0})$ such that $\det M = 0$ has three negative roots
  \begin{displaymath}
    \Lambda = - \eps^{2} + m \delta + O(\delta^{2}), \ \ 
    - \dis\frac{1 + \delta}{\tau} + O(\delta^{2}), \ \ 
    - \dis\frac{D^{2} + \delta}{\theta} + O(\delta^{2})
  \end{displaymath}
for $\delta \in (0, \delta_{0})$ and $\eps \in (0, \eps_{0})$. This means that $\det M = 0$ has three roots
  \begin{displaymath}
    \lambda = - \eps^{2} \xi^{2} + m + O(1/\xi^{2}), \ \ 
    - \dis\frac{\xi^{2} + 1}{\tau} + O(1/\xi^{2}), \ \ 
    - \dis\frac{D^{2} \xi^{2} + 1}{\theta} + O(1/\xi^{2})
  \end{displaymath}
for $\xi > 1/\sqrt{\delta_{0}}$ and $\eps \in (0, \eps_{0})$.

When $0 \leq \xi \leq 1/\sqrt{\delta_{0}}$, there exists a small $\eps_{1} = \eps_{1}(\delta_{0}) \ (< \eps_{0})$ such that $\det M = 0$ has three negative roots
  \begin{displaymath}
    \lambda = - \eps^{2} \xi^{2} + m + O(\eps), \ \ 
    - \dis\frac{\xi^{2}+1}{\tau} + O(\eps), \ \ 
    - \dis\frac{D^{2} \xi^{2}+1}{\theta} + O(\eps)
  \end{displaymath}
for $\eps \in (0, \eps_{1})$ and $\xi \in [0, 1/\delta_{0}]$. Therefore if we choose a constant $C_{1}$ as $C_{1} := \max(m, -1/\tau, -1/\theta)$, we can conclude that 
  \begin{displaymath}
    \sigma_{ess} \{ {\rm (4.1)}  \} \leq C_{1} < 0
  \end{displaymath}
for sufficiently small $\eps>0$.
\end{proof}

As was mentioned in section 2, it suffices to study the spectra of ${\rm (EP)}_{ev}$ and ${\rm (EP)}_{od}$. 
Recall that ${\rm (EP)}_{ev}$ (resp. ${\rm (EP)}_{od}$) stands for the eigenvalue problem (2.10), (2.11) and (2.12) (resp. (2.10), (2.11) and (2.13)).
Just to be sure, the proof of Proposition 2.4 is given here (see also \cite{IN}).

\vspace{2mm}

{\it Proof of Proposition 2.4.} For the eigenpair $\{ \lambda, (p(x), q(x), r(x)) \}$ of (4.1), we define $(p_{ev}(x), q_{ev}(x), r_{ev}(x))$ and $(p_{od}(x)$, $q_{od}(x)$, $r_{od}(x))$ by
  \begin{displaymath}
    p_{ev}(x) := \dis\frac{1}{2} \{ p(x) + p(-x) \}, \ 
    q_{ev}(x) := \dis\frac{1}{2} \{ q(x) + q(-x) \}, \ 
    r_{ev}(x) := \dis\frac{1}{2} \{ r(x) + r(-x) \},
  \end{displaymath}
and
  \begin{displaymath}
    p_{od}(x) := \dis\frac{1}{2} \{ p(x) - p(-x) \}, \ \ 
    q_{od}(x) := \dis\frac{1}{2} \{ q(x) - q(-x) \}, \ \ 
    r_{od}(x) := \dis\frac{1}{2} \{ r(x) - r(-x) \}.
  \end{displaymath}
We can easily show that $\{ \lambda, (p_{ev}(x), q_{ev}(x), r_{ev}(x)) \}$ and $\{ \lambda, (p_{od}(x), q_{od}(x), r_{od}(x)) \}$ are the eigenpairs of ${\rm (EP)}_{ev}$ and ${\rm (EP)}_{od}$ respectively.

If $\{\lambda, (p_{ev}(x), q_{ev}(x), r_{ev}(x)) \}$ is an eigenpair of ${\rm (EP)}_{ev}$, we define $(p(x), q(x), r(x))$ by
  \begin{displaymath}
    p(x):= p_{ev}(|x|), \ \ q(x):= q_{ev}(|x|), \ \ r(x):= r_{ev}(|x|) \ \ 
      x \in {\bf R}.
  \end{displaymath}
Then $\{\lambda, (p(x), q(x), r(x)) \}$ is an eigenpair of (3.1) and $(p(x), q(x), r(x))$ are even functions. On the other hand, let $(p(x), q(x), r(x))$ be
  \begin{displaymath}
    p(x):= {\rm sgn}(x) p_{od}(|x|), \ \ q(x):= {\rm sgn}(x) q_{od}(|x|), \ \ 
    r(x):= {\rm sgn}(x) r_{od}(|x|) \ \ x \in {\bf R}.
  \end{displaymath}
Then $\{\lambda, (p(x), q(x), r(x)) \}$ is an eigenpair of (4.1) and $(p(x), q(x), r(x))$ are odd functions. \hfill $\Box$

In order to study the spectral distribution of ${\rm (EP)}_{ev}$ and ${\rm (EP)}_{od}$, the following singular Sturm-Liouville problem plays a key role.
  \begin{equation}
    L^{\eps} \phi := \eps^{2} \phi_{xx} + f_{u}^{\eps} \phi 
      = \zeta \phi, \ \ x \in I = (0, \infty)
  \end{equation}
  \begin{equation}
    \dis\lim_{x \rightarrow \infty} \phi (x) = 0,
  \end{equation}
and
  \begin{displaymath}
    \phi_{x}(0) = 0 \ \ \mbox{or} \ \ \phi(0) = 0.
  \end{displaymath}
Let ${\rm (SL)}_{ev}$ (resp. ${\rm (SL)}_{od}$) be the singular Sturm-Liouville problem (4.2), (4.3) and $\phi_{x}(0) = 0$ (resp. $\phi(0) = 0$), and $\{ \zeta_{n,N}^{\eps}, \phi_{n,N}^{\eps} \}_{n \geq 0}$ 
(resp. $\{ \zeta_{n,D}^{\eps}, \phi_{n,D}^{\eps} \}_{n \geq 0}$) be the complete orthonormal set of eigenvalues and eigenfunctions of ${\rm (SL)}_{ev}$ (resp. ${\rm (SL)}_{od}$) in $L^{2}(I)$-sense.

We concentrate ${\rm (SL)}_{ev}$ for a while. The following two eigenvalue problems are useful for later use. One is the stretched Sturm-Liouville problem $(\widetilde{\rm SL})_{ev}$ for ${\rm (SL)}_{ev}$ centered at the layer position $x^*$:
  \begin{displaymath}
    \begin{array}{l}
      \dis\frac{d^{2}}{dy^{2}} \hat{\phi} + \tilde{f}_{u}^{\eps} \hat{\phi}
           = \zeta \hat{\phi}, \hspace{5mm} y \in (-x^{*}/\eps, \infty), \\
              \\
      \hat{\phi}_{y}(-x^{*}/\eps) = 0, \ \ \ 
      \dis\lim_{y \rightarrow \infty} \hat{\phi}(y) = 0.
    \end{array}
  \end{displaymath}
The other is the limiting eigenvalue problem ${\rm (SL)}_{ev}^{*}$:
  \begin{displaymath}
    \dis\frac{d^{2}}{dy^{2}} \hat{\phi} + \tilde{f}_{u}^{*} \hat{\phi}
       = \zeta \hat{\phi}, \hspace{5mm} y \in {\bf R}, 
  \end{displaymath}
where $f_{u}^{*} = f'(u^{*}(y))$.

It is easy to see the following lemma for $(\widetilde{\rm SL})_{ev}$ and ${\rm (SL)}_{ev}^{*}$.

\begin{lem}
{\rm (i)} The value $0$ is the principal eigenvalue of ${\rm (SL)}_{ev}^{*}$, which is simple in $L^{2}({\bf R})$, and the associated normalized positive eigenfunction $\hat{\phi}_{0,N}^{*}$ is given by
  \begin{displaymath}
    \hat{\phi}_{0,N}^{*} = \kappa^{*} \tilde{u}_{y}^{*},
  \end{displaymath}
where $\kappa^{*} := \| \tilde{u}_{y}^{*} \|_{L^{2}({\bf R})}^{-1}$.

\vspace{2mm}

{\rm (ii)} Let $\hat{\phi}_{0,N}^{\eps} (=\sqrt{\eps} \tilde{\phi}_{0,N}^{\eps}) $ be the normalized positive principal eigenfunction of $(\widetilde{\rm SL})_{ev}$. Then it follows that
  \begin{displaymath}
     \lim_{\eps \downarrow 0} \hat{\phi}_{0,N}^{\eps} = \hat{\phi}_{0,N}^{*} 
         \hspace{5mm} \mbox{in} \ C_{c.u.}^{2}({\rm \bf R}) \mbox{-sense}.
  \end{displaymath}
\end{lem}

\begin{proof}
See Lemma 1.2 and 1.3 in \cite{NF} for the proof.
\end{proof}

Now we summarize the spectral behavior of (4.2) and (4.3).

\begin{lem}[spectral properties of $L^{\eps}$] 

{\rm (i)} Let $\sigma_{ess} \{ {\rm(SL)}_{ev} \}$ be the essential spectrum of $L^{\eps}$. 
Then, there exists a positive constant $\mu_{0}$ independent of $\eps$ such that
  \begin{displaymath}
    \sigma_{ess} \{ {\rm(SL)}_{ev} \} < -\mu_{0} < 0
  \end{displaymath}
holds.

{\rm (ii)} It holds that 
  \begin{displaymath}
    \zeta_{0,N}^{\eps} > 0 > -\Delta^{*} 
         > \zeta_{1,N}^{\eps} > \zeta_{2,N}^{\eps} > \cdots,
  \end{displaymath}
for small $\eps>0$, where
  \begin{equation}
    \begin{array}{rcl}
      \hat{\zeta}_{0,N}^{*} &:=& 
         \dis\lim_{\eps \downarrow 0} \frac{\zeta_{0,N}^{\eps}}{\eps^{2}}
      = - 2 (\kappa^{*})^{2} 
            [\alpha V_{x}^{0}(x^{*}) + \beta W_{x}^{0}(x^{*})] \\
                  \\
      &=& \dis\frac{3 \sqrt{2}}{2} 
        \left[
           \alpha (1 - e^{-2 x^{*}}) + \dis\frac{\beta}{D} (1 - e^{-2 x^{*}/D})
        \right]
    \end{array}
  \end{equation}
and $\Delta^{*}$ is independent of $\eps$. Here we used the fact that
  \begin{displaymath}
    (\kappa^{*})^{2} = \dis\frac{3 \sqrt{2}}{4}, \ \ \ 
    V_{x}^{0}(x^{*}) = e^{-2 x^{*}} - 1, \ \ \ 
    W_{x}^{0}(x^{*}) = \dis\frac{1}{D} (e^{-2 x^{*}/D} - 1).
  \end{displaymath}
\end{lem}

\begin{proof}
See Proposition 3.1 of \cite{NMIF} and Lemma 1.4 of \cite{NF} for the proof. Here we show the outline for the derivation of (4.4). For the principal eigenfunction $\phi_{0,N}^{\eps}(x)$, we first define a new function:
  \begin{displaymath}
     \hat{\phi}_{0,N}^{\eps}(y) := 
            \sqrt{\eps} \phi_{0,N}^{\eps}(x^{*} + \eps y), \ \ \ 
            y \in (-x^{*}/\eps, \infty ).
  \end{displaymath}
Then $\hat{\phi}_{0,N}^{\eps}$ satisfies
  \begin{equation}
    \dis\frac{d^{2}}{dy^{2}} \hat{\phi}_{0,N}^{\eps} 
           +  \tilde{f}_{u}^{\eps} \hat{\phi}_{0,N}^{\eps}
          = \zeta_{0,N}^{\eps} \hat{\phi}_{0,N}^{\eps}.
  \end{equation}
On the other hand, the stretched equation of the first one to (2.1) becomes
  \begin{equation}
    \dis\frac{d^{2}}{dy^{2}} \tilde{u} + f(\tilde{u}) 
     - \eps (\alpha \tilde{v} + \beta \tilde{w} + \gamma) = 0,
  \end{equation}
where $\tilde{u} = \tilde{u}(y) := u(x^{*} + \eps y)$ and so on. Differentiating (4.6) with respect to $y$, we have
  \begin{equation}
    \dis\frac{d^{2}}{dy^{2}} \tilde{u}_{y} + \tilde{f}_{u}^{\eps} \tilde{u}_{y}
     - \eps (\alpha \tilde{v}_{y} + \beta \tilde{w}_{y}) = 0,
  \end{equation}

Multiplying $\tilde{u}_{y}$ to (4.5) and integrating with respect to $y$ from $0$ to $\infty$, we have
  \begin{displaymath}
    \left \langle
       \dis\frac{d^{2}}{dy^{2}} \hat{\phi}_{0,N}^{\eps}, \tilde{u}_{y}
    \right \rangle + 
    \left \langle
       \tilde{f}_{u}^{\eps} \hat{\phi}_{0,N}^{\eps}, \tilde{u}_{y}
    \right \rangle = \zeta_{0,N}^{\eps} 
    \left \langle
       \hat{\phi}_{0,N}^{\eps}, \tilde{u}_{y}
    \right \rangle.
  \end{displaymath}
Applying the integration by parts twice, and using (4.7), we have
  \begin{equation}
    - \left. \hat{\phi}_{0,N}^{\eps} \tilde{u}_{yy} 
      \right|_{-x^{*}/\eps}^{\infty} +
    \eps \left \langle 
        \hat{\phi}_{0,N}^{\eps}, \alpha \tilde{v}_{y} + \beta \tilde{w}_{y}
    \right \rangle = \zeta_{0,N}^{\eps} 
    \left \langle
       \hat{\phi}_{0,N}^{\eps}, \tilde{u}_{y}
    \right \rangle.
  \end{equation}
In the same manner as in the proof of Lemma 1.4 of \cite{NF}, we can show that
  \begin{displaymath}
    \lim_{\eps \downarrow 0} \left \langle \hat{\phi}_{0,N}^{\eps}, 
         \tilde{u}_{y} \right \rangle 
         = \dis\frac{1}{\kappa^{*}}, \hspace{5mm}
    \left| \mbox{the first term of (4.8)} \right| \leq C \exp (-c_{0}/\eps),
  \end{displaymath}
  \begin{displaymath}
    \lim_{\eps \downarrow 0} \dis\frac{1}{\eps} 
        \left \langle \hat{\phi}_{0,N}^{\eps}, \tilde{v}_{y} \right \rangle
    = \kappa^{*} V_{x}^{0}(x^{*}) 
        \int_{-\infty}^{\infty} \tilde{u}_{y}^{*} dy 
    = - 2 \kappa^{*} V_{x}^{0}(x^{*}),
  \end{displaymath}
  \begin{displaymath}
    \lim_{\eps \downarrow 0} \dis\frac{1}{\eps} 
        \left \langle \hat{\phi}_{0,N}^{\eps}, \tilde{w}_{y} \right \rangle
    = - 2 \kappa^{*} W_{x}^{0}(x^{*}),
    \nonumber
  \end{displaymath}
where $C>0$ and $c_{0}>0$ are constants independent of $\eps$. 
\end{proof}

In the following, we consider the spectra of (2.10) which lie in the region $\Lambda_{0}$ defined by
  \begin{displaymath}
     \Lambda_{0} := \{ \lambda \in {\bf \rm C} \ | \ {\rm Re} \lambda > 
                  \max \{ -\Delta^{*}, -\mu_{0} \} \}.
  \end{displaymath}

The resolvent operator $(L^{\eps} - \lambda)^{-1}$ plays an important role in the following subsection. It has the eigenfunction expansion
  \begin{displaymath}
    (L^{\eps} - \lambda)^{-1} h = 
    \dis \frac{\langle h, \ \phi_{0,N}^{\eps} \rangle}
            {\zeta_{0,N}^{\eps} - \lambda} \phi_{0,N}^{\eps} 
           + (L^{\eps} - \lambda)^{\dag}h
  \end{displaymath}
where
  \begin{displaymath}
    (L^{\eps} - \lambda)^{\dag}(\ \cdot \ ) := \sum_{n \geq 1} 
       \frac{\langle \ \cdot \ ,\ \phi_{n,N}^{\eps} \rangle}
           {\zeta_{n,N}^{\eps} - \lambda} \phi_{n,N}^{\eps}.
  \end{displaymath}
$(L^{\eps} - \lambda)^{\dag}$ is a uniformly bounded operator from $L^{2}(I)$ into itself, i.e.,
  \begin{displaymath}
    \pl (L^{\eps} - \lambda)^{\dag} \pl_{{\cal L}(L^{2}(I), L^{2}(I))} \leq M
  \end{displaymath}
for any small $\eps$ and $\lambda \in \Lambda_{0}$ with an appropriate $M>0$ independent of $\eps$. The next two lemmas are crucial to derive a singular limiting eigenvalue problem called the SLEP equation.

\begin{lem}
For small $\eps$, it holds that
  \begin{displaymath}
    \zeta_{0,N}^{\eps} \notin \sigma \{ {\rm (EP)}_{ev} \},
  \end{displaymath}
where $\sigma \{ {\rm (EP)}_{ev} \}$ denotes the set of spectra to ${\rm (EP)}_{ev}$. 
\end{lem}

\begin{lem}
Let $L^{\eps}$ be the Sturm-Liouville operator in ${\rm (SL)}_{ev}$. For $h \in L^{2}(I) \cap L^{\infty}(I)$ and $\lambda \in {\bf C} \backslash \sigma(L^{\eps})$,
  \begin{displaymath}
    \lim_{\eps \downarrow 0}(L^{\eps} - \lambda)^{\dag}(h) = 
    \frac{h}{f_{u}^{*} - \lambda},
    \ \ \mbox{strongly in } L^{2} \mbox{-sense},
  \end{displaymath}
where $f_{u}^{*} := \lim_{\eps \downarrow 0} f_u^\eps (u^\eps(x))= -2$.
\end{lem}

\begin{lem}
  \begin{displaymath}
    \lim_{\eps \downarrow 0} \frac{1}{\sqrt{\eps}} \phi_{0,N}^{\eps} 
    = 2 \kappa^{*} \delta_{*} \ \ \ H^{-1}(I) \mbox{-sense}
  \end{displaymath}
where $\delta_{*}$ is a Dirac's $\delta$-function at $x^{*}$.
\end{lem}

{\it Proof of Lemmas 4.4, 4.5 and 4.6}. This can be done in a similar way to those of Lemmas 2.1, 2.2 and 2.3 of \cite{NF}.  \hfill $\Box$

\vspace{2mm}

\begin{rmk}
Lemmas 4.2, 4.3, 4.4, 4.5 and 4.6 hold for ${\rm (SL)}_{od}$. Especially the asymptotic limit of the coefficient $\hat{\zeta}_{0,D}^{*}$ of $\eps^{2}$ does not depend on the boundary condition. That is
  \begin{displaymath}
    \hat{\zeta}_{0,D}^{*} := 
      \dis\lim_{\eps \downarrow 0} \frac{\zeta_{0,D}^{\eps}}{\eps^{2}} 
    = \hat{\zeta}_{0,N}^{*},
  \end{displaymath}
where $\zeta_{0,D}^{\eps}$ is the principal eigenvalue of ${\rm (SL)}_{od}$. For more details, see \cite{NF}.
\end{rmk}

\vspace{2mm}

Since the computation below has no essential difference between even and odd cases, we suppress the subscript "$N$" or "$D$" for a while. We solve the first equation of (2.10) with respect to $p$ for $\lambda \in \Lambda_{0}$ as
  \begin{equation}
    \begin{array}{rcl}
      p &=& \eps (L^{\eps}-\lambda)^{-1} (\alpha q + \beta r) \\
              \\
        &=& \dis \frac{\eps \langle \alpha q 
           + \beta r,\ \phi_{0}^{\eps} \rangle}
            {\zeta_{0}^{\eps} - \lambda} \phi_{0}^{\eps} 
           + \eps (L^{\eps} - \lambda)^{\dag}(\alpha q + \beta r).
    \end{array}
  \end{equation}
Taking account of Lemmas 4.3, 4.5 and 4.6, we rewrite (4.9) as
  \begin{equation}
    p = \dis \frac{\langle \alpha q 
         + \beta r,\ \phi_{0}^{\eps}/\sqrt{\eps} \rangle}
         {\zeta_{0}^{\eps}/\eps^{2} - \lambda/\eps^{2}} 
          \phi_{0}^{\eps}/\sqrt{\eps} 
       + \eps (L^{\eps} - \lambda)^{\dag}(\alpha q + \beta r).
  \end{equation}
Noting that $\phi_{0}^{\eps}/\sqrt{\eps}$ is a nice function that converges to a constant-multiple of Dirac-$\delta$'s function in weak sense and $\zeta_{0}^{\eps} = O(\eps^{2})$, it is not difficult to see that (4.10) has an nontrivial well-defined limit when $|\lambda|=O(\eps^{2})$ as $\eps \downarrow 0$. In fact, all the other cases are not important as shown below. Also $\lambda$ may be a complex value, however it turns out that there do not appear complex eigenvalues that affect the stability in the regime of $\tau=O(1)$ and $\theta = O(1)$. In fact the following argument is valid for complex 
$\lambda$, however the relevant ones become real as a conclusion. We classify the problem into five cases.

(i) $\lambda=O(1)$, 

(ii) $\lambda = o(1)$ with $\dis\lim_{\eps \downarrow 0} |\lambda|/\eps^{2} = \infty$,

(iii) $\lambda=O(\eps^{2})$,

(iv) $\dis\lim_{\eps \downarrow 0} |\lambda|/\eps^{2} = 0$ as $\eps \downarrow 0$,

(v) $\dis\lim_{\eps \downarrow 0} |\lambda| = \infty$.

\vspace{2mm}

Case (i) $\lambda=O(1)$: Since the denominator of the first term in (4.10) diverges as $\eps \downarrow 0$, we obtain the limiting function (i.e. the principal part of $p$) $p \equiv 0$. Then the limiting system for $(q,r)$ becomes
  \begin{displaymath}
    \left \{
    \begin{array}{rcl}
       - q_{xx} + q + \tau \lambda q = 0  \\
                    \\
       - D^{2} r_{xx} + r + \theta \lambda r = 0
    \end{array}
    \right.
  \end{displaymath}
If ${\bf Re} \ \lambda > \max(-1/\tau, -1/\theta)$, then $(q, r) \equiv (0, 0)$. This means that $\lambda$ such that ${\bf Re} \ \lambda >0$ belongs to the resolvent set. 

\vspace{2mm}

Case (ii) $\lambda = o(1)$ with $\dis\lim_{\eps \downarrow 0} |\lambda|/\eps^{2} = \infty$: $\lambda$ tends to zero slower than $\eps^{2}$. The denominator of the first term of (4.10) diverges and $\lambda=o(1)$ as $\eps \downarrow 0$. Then the limiting system for the principal part of $(q,r)$ becomes
  \begin{displaymath}
    \left \{
    \begin{array}{rcl}
       - q_{xx} + q = 0  \\
                    \\
       - D^{2} r_{xx} + r = 0
    \end{array}
    \right.
  \end{displaymath}
This yields that $(q, r) \equiv (0, 0)$, hence $p \equiv 0$. Such an eigenpair does not exist.

\vspace{2mm}

Case (iii) $\lambda=O(\eps^{2})$ : 
We set the asymptotic form of $\lambda(\eps)$ as $\lambda(\eps) = \eps^{2} \hat{\lambda}(\eps)$. First we consider when the limiting value $\hat{\lambda}^{*}$ of $\hat{\lambda}(\eps)$ exists, that is
  \begin{displaymath}
    \hat{\lambda}^{*} = \dis\lim_{\eps \downarrow 0} \hat{\lambda}(\eps).
  \end{displaymath}
Substitute (4.10) into the second and the third equations of (2.10), we have a closed eigenvalue problem with respect to $(q, r)$:
  \begin{equation}
    \left \{
    \begin{array}{l}
       - q_{xx} + q + \tau \lambda q 
              - \eps \alpha (L^{\eps} - \lambda)^{\dag}(q)  \\
                    \\ \hspace{10mm}
       =\dis \frac{\langle \alpha q + \beta r, 
          \ \phi_{0,N}^{\eps}/\sqrt{\eps} \rangle}
          {\zeta_{0,N}^{\eps}/\eps^{2} - \lambda/\eps^{2}} 
          \phi_{0,N}^{\eps}/\sqrt{\eps} 
           + \eps \beta (L^{\eps} - \lambda)^{\dag}(r),  \\
               \\
       - D^{2} r_{xx} + r + \theta \lambda r 
              - \eps \beta (L^{\eps} - \lambda)^{\dag}(r) \\
                     \\ \hspace{10mm}
       = \dis \frac{\langle \alpha q + \beta r, 
          \ \phi_{0,N}^{\eps}/\sqrt{\eps} \rangle}
          {\zeta_{0,N}^{\eps}/\eps^{2} - \lambda/\eps^{2}} 
          \phi_{0,N}^{\eps}/\sqrt{\eps} 
           + \eps \alpha (L^{\eps} - \lambda)^{\dag}(q).
    \end{array}
    \right.
  \end{equation}
Note that (4.11) should be written in a weak form. But for notational simplicity, we write them in a classical form. Now we define two operators $K_{q}^{\eps}$ and $K_{r}^{\eps}$ by
  \begin{displaymath}
    K_{q}^{\eps} := \left \{ -\dis\frac{d^{2}}{d x^{2}} + 1 + \tau \lambda
          - \eps \alpha (L^{\eps} - \lambda)^{\dag}(\cdot) \right \}^{-1},
  \end{displaymath}
  \begin{displaymath}
    K_{r}^{\eps} := \left \{ - D \dis\frac{d^{2}}{d x^{2}} + 1 + \theta \lambda
          - \eps \beta (L^{\eps} - \lambda)^{\dag}(\cdot) \right \}^{-1}
  \end{displaymath}
with suitable boundary conditions at $x=0$. Then (4.11) is recast as
  \begin{equation}
    \left \{
    \begin{array}{rcl}
       q &=& 
       \dis \frac{\langle \alpha q + \beta r, 
          \ \phi_{0}^{\eps}/\sqrt{\eps} \rangle}
          {\zeta_{0}^{\eps}/\eps^{2} - \lambda/\eps^{2}} 
          K_{q}^{\eps}(\phi_{0}^{\eps}/\sqrt{\eps}) 
           + \eps \beta K_{q}^{\eps}((L^{\eps} - \lambda)^{\dag}(r))  \\
               \\
       r &=& 
       \dis \frac{\langle \alpha q + \beta r, 
          \ \phi_{0}^{\eps}/\sqrt{\eps} \rangle}
          {\zeta_{0}^{\eps}/\eps^{2} - \lambda/\eps^{2}} 
          K_{r}^{\eps}(\phi_{0}^{\eps}/\sqrt{\eps}) 
           + \eps \alpha K_{r}^{\eps}((L^{\eps} - \lambda)^{\dag}(q)).
    \end{array}
    \right.
  \end{equation}
Making use of Lemmas 4.4, 4.5 and 4.6, we can take the limit of (4.12) as $\eps \downarrow 0$. That is
  \begin{equation}
    \left \{
    \begin{array}{rcl}
      q^{*} &=& 
      4 (\kappa^{*})^{2} \alpha \dis \frac{\langle q^{*}, \ \delta_{*} \rangle}
          {\hat{\zeta}_{0}^{*} - \hat{\lambda}^{*}} K_{q}^{*}(\delta_{*}) + 
      4 (\kappa^{*})^{2} \beta \dis \frac{\langle r^{*}, \ \delta_{*} \rangle}
          {\hat{\zeta}_{0}^{*} - \hat{\lambda}^{*}} K_{q}^{*}(\delta_{*})    \\
               \\
      r^{*} &=& 
      4 (\kappa^{*})^{2} \alpha \dis \frac{\langle q^{*}, \ \delta_{*} \rangle}
          {\hat{\zeta}_{0}^{*} - \hat{\lambda}^{*}} K_{r}^{*}(\delta_{*}) + 
      4 (\kappa^{*})^{2} \beta \dis \frac{\langle r^{*}, \ \delta_{*} \rangle}
          {\hat{\zeta}_{0}^{*} - \hat{\lambda}^{*}} K_{r}^{*}(\delta_{*})
    \end{array}
    \right.
  \end{equation}
where
  \begin{displaymath}
    K_{q}^{*} := \left \{ -\dis\frac{d^{2}}{d x^{2}} + 1 \right \}^{-1}, \ \ 
    K_{r}^{*} := \left \{ -D\dis\frac{d^{2}}{d x^{2}} + 1 \right \}^{-1}.
  \end{displaymath}
(4.13) shows that $q^{*}$ (resp. $r^{*}$) is a constant multiple of $K_{q}^{*}(\delta_{*})$ (resp. $K_{r}^{*}(\delta_{*})$). So we set
  \begin{equation}
    q^{*} = A K_{q}^{*}(\delta_{*}), \ \ r^{*} = B K_{r}^{*}(\delta_{*})
  \end{equation}
and substitute (4.14) into (4.13). Then we have
  \begin{displaymath}
    \left \{
    \begin{array}{rcl}
       A K_{q}^{*}(\delta_{*}) &=& 
       4 (\kappa^{*})^{2} A \alpha \dis \frac{\langle K_{q}^{*}(\delta_{*}), \ 
           \delta_{*} \rangle}
          {\hat{\zeta}_{0}^{*} - \hat{\lambda}^{*}} K_{q}^{*}(\delta_{*}) + 
       4 (\kappa^{*})^{2} B \beta \dis \frac{\langle K_{r}^{*}(\delta_{*}), \ 
           \delta_{*} \rangle}
          {\hat{\zeta}_{0}^{*} - \hat{\lambda}^{*}} K_{q}^{*}(\delta_{*})    \\
               \\
       B K_{r}^{*}(\delta_{*}) &=& 
       4 (\kappa^{*})^{2} A \alpha \dis \frac{\langle K_{q}^{*}(\delta_{*}), \ 
           \delta_{*} \rangle}
          {\hat{\zeta}_{0}^{*} - \hat{\lambda}^{*}} K_{r}^{*}(\delta_{*}) + 
       4 (\kappa^{*})^{2} B \beta \dis \frac{\langle K_{r}^{*}(\delta_{*}), \ 
           \delta_{*} \rangle}
          {\hat{\zeta}_{0}^{*} - \hat{\lambda}^{*}} K_{r}^{*}(\delta_{*}).
    \end{array}
    \right.
  \end{displaymath}
This is equivalent to
  \begin{equation}
    \left(
    \begin{array}{cc}
      \hat{\zeta}_{0}^{*} - \hat{\lambda}^{*} 
      - 4 (\kappa^{*})^{2} \alpha \langle K_{q}^{*}(\delta_{*}), \ \delta_{*} 
           \rangle &
      - 4 (\kappa^{*})^{2} \beta \langle K_{r}^{*}(\delta_{*}), \ \delta_{*} 
           \rangle \\
                \\
      - 4 (\kappa^{*})^{2} \alpha \langle K_{q}^{*}(\delta_{*}), \ \delta_{*} 
           \rangle &
      \hat{\zeta}_{0}^{*} - \hat{\lambda}^{*} 
      - 4 (\kappa^{*})^{2} \beta \langle K_{r}^{*}(\delta_{*}), \ \delta_{*} 
           \rangle
    \end{array}
    \right)
    \left(
    \begin{array}{c}
       A \\
         \\
       B
    \end{array}
    \right)
  \end{equation}
  \begin{displaymath}
    = \left(
    \begin{array}{c}
       0 \\
         \\
       0
    \end{array}
    \right).
  \end{displaymath}
In order to have a nontrivial solution for $(A,B)^{T}$, the determinant must be equal to zero, i.e., 
  \begin{displaymath}
    (\hat{\zeta}_{0}^{*} - \hat{\lambda}^{*}) 
    (\hat{\zeta}_{0}^{*} - \hat{\lambda}^{*} - 4 (\kappa^{*})^{2} [
        \alpha \langle K_{q}^{*}(\delta_{*}), \ \delta_{*} \rangle
      + \beta \langle K_{r}^{*}(\delta_{*}), \ \delta_{*} \rangle ]) = 0.
  \end{displaymath}
Using the fact that $\hat{\zeta}_{0}^{*} - \hat{\lambda}^{*} \neq 0$ (see Lemma 3.4), we obtain
  \begin{equation}
    \hat{\zeta}_{0}^{*} - \hat{\lambda}^{*} = 4 (\kappa^{*})^{2} [
        \alpha \langle K_{q}^{*}(\delta_{*}), \ \delta_{*} \rangle
      + \beta \langle K_{r}^{*}(\delta_{*}), \ \delta_{*} \rangle ].
  \end{equation}
We call (4.16) {\it the} SLEP {\it equation} of (2.10). Once we find a solution $\hat{\lambda}^{*}$ of the SLEP equation (4.16), we can show that the existence of the eigenvalues $\lambda^{\eps}$ of the original eigenvalue problem (2.10) for positive $\eps$ such that
$\dis\lim_{\eps \downarrow 0} \lambda^{\eps} / \eps^{2} = \hat{\lambda}^{*}$ 
by applying the implicit function theorem (see \cite{N}, \cite{NF}). 
Computing the difference between the first and the second equations of (4.15), we have
  \[
     (\hat{\zeta}_{0}^{*} - \hat{\lambda}^{*})(A - B) = 0.
  \]
This yields $A=B$.
Then the limiting function $p^{*}$ of $p$ as $\eps \downarrow 0$ is given by
  \begin{displaymath}
    p^{*} = A \dis \frac{\langle \alpha q^{*} 
       + \beta r^{*} ,\ \delta_{*} \rangle}
         {\hat{\zeta}_{0}^{*} - \hat{\lambda}^{*}} \delta_{*}
       = - 4 (\kappa^{*})^{2} A \ \delta_{*}
  \end{displaymath}
That is, the asymptotic eigenspace associated with the eigenvalue $\lambda \approx \eps^{2} \hat{\lambda}^{*}$ is spanned by
  \begin{displaymath}
    (- 4 (\kappa^{*})^{2} \delta_{*}, \ K_{q}^{*}(\delta_{*}), 
      \ K_{r}^{*}(\delta_{*}))
  \end{displaymath}
and its dimension is one. 

Now we are ready to study the precise behavior of eigenvalues for the case (iii). Let $K_{q}^{e,*}$ (resp. $K_{r}^{e,*}$) be $K_{q}^{*}$ (resp.  $K_{r}^{*}$) with Neumann boundary condition at $x=0$ (i.e. the symmetric mode), and $K_{q}^{o,*}$ (resp. $K_{r}^{o,*}$) be $K_{q}^{*}$ (resp.  $K_{r}^{*}$) with Dirichlet boundary condition at $x=0$ (i.e. the odd symmetric mode). Then we have the following lemma.

\begin{lem}
  \begin{equation}
    \langle K_{q}^{e,*}(\delta_{*}), \ \delta_{*} \rangle 
            = \dis\frac{\cosh (x^{*})}{e^{x^{*}}}, \ \ \ 
    \langle K_{r}^{e,*}(\delta_{*}), \ \delta_{*} \rangle 
            = \dis\frac{\cosh (x^{*}/D)}{D e^{x^{*}/D}}
  \end{equation}
  \begin{equation}
    \langle K_{q}^{o,*}(\delta_{*}), \ \delta_{*} \rangle 
            = \dis\frac{\sinh (x^{*})}{e^{x^{*}}}, \ \ \ 
    \langle K_{r}^{o,*}(\delta_{*}), \ \delta_{*} \rangle 
            = \dis\frac{\sinh (x^{*}/D)}{D e^{x^{*}/D}}
  \end{equation}
\end{lem}

\begin{proof} \ See Appendix A. 
\end{proof}

\begin{rmk}
It is easy to see $\hat{\zeta}_{0}^{*}$ has the following form:
  \begin{equation}
    \hat{\zeta}_{0}^{*} = 4 (\kappa^{*})^{2} [
    \alpha \langle K_{q}^{o,*}(\delta_{*}), \ \delta_{*} \rangle 
    + \beta \langle K_{r}^{o,*}(\delta_{*}), \ \delta_{*} \rangle ].
  \end{equation}
\end{rmk}

For the symmetric mode, we have
  \begin{displaymath}
    \begin{array}{rcl}
      \hat{\lambda}^{*} &=& \hat{\zeta}_{0}^{*} - 4 (\kappa^{*})^{2} 
         [\alpha \langle K_{q}^{e,*}(\delta_{*}), \ \delta_{*} \rangle
          + \beta \langle K_{r}^{e,*}(\delta_{*}), \ \delta_{*} \rangle ] \\
                   \\
      &=& 2 (\kappa^{*})^{2} \left[ 
       \alpha (1 - e^{-2 x^{*}}) + \dis\frac{\beta}{D} (1 - e^{-2 x^{*}/D})
       \right]   \\
                   \\
       & & \hspace{20mm} - 4 (\kappa^{*})^{2}
         \left[ \alpha \dis\frac{\cosh (x^{*})}{e^{x^{*}}}
         + \beta \dis\frac{\cosh (x^{*}/D)}{D e^{x^{*}/D}} \right] \\
                   \\
      &=& - 3 \sqrt{2} 
       \left[ \alpha e^{-2 x^{*}} + \dis\frac{\beta}{D} e^{-2 x^{*}/D} \right]
    \end{array}
  \end{displaymath}
Here we used the fact that $(\kappa^{*})^{2}=3\sqrt{2}/4$ (see Lemma 4.3). 

Similarly the odd symmetric case yields $\hat{\lambda}^{*}=0$. That is, 
  \begin{displaymath}
    \begin{array}{rcl}
      \hat{\lambda}^{*} &=& \hat{\zeta}_{0}^{*} - 4 (\kappa^{*})^{2} 
         [\alpha \langle K_{q}^{o,*}(\delta_{*}), \ \delta_{*} \rangle
         + \beta \langle K_{r}^{o,*}(\delta_{*}), \ \delta_{*} \rangle ] \\
                   \\
      &=& 2 (\kappa^{*})^{2} \left[ 
       \alpha (1 - e^{-2 x^{*}}) + \dis\frac{\beta}{D} (1 - e^{-2 x^{*}/D})
       \right] \\
                    \\
       & & \hspace{20mm} - 4 (\kappa^{*})^{2}
         \left[ \alpha \dis\frac{\sinh (x^{*})}{e^{x^{*}}}
         + \beta \dis\frac{\sinh (x^{*}/D)}{D e^{x^{*}/D}} \right] \\
      &=& 0
    \end{array}
  \end{displaymath}
This recovers the zero eigenvalue coming from the translation invariance.

So far we assume that $\hat{\lambda}(\eps)$ has a definite limit, however, in general, it is bounded but may not have the limit. It turns out that even such a case can be reduced to the above discussion. Namely, for any convergent subsequence $\{ \eps_{k} \}$ with $\lim_{k \rightarrow \infty} \eps_{k} = 0$ such that $\lim_{k \rightarrow \infty} \hat{\lambda}(\eps_{k})$ exists, and the limit is a solution of (4.16). It is one of values discussed above. Any subsequence converges to the same value so that uniqueness result implies that this case does not happen.

\vspace{2mm}

Case (iv): $\dis\lim_{\eps \downarrow 0} |\lambda| /\eps^{2} = 0$: $\lambda$ tends to zero faster that $\eps^{2}$. In this case, the limiting eigenvalue problem is given by (4.13) with $\hat{\lambda}^{*}=0$. Therefore we recover the zero eigenvalue exactly same discussion as in the case (iii).

\vspace{2mm}

Case (v): It is easy to check that the solution 
of the limiting system is $(p,q,r) =(0,0,0)$, which implies that such an eigenvalue does not exist.

\vspace{2mm}

This completes the proof of Theorem 2.5.

\section{Stability of standing pulse for $\tau=O(1/\eps^{2})$ and $\theta=O(1/\eps^{2})$}

In this section, we consider the case $\tau=O(1/\eps^{2})$ and $\theta=O(1/\eps^{2})$. Namely the relaxation time of two inhibitors is very slow. After setting $\tau = \hat{\tau}/\eps^{2}$, $\theta = \hat{\theta}/\eps^{2}$, the linearized eigenvalue problem (4.1) is recast as
  \begin{equation}
    \left \{
    \begin{array}{l}
       \lambda p = \eps^{2} p_{xx} + f_{u} p 
            - \eps \alpha q - \eps \beta r,  \\
               \\
       \dis\frac{1}{\eps^2} \hat{\tau} \lambda q = q_{xx} + p - q,
                    \hspace{10mm} x \in {\bf R} \\
               \\
       \dis\frac{1}{\eps^2} \hat{\theta} \lambda r = D^2 r_{xx} + p - r,  \\
               \\
       \dis\lim_{x \rightarrow \pm \infty} p(x) = 0, \ \ 
               \dis\lim_{x \rightarrow \pm \infty} q(x) = 0, \ \ 
               \dis\lim_{x \rightarrow \pm \infty} r(x) = 0.
    \end{array}
    \right.
  \end{equation}

The essential spectrum $\sigma_{ess} \{ {\rm (5.1)} \}$ of (5.1) has the following property:

\begin{pro} 
There exists a negative constant $C_{2}=C_{2}(\hat{\tau},\hat{\theta})$ independent of $\eps$ such that
  \begin{equation}
    \sigma_{ess} \{ {\rm (5.1)} \} \leq \eps^{2}  C_{2}< 0
  \end{equation}
holds for small $\eps>0$.
\end{pro}

\begin{proof}
The proof can be done in a similar way to that of Proposition 4.1.
\end{proof}

As in the previous section, we decompose the eigenvalue problem (5.1) into an equivalent pair of the problems on $I=(0, \infty )$ with the Neumann and Dirichlet boundary conditions at $x=0$. Both eigenvalue problems have the same boundary conditions at $x=\infty$.
  \begin{equation}
    \lim_{x \rightarrow \infty} p(x) = 0, \ \ 
    \lim_{x \rightarrow \infty} q(x) = 0, \ \ 
    \lim_{x \rightarrow \infty} r(x) = 0.
  \end{equation}

Let ${\rm (\widehat{EP})}_{ev}$ be the eigenvalue problem (5.1), (5.3) and 
  \begin{displaymath}
    p_{x}(0) = 0, \ \ q_{x}(0) = 0, \ \ r_{x}(0) = 0,
  \end{displaymath}
and ${\rm (\widehat{EP})}_{od}$ be the eigenvalue problem (5.1), (5.3) and 
  \begin{displaymath}
    p(0) = 0, \ \ q(0) = 0, \ \ r(0) = 0.
  \end{displaymath}
Applying the same procedures in section 4, (5.1) can be rewritten as follows.
  \begin{equation}
    \left \{
    \begin{array}{l}
       - q_{xx} + q + \dis\frac{\hat{\tau}}{\eps^{2}} \lambda q 
              - \eps \alpha (L^{\eps} - \lambda)^{\dag}(q) \\
              \\ \hspace{15mm}
       = \dis \frac{\langle \alpha q + \beta r, 
          \ \phi_{0}^{\eps}/\sqrt{\eps} \rangle}
          {\zeta_{0}^{\eps}/\eps^{2} - \lambda/\eps^{2}} 
          \phi_{0}^{\eps}/\sqrt{\eps} 
           + \eps \beta (L^{\eps} - \lambda)^{\dag}(r),  \\
               \\
       - D^{2} r_{xx} + r + \dis\frac{\hat{\theta}}{\eps^{2}} \lambda r 
              - \eps \beta (L^{\eps} - \lambda)^{\dag}(r) \\
               \\ \hspace{15mm}
       = \dis \frac{\langle \alpha q + \beta r, 
          \ \phi_{0}^{\eps}/\sqrt{\eps} \rangle}
          {\zeta_{0}^{\eps}/\eps^{2} - \lambda/\eps^{2}} 
          \phi_{0}^{\eps}/\sqrt{\eps} 
           + \eps \alpha (L^{\eps} - \lambda)^{\dag}(q).
    \end{array}
    \right.
  \end{equation}
Here again, as shown below, we classify the eigenvalues depending on their behaviors as $\eps \downarrow 0$. We only consider the first three cases, especially focus on the case (iii), and omit the remaining two cases, since the same arguments as before show that those are not relevant to the stability.

(i) $\lambda=O(1)$, 

(ii) $\lambda = o(1)$ with $\dis\lim_{\eps \downarrow 0} |\lambda|/\eps^{2}=\infty$,

(iii) $\lambda=O(\eps^{2})$,

(iv) $\dis\lim_{\eps \downarrow 0} |\lambda|/\eps^{2}=0$ as $\eps \downarrow 0$.

(v) $\dis\lim_{\eps \downarrow 0} |\lambda| = \infty$.

\vspace{2mm}

Case (i) $\lambda=O(1)$: When $\eps \downarrow 0$, the first term of the righthand side of (5.4) goes to zero so that we have 
  \begin{equation}
    \hat{\tau} \lambda q = 0 \ \ \ \mbox{and} \ \ \ 
    \hat{\theta} \lambda r = 0
  \end{equation}
as $\eps \downarrow 0$. First we consider the case where $\lambda$ is real. If $\lambda = 0$, and $q$ and $r$ are non-zero functions, (5.4) becomes an eigenvalue problem with $\lambda=0$. Then the eigenfunctions are the same appeared in the case (iii) of section 4. On the other hand, if $\lambda \neq 0$ and $q \equiv 0 \equiv r$, then $p \equiv 0$. They are not eigenfunctions.

Next, when $\lambda$ is complex, we set
  \begin{equation}
    \lambda = \lambda_{R} + i \lambda_{I} \ \ 
    (\lambda_{I} \neq 0), \ \ 
    q = q_{R} + i q_{I}, \ \ r = r_{R} + i r_{I}.
  \end{equation}
Substituting (5.6) into (5.5), we have
  \begin{displaymath}
    \left(
      \begin{array}{cc}
         \lambda_{R} & -\lambda_{I}  \\
               \\
         \lambda_{I} & \lambda_{R}
      \end{array}
    \right) 
    \left(
      \begin{array}{c}
         q_{R}  \\
               \\
         q_{I}
      \end{array}
    \right) = 
    \left(
      \begin{array}{c}
         0  \\
               \\
         0
      \end{array}
    \right), \hspace{5mm}
    \left(
      \begin{array}{cc}
         \lambda_{R} & -\lambda_{I}  \\
               \\
         \lambda_{I} & \lambda_{R}
      \end{array}
    \right) 
    \left(
      \begin{array}{c}
         r_{R}  \\
               \\
         r_{I}
      \end{array}
    \right) = 
    \left(
      \begin{array}{c}
         0  \\
               \\
         0
      \end{array}
    \right)
  \end{displaymath}
Since $\lambda_{I} \neq 0$, we obtain $(q_{R}, q_{I}) \equiv (0, 0)$ and $(r_{R}, r_{I}) \equiv (0, 0)$, that is, $q \equiv 0$ and $r \equiv 0$, and hence $p \equiv 0$, which implies that there are no complex eigenvalues.

\vspace{2mm}

Case (ii) $\lambda = o(1)$ with $\dis\lim_{\eps \downarrow 0} |\lambda|/\eps^{2}=\infty$: When $\eps \downarrow 0$, it is easy to see that (5.4) becomes
  \begin{displaymath}
    \hat{\tau} q = 0 \ \ \ \mbox{and} \ \ \ \hat{\theta} r = 0.
  \end{displaymath}
In this case, we obtain $q \equiv 0$ and $r \equiv 0$, hence $p \equiv 0$, hence this is not the eigenvalue.

Note that the discussion for the above two cases (i) and (ii) do not depend on the boundary condition at $x=0$. The following case (iii) is most interesting and the drift and Hopf instabilities come out through the analysis of the associated SLEP equation similar to (4.16).

Case (iii) $\lambda=O(\eps^{2})$ : The associated SLEP equation (5.7) similar to (4.16) can be obtained by applying the same computation to (5.4) as in section 4. 
Here we set the asymptotic form of $\lambda(\eps)$ as $\lambda(\eps) = \eps^{2} \hat{\lambda}(\eps)$ and $\hat{\lambda}^{*} = \lim_{\eps \downarrow 0} \hat{\lambda}(\eps)$. If $\hat{\lambda}^{*}$ does not exist, we choose a subsequence $\{ \eps_{k} \}$ with $\lim_{k \rightarrow \infty} \eps_{k} = 0$ such that $\lim_{k \rightarrow \infty} \hat{\lambda}(\eps_{k})$ exists and apply the same argument to it. The resulting equation reads
  \begin{equation}
    \hat{\zeta}_{0}^{*} - \hat{\lambda}^{*} = 4 (\kappa^{*})^{2} [
        \alpha \langle \widehat{K}_{q}^{*}(\delta_{*}), \ \delta_{*} \rangle
      + \beta  \langle \widehat{K}_{r}^{*}(\delta_{*}), \ \delta_{*} \rangle ],
  \end{equation}
where
  \begin{displaymath}
    \widehat{T}_{q}^{*} := - \dis\frac{d^{2}}{dx^{2}} 
              + 1 + \hat{\tau} \hat{\lambda}^{*}, \ \ \ 
    \widehat{T}_{r}^{*} := - D^{2} \dis\frac{d^{2}}{dx^{2}} 
             + 1 + \hat{\theta} \hat{\lambda}^{*},
  \end{displaymath}
  \begin{displaymath}
    \widehat{K}_{q}^{*} := (\widehat{T}_{q}^{*})^{-1}, \ \ \ 
    \widehat{K}_{r}^{*} := (\widehat{T}_{r}^{*})^{-1}, \ \ \ 
    \hat{\lambda}^{*} = \dis\lim_{\eps \downarrow 0} \lambda(\eps)/\eps^{2}.
  \end{displaymath}
We should note that $\widehat{K}_{q}^{*}$ (resp. $\widehat{K}_{r}^{*}$) depends strongly on $\hat{\lambda}^{*}$ and $\hat{\tau}$ (resp. $\hat{\theta}$). 

In the following, we consider (5.7) for the even symmetric case and the odd symmetric one separately. Let $\widehat{K}_{q}^{e,*}$ (resp. $\widehat{K}_{r}^{e,*}$) be $\widehat{K}_{q}^{*}$ (resp.  $\widehat{K}_{r}^{*}$) with Neumann boundary condition at $x=0$, and $\widehat{K}_{q}^{o,*}$ (resp. $\widehat{K}_{r}^{o,*}$) be $\widehat{K}_{q}^{*}$ (resp.  $\widehat{K}_{r}^{*}$) with Dirichlet boundary condition at $x=0$. Then we have the following lemma.

\begin{lem}
  \begin{displaymath}
    \langle \widehat{K}_{q}^{e,*}(\delta_{*}), \delta_{*} \rangle
      = x^{*} g_{+}(2 \omega_{q} x^{*}), \ \ \ 
    \langle \widehat{K}_{q}^{o,*}(\delta_{*}), \delta_{*} \rangle
      = x^{*} g_{-}(2 \omega_{q} x^{*}), 
  \end{displaymath}
  \begin{displaymath}
    \langle \widehat{K}_{r}^{e,*}(\delta_{*}), \delta_{*} \rangle
      = \dis\frac{x^{*}}{D^{2}} g_{+}(2 \omega_{r} x^{*}), \ \ \ 
    \langle \widehat{K}_{r}^{o,*}(\delta_{*}), \delta_{*} \rangle
      = \dis\frac{x^{*}}{D^{2}} g_{-}(2 \omega_{r} x^{*}),
  \end{displaymath}
where
  \begin{displaymath}
    g_{\pm}(y):=\dis\frac{1}{y} (1 \pm e^{-y}), \ \ \ 
    \omega_{q} := \sqrt{1 + \hat{\tau} \hat{\lambda}^{*}}, \ \ \ 
    \omega_{r} := \dis\frac{1}{D} \sqrt{1 + \hat{\theta} \hat{\lambda}^{*}}.
  \end{displaymath}
\end{lem}

\begin{proof} See Appendix B.
\end{proof}

\noindent
The next lemma is key to study the dependency on $\hat{\lambda}^{*}$, $\hat{\tau}$ and $\hat{\theta}$ of (5.7).

\begin{lem}
For $y>0$, $g_{\pm}(y)$ have the following properties:
  \begin{displaymath}
    \dis\frac{d g_{\pm}}{dy} < 0, \ \ \dis\frac{d^{2} g_{\pm}}{dy^{2}}>0, \ \ 
    \dis\frac{d g_{+}}{dy} = O(1/y^{2}) \ \ 
          (y \rightarrow +0, \ y \rightarrow \infty).
  \end{displaymath}
\end{lem}

\begin{proof}
We can prove them by computing directly, for example,
  \begin{displaymath}
    \dis\frac{d g_{+}}{dy} = - \dis\frac{1}{y^{2} e^{y}}(e^{y} + y + 1)<0, \ \ 
    \dis\frac{d^{2} g_{+}}{dy^{2}} = \dis\frac{1}{y^{3} e^{y}}
         (2 e^{y} + y^{2} + 2y + 2) > 0.
  \end{displaymath}
\end{proof}

\vspace{2mm} 


{\bf 5.1. Existence of the drift bifurcation for $\tau=O(1/\eps^{2})$ and $\theta=O(1/\eps^{2})$ : Proof of Theorem 2.8}

\vspace{2mm} 

First we consider ${\rm (\widehat{EP})}_{od}$. Then SLEP equation (5.7) is recast as
  \begin{equation}
    \begin{array}{rcl}
      G_{od}(\hat{\lambda}^{*}; \hat{\tau}, \hat{\theta}) &:=& 
       \hat{\lambda}^{*} - \hat{\zeta}_{0}^{*} + 4 (\kappa^{*})^{2} [
       \alpha \langle \widehat{K}_{q}^{o,*}(\delta_{*}), \ \delta_{*} \rangle
       + \beta  \langle \widehat{K}_{r}^{o,*}(\delta_{*}), \ 
                        \delta_{*} \rangle ] \\
                         \\
       & = &
       \hat{\lambda}^{*}  - \hat{\zeta}_{0}^{*}
       + 4 (\kappa^{*})^{2} \left[ \alpha x^{*} g_{-}(2 \omega_{q} x^{*})
       + \beta  \dis\frac{x^{*}}{D^{2}} g_{-}(2 \omega_{r} x^{*}) \right]  \\
                         \\
       & = & 0,
    \end{array}
  \end{equation}
where
  \begin{displaymath}
    \hat{\zeta}_{0}^{*} = 
       2 (\kappa^{*})^{2} \left[ 
       \alpha (1 - e^{-2 x^{*}}) + \dis\frac{\beta}{D} (1 - e^{-2 x^{*}/D})
       \right].
  \end{displaymath}

$G_{od}(\hat{\lambda}^{*}; \hat{\tau}, \hat{\theta})$ has the following properties.

\begin{lem}
For $\hat{\tau}>0$ and $\hat{\theta}>0$, $G_{od}(\hat{\lambda}^{*}; \hat{\tau}, \hat{\theta})$ satisfies

\vspace{2mm}

\noindent
{\rm (i)} \ \ \ $G_{od}(0; \hat{\tau}, \hat{\theta}) = 0$,

\vspace{2mm}

\noindent
{\rm (ii)} \ \ \ $\dis\frac{\pd}{\pd \hat{\lambda}^{*}} 
          G_{od}(0; \hat{\tau}, \hat{\theta}) = 1 +
          4 (\kappa^{*} x^{*})^{2} \left[ \alpha g'_{-}(2x^{*}) \hat{\tau}
          + \dis\frac{\beta}{D^{3}} g'_{-}(2x^{*}/D) \hat{\theta} \right]$,

\vspace{2mm}

\noindent
{\rm (iii)} \ \ \ $\dis\frac{\pd^{2}}{\pd \hat{\lambda}^{* 2}} 
          G_{od}(0; \hat{\tau}, \hat{\theta})>0$.
\end{lem}

\begin{proof}
We can show them by the direct calculations, so we omit the proof.
\end{proof}

{\it Proof of Theorem 2.8.} (i) is clear from Lemmas 5.3 and 5.4. In fact, the constants $C_{1}$ and $C_{2}$ are given by
  \begin{displaymath}
    C_{1} = - 4 (\kappa^{*} x^{*})^{2} \alpha g'_{-}(2x^{*})>0, \ \ \ 
    C_{2} = - \dis\frac{4 (\kappa^{*} x^{*})^{2} \beta}{D^{3}} 
              g'_{-}(2x^{*}/D)>0.
  \end{displaymath}

As for the second part (ii),  we expand $G_{od}(\hat{\lambda}^{*}; \hat{\tau}, \hat{\theta})$ into double Taylor series at $(\hat{\tau}_{0}, \hat{\theta}_{0})$, that is,
  \begin{displaymath}
    \begin{array}{l}
      G_{od}(\hat{\lambda}^{*}; \hat{\tau}, \hat{\theta}) 
      = A_{1} \hat{\lambda}^{* 2} 
        - A_{2} \hat{\lambda}^{*} (\hat{\tau} - \hat{\tau}_{0})
        - A_{3} \hat{\lambda}^{*} (\hat{\theta} - \hat{\theta}_{0}) + h.o.t. \\
                 \\ \hspace{20mm}
      = \hat{\lambda}^{*} \{ A_{1} \hat{\lambda}^{*} 
        - A_{2} (\hat{\tau} - \hat{\tau}_{0} )
        - A_{3} (\hat{\theta} - \hat{\theta}_{0}) \} + h.o.t.,
    \end{array}
  \end{displaymath}
where
  \begin{displaymath}
    A_{1} = \dis\frac{1}{2} \dis\frac{\pd^{2}}{\pd \hat{\lambda}^{* 2}} 
        G_{od}(0; \hat{\tau}_{0}, \hat{\theta}_{0})>0, \ \ \ 
    A_{2} = - \dis\frac{\pd^{2}}{\pd \hat{\lambda}^{*} \pd \hat{\tau}}
         G_{od}(0; \hat{\tau}_{0}, \hat{\theta}_{0})>0,
  \end{displaymath}
  \begin{displaymath}
    A_{3} = - \dis\frac{\pd^{2}}{\pd \hat{\lambda}^{*} \pd \hat{\theta}}
         G_{od}(0; \hat{\tau}_{0}, \hat{\theta}_{0})>0
  \end{displaymath}
and used the fact that
  \begin{displaymath}
    \dis\frac{\pd}{\pd \hat{\tau}} 
         G_{od} (0; \hat{\tau}_{0}, \hat{\theta}_{0}) = 
    \dis\frac{\pd^{2}}{\pd \hat{\tau}^{2}} 
         G_{od} (0; \hat{\tau}_{0}, \hat{\theta}_{0}) = 0, \ 
    \dis\frac{\pd}{\pd \hat{\theta}} 
         G_{od} (0; \hat{\tau}_{0}, \hat{\theta}_{0}) = 
    \dis\frac{\pd^{2}}{\pd \hat{\theta}^{2}} 
         G_{od} (0; \hat{\tau}_{0}, \hat{\theta}_{0}) = 0.
  \end{displaymath}
Then, we see that $\hat{\lambda}^{*}$ is expanded as
  \begin{displaymath}
    \hat{\lambda}^{*} \approx 
      \hat{\lambda}_{d}^{*}(\hat{\tau}, \hat{\theta}) := \dis\frac{1}{A_{1}} [
      A_{2} (\hat{\tau} - \hat{\tau}_{0}) 
        + A_{3} (\hat{\theta} - \hat{\theta}_{0})]
  \end{displaymath}
in the neighborhood of $(\hat{\tau}, \hat{\theta})=(\hat{\tau}_{0}, \hat{\theta}_{0})$. Therefore we have
  \begin{displaymath}
    \dis\frac{\pd}{\pd \hat{\tau}} 
          \hat{\lambda}_{d}^{*} (\hat{\tau}_{0}, \hat{\theta}_{0}) > 0, \ \ \ 
    \dis\frac{\pd}{\pd \hat{\theta}} 
          \hat{\lambda}_{d}^{*} (\hat{\tau}_{0}, \hat{\theta}_{0}) > 0
  \end{displaymath}
for $(\hat{\tau}_{0},\hat{\theta}_{0}) \in \Gamma_{d}$. Then, $\Gamma_{d}$ is a set of the drift bifurcation points in the parameter space $(\hat{\tau}, \hat{\theta})$.

The discussion for the algebraic multiplicity of the zero eigenvalue, see Appendix C. This completes the proof of Theorem 2.8. \hfill $\Box$

\vspace{2mm}

\noindent
{\bf 5.2. Existence of the Hopf bifurcation for $\tau=O(1/\eps^{2})$ and $\theta=O(1/\eps^{2})$ : Proofs of Theorem 2.9 and Proposition 2.11}

\vspace{2mm}

In this subsection, we will show the existence of the complex eigenvalues of 
${\rm (\widehat{EP})}_{ev}$ and their dependency on the parameters $\hat{\tau}$ and $\hat{\theta}$. The associated SLEP equation reads
  \begin{equation}
    G_{ev}(\hat{\lambda}^{*}; \hat{\tau}, \hat{\theta}) := 
       \hat{\lambda}^{*}  - \hat{\zeta}_{0}^{*}
       + 4 (\kappa^{*})^{2} \left[ \alpha x^{*} g_{+}(2 \omega_{q} x^{*})
       + \beta \dis\frac{x^{*}}{D^{2}} g_{+}(2 \omega_{r} x^{*}) \right] = 0
  \end{equation}
where
  \begin{displaymath}
    g_{+}(y):=\dis\frac{1}{y} (1 + e^{-y}), \ \ 
    \omega_{q} := \sqrt{1 + \hat{\tau} \hat{\lambda}^{*}}, \ \ 
    \omega_{r} := \dis\frac{1}{D} \sqrt{1 + \hat{\theta} \hat{\lambda}^{*}}.
  \end{displaymath}
It should be noted that $G_{ev}(\hat{\lambda}^{*}; \hat{\tau}, \hat{\theta})$ is analytic with respect to $\hat{\lambda}^{*}$, $\hat{\tau}$ and $\hat{\theta}$. 

\vspace{2mm}

{\it Proof of Theorem 2.9.} \ 
First we shall show the existence of a pure imaginary complex eigenvalue $\hat{\lambda^{*}} = i \xi$ ($\xi>0$). We prepare several preliminary things before analyzing the SLEP equation (5.9). First note that
  \[
    2 \omega_{q} x^{*} = 2 x^{*} \sqrt{1 + i \hat{\tau} \xi}, \ \ 
    2 \omega_{r} x^{*} = \dis\frac{2 x^{*}}{D} \sqrt{1 + i \hat{\theta} \xi}
  \]
are complex, we introduce the notation $X(z; d)$, $Y(z; d)$ and $z(\zeta)$ as
  \begin{displaymath}
    X(z; d) + i Y(z; d) := d \sqrt{1 + i \zeta},
  \end{displaymath}
  \begin{displaymath}
    z = z(\zeta) := {\rm arg} \sqrt{1 + i \zeta} = \frac{1}{2} \tan^{-1} \zeta
  \end{displaymath}
for $d>0$ and $\zeta \geq 0$. Then, $2 \omega_{q} x^{*}$ and $2 \omega_{r} x^{*}$ can be  represented as
  \begin{displaymath}
    2 \omega_{q} x^{*} 
     = X(z(\hat{\tau} \xi); 2 x^{*}) + i Y(z(\hat{\tau} \xi); 2 x^{*}),
  \end{displaymath}
  \begin{displaymath}
    2 \omega_{r} x^{*} 
     = X \left(z(\hat{\theta} \xi); \dis\frac{2 x^{*}}{D} \right) 
              + i Y \left(z(\hat{\theta} \xi); \dis\frac{2 x^{*}}{D} \right).
  \end{displaymath}
Moreover, noting that
  \begin{displaymath}
    |d \sqrt{1 + i \zeta}| = d(1 + \zeta^2)^{1/4} 
       = \dis\frac{d}{\sqrt{\cos 2z}},
  \end{displaymath}
we see that
  \begin{equation}
    X = X(z; d) = \dis\frac{d \cos z}{\sqrt{\cos 2z}}, \ \ \ 
    Y = Y(z; d) = \dis\frac{d \sin z}{\sqrt{\cos 2z}}.
  \end{equation}
Also we have
  \begin{equation}
    X(z; d)^{2} + Y(z; d)^{2} = \dis\frac{d^{2}}{\cos 2z}, \ \ 
    X(z; d)^{2} - Y(z; d)^{2} = d^{2},
  \end{equation}
$X \geq d$ and $Y \geq 0$ since $0 \leq z <\pi/4$.

Next we decompose $g_{+}(X + i Y)$ into the real and the imaginary parts. Noting that
  \begin{displaymath}
    \begin{array}{l}
       g_{+}(X + i Y) = \dis\frac{1 + \exp(-X - i Y)}{X + i Y} \\
             \\ \hspace{18mm}
       = \dis\frac{1}{X^{2} + Y^{2}} [X + e^{-X} (X \cos Y - Y \sin Y)] \\
             \\ \hspace{20mm}
           - i \dis\frac{1}{X^{2} + Y^{2}} [Y + e^{-X} (Y \cos Y + X \sin Y)],
    \end{array}
  \end{displaymath}
we define $R(z; d)$ and $I(z; d)$ by
  \begin{displaymath}
    \begin{array}{l}
      R(z;d):= \dis\frac{1}{d} 
             \sqrt{\cos 2z} \ [\cos z + e^{-X(z;d)} \cos(Y(z;d) + z) ], \\
                \\
      I(z;d):= - \dis\frac{1}{d} 
               \sqrt{\cos 2z} \ [\sin z + e^{-X(z;d)} \sin(Y(z;d) + z) ].
    \end{array}
  \end{displaymath}
Then we have 
  \begin{displaymath}
    g_{+}(2 x^{*} \sqrt{1 + i \hat{\tau} \xi}) 
      = R \left (
          \dis\frac{1}{2} \tan^{-1} \hat{\tau} \xi,2 x^{*} 
      \right)
      + i I \left (
           \dis\frac{1}{2} \tan^{-1} \hat{\tau} \xi,2 x^{*}
      \right)
  \end{displaymath}
and
  \begin{displaymath}
    g_{+} \left( \dis\frac{2 x^{*}}{D} \sqrt{1 + i \hat{\theta} \xi} \right)
      = R \left ( 
          \dis\frac{1}{2} \tan^{-1} \hat{\theta} \xi, \dis\frac{2 x^{*}}{D}
      \right )
       + i I \left ( 
          \dis\frac{1}{2} \tan^{-1} \hat{\theta} \xi, \dis\frac{2 x^{*}}{D}
      \right ),
  \end{displaymath}
Here we used the fact (5.10) and (5.11).

Now we can formulate the SLEP equation (5.9) into a pair of the real forms:
  \begin{equation}
      \hat{\zeta}_{0}^{*} = 4 (\kappa^{*})^{2} \left [
      \alpha x^{*} R 
        \left ( \dis\frac{1}{2} \tan^{-1} \hat{\tau} \xi,2 x^{*} \right)
      + \dis\frac{\beta x^{*}}{D^{2}} R 
      \left  ( 
          \dis\frac{1}{2} \tan^{-1} \hat{\theta} \xi, \dis\frac{2 x^{*}}{D}
        \right )
      \right ]
  \end{equation}
  \begin{equation}
      \xi = - 4 (\kappa^{*})^{2} \left [
      \alpha x^{*} I 
        \left ( \dis\frac{1}{2} \tan^{-1} \hat{\tau} \xi,2 x^{*} \right)
      + \dis\frac{\beta x^{*}}{D^{2}} I 
        \left  ( 
           \dis\frac{1}{2} \tan^{-1} \hat{\theta} \xi, \dis\frac{2 x^{*}}{D}
        \right )
      \right ].
  \end{equation}
In the following, we try to find $(\hat{\tau}, \hat{\theta})$ and $\xi>0$ satisfying (5.12) and (5.13). Here we prepare several key lemmas for the proof. Especially, the monotonicity of $R(z;d)$ with respect to $z$ is important. The following two lemmas can be obtained by direct computations.

\begin{lem}
$X(z;d)$ and $Y(z;d)$ $(0 \leq z < \pi/4)$ have the following properties.
  \begin{displaymath}
    \frac{dX}{dz} = \dis\frac{d \sin z}{\cos 2z \sqrt{\cos 2z}} \geq 0, \ \ \ 
    \frac{dY}{dz} = \dis\frac{d \cos z}{\cos 2z \sqrt{\cos 2z}} > 0, \ \ \ 
  \end{displaymath}
  \begin{displaymath}
    \lim_{z \rightarrow \pi/4 - 0} X(z) = \infty, \ \ 
    \lim_{z \rightarrow \pi/4 - 0} Y(z) = \infty
  \end{displaymath}
\end{lem}

\begin{lem}
  \begin{displaymath}
      \sin z = \dis\frac{Y}{\sqrt{2 Y^{2} + d^{2}}}, \ \ 
      \cos z = \sqrt{\dis\frac{Y^{2} + d^{2}}{2 Y^{2} + d^{2}}},
  \end{displaymath}
  \begin{displaymath}
      \cos 2z = \dis\frac{d^{2}}{2 Y^{2} + d^{2}}, \ \
      \sin 2z = \dis\frac{2y \sqrt{Y^{2} + d^{2}}}{2 Y^{2} + d^{2}}
  \end{displaymath}
\end{lem}

\begin{lem}
{\rm (i)} For $0 \leq z < \pi/4$ and $d>0$,
  \begin{displaymath}
    R(z; d)>0, \ \ \ I(z; d) \leq 0.
  \end{displaymath}
It holds that $I(z; d)=0$ if and only if $z=0$.

\vspace{2mm}

{\rm (ii)} For $0 \leq z < \pi/4$ and $d>0$,
  \begin{displaymath}
    \dis\frac{dR}{dz}(z;d) \leq 0.
  \end{displaymath}
It holds that $(dR/dz)(z;d) = 0$ if and only if $z=0$.

\vspace{2mm}

{\rm (iii)}
  \begin{displaymath}
    \lim_{z \rightarrow \pi/4 - 0} R(z;d) = 0, \ \ \ 
    \lim_{z \rightarrow \pi/4 - 0} I(z;d) = 0.
  \end{displaymath}
\end{lem}

\begin{proof}
See Appendix D.
\end{proof}

As was mentioned in section 2, we introduce the polar coordinate $(s, \psi)$ in the parameter space $(\hat{\tau}, \hat{\theta})$:
  \begin{displaymath}
    \hat{\tau} = s \cos \psi, \ \ \ \hat{\theta} = s \sin \psi, \ \ \ 
    s>0, \ \ \ 0 < \psi < \dis\frac{\pi}{2}.
  \end{displaymath}
Then the SLEP equation (5.9) becomes
  \begin{equation}
    \widehat{G}_{ev}(\hat{\lambda}^{*}; s, \psi)
       := G_{ev}(\hat{\lambda}^{*}(s,\psi); s \cos \psi, s \sin \psi) = 0.
  \end{equation}
Also (5.12) and (5.13) are recast as
  \begin{displaymath}
    \hat{\zeta}_{0}^{*} = 4 (\kappa^{*})^{2} \left [
      \alpha x^{*} R 
         \left ( \dis\frac{1}{2} \tan^{-1} (s \xi \cos \psi), 2 x^{*} \right)
      + \dis\frac{\beta x^{*}}{D^{2}} R 
        \left  ( 
          \dis\frac{1}{2} \tan^{-1} (s \xi \sin \psi), \dis\frac{2 x^{*}}{D}
      \right )
      \right ]
  \end{displaymath}
and
  \begin{equation}
    \begin{array}{l}
    \xi = - 4 (\kappa^{*})^{2} \left [
      \alpha x^{*} I 
         \left ( \dis\frac{1}{2} \tan^{-1} (s \xi \cos \psi),2 x^{*} \right)
      \right.   \\
                \\ \hspace{40mm}
      \left. + \dis\frac{\beta x^{*}}{D^{2}} I 
        \left  ( 
           \dis\frac{1}{2} \tan^{-1} (s \xi \sin \psi), \dis\frac{2 x^{*}}{D}
       \right )
      \right ],
    \end{array}
  \end{equation}
respectively. Set $\eta = s \xi$ and introduce a new function $\hat{R}(\eta, \psi)$:
  \begin{displaymath}
    \begin{array}{l}
      \hat{R}(\eta, \psi) := 4 (\kappa^{*})^{2} \left [
      \alpha x^{*} R 
        \left ( \dis\frac{1}{2} \tan^{-1} (\eta \cos \psi), 2 x^{*} \right)
        \right.   \\
                  \\ \hspace{40mm}
      \left. + \dis\frac{\beta x^{*}}{D^{2}} R 
      \left  ( 
          \dis\frac{1}{2} \tan^{-1} (\eta \sin \psi), \dis\frac{2 x^{*}}{D}
      \right )
      \right ].
    \end{array}
  \end{displaymath}
For any $\psi \in (0, \pi/2)$, we see from Lemma 5.7 that
  \begin{displaymath}
    \hat{R}(0, \psi) = 2 (\kappa^{*})^{2} \left[
       \alpha (1+e^{-2 x^{*}}) + \dis\frac{\beta}{D} (1+e^{-2 x^{*}/D})
       \right] > \hat{\zeta}_{0}^{*},
  \end{displaymath}
  \begin{displaymath}
    \lim_{\eta \rightarrow \infty} \hat{R}(\eta, \psi) = 0, \ \ \
    \dis\frac{d \hat{R}}{d \eta} < 0 \ \ (\eta>0).
  \end{displaymath}
Therefore, by the intermediate value theorem, we conclude that there exists a unique positive $\eta = \eta^{*} (\psi)$ such that
  \begin{displaymath}
    \hat{\zeta}_{0}^{*} = \hat{R}(\eta^{*} (\psi), \psi).
  \end{displaymath}
Now we set $s \xi = \eta^{*} (\psi)$ and substitute this into (5.15). Then we obtain the expression for $\xi$ as
  \begin{displaymath}
    \begin{array}{l}
    \xi = \xi^{*} (\psi) := - 4 (\kappa^{*})^{2} \left [
    \alpha x^{*} I 
       \left ( \dis\frac{1}{2} \tan^{-1} (\eta^{*} (\psi) \cos \psi), 
                              2 x^{*} \right) \right.  \\
                 \\ \hspace{30mm}
    \left. + \dis\frac{\beta x^{*}}{D^{2}} I
      \left  ( 
          \dis\frac{1}{2} \tan^{-1} (\eta^{*} (\psi) \sin \psi), 
                              \dis\frac{2 x^{*}}{D}
      \right )
    \right ].
    \end{array}
  \end{displaymath}
Finally, substituting $\xi^{*} (\psi)$ into $s \xi = \eta^{*} (\psi)$, we have 
$s^{*}(\psi) := \eta^{*} (\psi)/\xi^{*} (\psi)$ and the principal part of a pure imaginary eigenvalue is given by
  \begin{displaymath}
     \hat{\lambda}_{h}^{*} = i \xi^{*}(\psi) \ \ \mbox{at} \ \ 
     (\hat{\tau}, \hat{\theta}) 
         = (s^{*}(\psi) \cos \psi, s^{*}(\psi) \sin \psi).
  \end{displaymath}
We have completed the proof of the existence to pure imaginary eigenvalues. 
\vspace{2mm}

Suppose that there exists a complex solution $\hat{\lambda}^{*}$ of (5.9) for some $(\hat{\tau}, \hat{\theta})$. We shall show that such a solution is isolated in ${\bf C}$ and uniquely parameterized by $(\hat{\tau}, \hat{\theta})$ locally.

\begin{pro} 
Let $(\hat{\lambda}_{0}^{*}, \hat{\tau}_{0}, \hat{\theta}_{0}) \in {\bf C} \times {\bf R}_{+} \times {\bf R}_{+}$ be a solution of $G_{ev}(\hat{\lambda}^{*}; \hat{\tau}, \hat{\theta}) =0$. If ${\rm Im} \ \hat{\lambda}_{0}^{*} \neq 0$, then
  \begin{displaymath}
    \dis\frac{d}{d \hat{\lambda}^{*}} 
       G_{ev}(\hat{\lambda}_{0}^{*}; \hat{\tau}_{0}, \hat{\theta}_{0}) \neq 0.
  \end{displaymath}
\end{pro}
This proposition implies that such a solution can be extended at least locally, and depends (real-) analytically on $(\hat{\tau}, \hat{\theta})$ as $\hat{\lambda}^{*} = \hat{\lambda}^{*}(\hat{\tau}, \hat{\theta})$ with $\hat{\lambda}^{*}(\hat{\tau}_{0}, \hat{\theta}_{0}) = \hat{\lambda}_{0}^{*}$. 
Therefore it is also possible to extend the pure imaginary solution 
$i \xi^{*}(\psi)$ of (5.9) (respectively, $-i \xi^{*}(\psi)$) to a locally isolated complex solution $\hat{\lambda}_{h}^{*}(s,\psi)$ of (5.9) (respectively, $\overline{\hat{\lambda}_{h}^{*}(s,\psi)}$) that depends smoothly on $s$ near $s=s^{*}(\psi)$ with $\hat{\lambda}_{h}^{*}(s^{*}(\psi),\psi) = i \xi^{*}(\psi)$ (respectively, $\overline{\hat{\lambda}_{h}^{*}(s^{*}(\psi),\psi)} = -i \xi^{*}(\psi)$).

\vspace{2mm}

{\it Proof of Proposition 5.8}. 
First we note that
  \begin{equation}
    \dis\frac{d}{d \hat{\lambda}^{*}} 
            G_{ev}(\hat{\lambda}^{*}; \hat{\tau}, \hat{\theta}) = 
      1 + 4 (\kappa^{*})^{2} x^{*} 
      \dis\frac{d}{d \hat{\lambda}^{*}} \left[ \alpha g_{+}(2 x^{*} \omega_{q})
          + \dis\frac{\beta}{D^{2}} g_{+}(2 x^{*} \omega_{r}) \right].
  \end{equation}
We shall prove that the imaginary part of (5.16) is not equal to zero. We shall examine the imaginary part of
  \begin{displaymath}
    \dis\frac{d}{d \lambda} g_{+}(d \sqrt{1 + c \lambda})
  \end{displaymath}
for $d>0$ and $c>0$. Noting that
  \begin{displaymath}
    \dis\frac{d}{d \lambda} g_{+}(d\sqrt{1+c \lambda}) = 
    g'_{+}(d\sqrt{1+c \lambda}) \times \dis\frac{c d}{2 \sqrt{1 + c \lambda}},
        \ \ \ 
    g'_{+}(y) = - \dis\frac{1 + e^{-y}}{y^{2}} - \dis\frac{e^{-y}}{y},
  \end{displaymath}
we have
  \begin{equation}
    \begin{array}{rcl}
      \dis\frac{d}{d \lambda} g_{+}(d \sqrt{1 + c \lambda}) &=& 
        \left \{
            - \dis\frac{1 + e^{-d \sqrt{1 + c \lambda}}}
                 {d^{2} (1 + c \lambda)}
      - \dis\frac{e^{-d \sqrt{1 + c \lambda}}}{d \sqrt{1 + c \lambda}}
        \right \}
      \times \dis\frac{c d}{2 \sqrt{1 + c \lambda}}  \\
                       \\
      &=& - \dis\frac{(1 + e^{-d \sqrt{1 + c \lambda}}) c}
                 {2d (1 + c \lambda)^{3/2}}
      - \dis\frac{e^{-d \sqrt{1 + c \lambda}} c}{2 (1 + c \lambda)} \\
                       \\
      &=& - \dis\frac{1 + e^{-d \sqrt{\rho + i \zeta}} }
                             {2d (\rho + i \zeta)^{3/2}} 
          - \dis\frac{e^{-d \sqrt{\rho + i \zeta}} }{2 (\rho + i \zeta)},
    \end{array}
  \end{equation}
where $\lambda = \lambda_{R} + i \lambda_{I}$, $\rho := 1 + c \lambda_{R}$, 
$\zeta := c \lambda_{I}$. For $d>0$, $\rho>0$, $\zeta>0$ and $c>0$, we define $x$, $y$ and $z$ as
  \begin{displaymath}
    x + i y = d \sqrt{\rho + i \zeta}, \ \ \ 
    z = \dis\frac{1}{2} \tan^{-1} \dis\frac{\zeta}{\rho}.
  \end{displaymath}
Then we see that
  \begin{displaymath}
    x = \dis\frac{d \sqrt{\rho} \cos z}{\sqrt{\cos 2z}}, \ \ \ 
    y = \dis\frac{d \sqrt{\rho} \sin z}{\sqrt{\cos 2z}}, \ \ \ 
    x^{2} - y^{2} = d^{2} \rho.
  \end{displaymath}
Each term of (5.17) can be computed as
  \begin{displaymath}
    \begin{array}{l}
      \dis\frac{1 + e^{-d \sqrt{\rho + i \zeta}} }{2d (\rho + i \zeta)^{3/2}}
       = \dis\frac{\cos 3z}{2d (\rho^{2} + \zeta^{2})^{3/4}}
       + \dis\frac{ e^{-x}}{2d (\rho^{2} + \zeta^{2})^{3/4}}
         ( \cos 3z \cos y - \sin 3z \sin y )  \\
                    \\ \hspace{20mm}
        + i \left[ - \dis\frac{\sin 3z}{2d (\rho^{2} + \zeta^{2})^{3/4}}
        - \dis\frac{ e^{-x}}{2d (\rho^{2} + \zeta^{2})^{3/4}}
          ( \cos y \sin 3z + \cos 3z \sin y ) \right] \\
    \end{array}
  \end{displaymath}
  \begin{displaymath}
    \begin{array}{l}
       = \dis\frac{\cos 3z}{2d (\rho^{2} + \zeta^{2})^{3/4}}
       + \dis\frac{ e^{-x}}{2d (\rho + i \zeta)^{3/4}} \cos (3z + y) \\
                    \\ \hspace{20mm}
        - i \left[ \dis\frac{\sin 3z}{2d (\rho^{2} + \zeta^{2})^{3/4}}
           + \dis\frac{ e^{-x}}{2d (\rho^{2} + \zeta^{2})^{3/4}} \sin (3z + y)
       \right]
    \end{array}
  \end{displaymath}
  \begin{displaymath}
    \begin{array}{rcl}
      \dis\frac{e^{-d \sqrt{\rho + i \zeta}} }{2 (\rho + i \zeta)}
       &=& \dis\frac{ e^{-x}}{2 (\rho^{2} + \zeta^{2})}
           (\rho \cos y - \zeta \sin y)
           - i \dis\frac{ e^{-x}}{2 (\rho^{2} + \zeta^{2})}
           (\zeta \cos y + \rho \sin y)  \\
                    \\
       &=& \dis\frac{ e^{-x}}{2 (\rho^{2} + \zeta^{2})^{1/2}}
           (\cos 2z \cos y - \sin2 z \sin y)  \\
                   \\
       & & - i \dis\frac{ e^{-x}}{2 (\rho^{2} + \zeta^{2})^{1/2}}
           (\sin2 z \cos y + \cos 2z \sin y)  \\
                    \\
       &=& \dis\frac{ e^{-x}}{2 (\rho^{2} + \zeta^{2})^{1/2}} \cos (2z + y) 
         - i \dis\frac{ e^{-x}}{2 (\rho^{2} + \zeta^{2})^{1/2}} \sin (2z + y).
    \end{array}
  \end{displaymath}
Therefore we have
  \begin{displaymath}
    \begin{array}{l}
      - \dis\frac{1}{c} \dis\frac{d}{d \lambda} g_{+}(d \sqrt{1 + c \lambda}) =
      \dis\frac{1 + e^{-d \sqrt{1 + c \lambda}}}
                 {2d (1 + c \lambda)^{3/2}} 
      + \dis\frac{e^{-d \sqrt{1 + c \lambda}}}{2 (1 + c \lambda)}  \\
                   \\
      = \dis\frac{\cos 3z}{2d (\rho^{2} + \zeta^{2})^{3/4}}
       + \dis\frac{ e^{-x}}{2d (\rho^{2} + \zeta^{2})^{3/4}} \cos (3z + y)
       + \dis\frac{ e^{-x}}{2 (\rho^{2} + \zeta^{2})^{1/2}} \cos (2z + y)  \\
                   \\ \hspace{5mm}
      - i \left[
        \dis\frac{\sin 3z}{2d (\rho^{2} + \zeta^{2})^{3/4}}
           + \dis\frac{ e^{-x}}{2d (\rho^{2} + \zeta^{2})^{3/4}} \sin (3z + y)
           + \dis\frac{ e^{-x}}{2 (\rho^{2} + \zeta^{2})^{1/2}} \sin (2z + y)
      \right].
    \end{array}
  \end{displaymath}
Now the imaginary part of (5.17) becomes
  \begin{equation}
    \begin{array}{l}
      \dis\frac{2d (\rho^{2} + \zeta^{2})^{3/4}}{c} \cdot 
      {\rm Im} \left \{ \dis\frac{d}{d \lambda} g_{+}(d \sqrt{1 + c \lambda}) 
      \right \}   \\
                  \\ \hspace{15mm}
      = \sin 3z + \dis\frac{1}{e^{x}} \sin (3z + y)
        + \dis\frac{d \sqrt{\rho}}{e^{x} \sqrt{\cos 2z}} \sin (2z + y).
    \end{array}
  \end{equation}
The following relations are analogous to Lemma 5.6. 

\begin{lem}
Let $x$ and $y$ be
  \begin{displaymath}
    x = \dis\frac{d \sqrt{\rho} \cos z}{\sqrt{\cos 2z}}, \ \ 
    y = \dis\frac{d \sqrt{\rho} \sin z}{\sqrt{\cos 2z}}
  \end{displaymath}
for $z \in (0, \pi/4)$, $d>0$ and $\rho>0$, then
  \begin{displaymath}
      \sin z = \dis\frac{y}{\sqrt{2 y^{2} + d^{2} \rho}}, \ \ 
      \cos z = \sqrt{\dis\frac{y^{2} + d^{2} \rho}{2 y^{2} + d^{2} \rho}},
  \end{displaymath}
  \begin{displaymath}
      \cos 2z = \dis\frac{d^{2} \rho}{2 y^{2} + d^{2} \rho}, \ \
      \sin 2z = \dis\frac{2y \sqrt{y^{2} + d^{2} \rho}}{2 y^{2} + d^{2} \rho}.
  \end{displaymath}
\end{lem}

By using the results of Lemma 5.9, we can show that the right hand side of (5.18) is positive in the same manner as in the proof of $R(z;d)>0$ in Lemma 5.7 (see (D.1) in Appendix D). This yields that
  \begin{displaymath}
    {\rm Im} \left \{ \dis\frac{d}{d \lambda} g_{+}(d \sqrt{1 + c \lambda}) 
      \right \} \neq 0
  \end{displaymath}
 \hfill $\Box$

\vspace{2mm}

Finally we shall prove the transversality:
  \begin{displaymath}
    \dis\frac{\pd}{\pd s} {\rm Re} \hat{\lambda}^{*}(s^{*}(\psi),\psi) \neq 0.
  \end{displaymath}
Noting that
  \begin{displaymath}
    \begin{array}{l}
      \dis\frac{d}{d \lambda} g_{+}(d\sqrt{1+c \lambda}) = 
             g'_{+}(d\sqrt{1+c \lambda}) 
             \cdot \dis\frac{c d}{2 \sqrt{1 + c \lambda}}
        = g'_{+} \cdot \dis\frac{d}{2 \sqrt{1 + c \lambda}} \cdot s
        \left \{
           \begin{array}{l}
              \cos \psi  \\
              \sin \psi
           \end{array}
        \right.,
                        \\
                        \\
      \dis\frac{d}{d s} g_{+}(d\sqrt{1+c \lambda}) = 
        g'_{+}(d\sqrt{1+c \lambda}) 
             \cdot \dis\frac{d \lambda}{2 \sqrt{1 + c \lambda}}
             \cdot \dis\frac{\pd c}{\pd s}
        = g'_{+} \cdot \dis\frac{d}{2 \sqrt{1 + c \lambda}} 
               \cdot \lambda
        \left \{
           \begin{array}{l}
              \cos \psi  \\
              \sin \psi
           \end{array}
        \right.,
    \end{array}
  \end{displaymath}
we see that
  \begin{displaymath}
    \dis\frac{d}{d s} g_{+}(d\sqrt{1+c \lambda}) 
       = \dis\frac{\lambda}{s}
         \dis\frac{d}{d \lambda} g_{+}(d\sqrt{1+c \lambda}),
  \end{displaymath}
where $c=\hat{\tau}$ or $c=\hat{\theta}$. Differentiating the SLEP equation (5.14) with respect to $s$, we have
  \begin{equation}
    \begin{array}{l}
      \dis\frac{d}{d s} 
        G_{ev}(\hat{\lambda}^{*}(s^{*}(\psi),\psi);s \cos \psi, s \sin \psi) \\
                \\ \hspace{18mm}
        = \dis\frac{\pd \hat{\lambda}^{*}}{\pd s} 
          + 4 (\kappa^{*})^{2} x^{*} \alpha 
              \dis\frac{d}{d \hat{\lambda}^{*}} g_{+}(2 x^{*} \omega_{q})
              \left ( \dis\frac{d \hat{\lambda}^{*}}{d s}  
              + \dis\frac{\hat{\lambda}^{*}}{s} \right )  \\
                   \\ \hspace{25mm}
      + \dis\frac{4 (\kappa^{*})^{2} x^{*} \beta}{D^{2}} 
            \dis\frac{d}{d \hat{\lambda}^{*}} g_{+}(2 x^{*} \omega_{r})
              \left ( \dis\frac{d \hat{\lambda}^{*}}{d s}  
              + \dis\frac{\hat{\lambda}^{*}}{s} \right ) = 0.
    \end{array}
  \end{equation}
In the following, we set $s=s^{*}(\psi)$ on (5.19). Then
  \begin{displaymath}
    \hat{\lambda}^{*}(s^{*}(\psi),\psi) = i \hat{\lambda}_{I}^{*}
        (= i \xi^{*}(\psi))
  \end{displaymath}
for some $\hat{\lambda}_{I}^{*} >0$ since ${\rm Re} \hat{\lambda}^{*}(s^{*}(\psi),\psi)=0$. For notational simplicity, we set
  \begin{displaymath}
    \begin{array}{l}
      E + i F := 
      \dis\frac{\pd \hat{\lambda}^{*}}{\pd s} = 
        \dis\frac{\pd}{\pd s} {\rm Re} \hat{\lambda}^{*}(s^{*}(\psi),\psi) + 
       i \dis\frac{\pd}{\pd s} {\rm Im} \hat{\lambda}^{*}(s^{*}(\psi),\psi), \\
                    \\
      G_{q} + i H_{q} := 4 (\kappa^{*})^{2} x^{*} \alpha 
        \dis\frac{d}{d \hat{\lambda}^{*}} g_{+}(2 x^{*} \omega_{q}), \ \ \ 
      G_{r} + i H_{r} := 4 (\kappa^{*})^{2} x^{*} \beta 
        \dis\frac{d}{d \hat{\lambda}^{*}} g_{+}(2 x^{*} \omega_{r}).
    \end{array}
  \end{displaymath}
When $s=s^{*}(\psi)$, (5.19) is written as
  \begin{displaymath}
    \begin{array}{l}
      E + i F + \{ (G_{q} + i H_{q}) + (G_{r} + i H_{r}) \}
        \left( E + i F + \dis\frac{ i \hat{\lambda}_{I}^{*} }{s^{*}(\psi)} 
        \right) = 0.
    \end{array}
  \end{displaymath}
This is equivalent to
  \begin{displaymath}
    \left(
    \begin{array}{cc}
       1 + G_{q} + G_{r}  & -(H_{q} + H_{r})  \\
                    \\
       H_{q} + H_{r}      & 1 + G_{q} + G_{r}
    \end{array}
    \right)
    \left(
    \begin{array}{c}
       E  \\
          \\
       F
    \end{array}
    \right) = \dis\frac{\hat{\lambda}_{I}^{*}}{s^{*}(\psi)}
    \left(
    \begin{array}{c}
       H_{q} + H_{r}  \\
              \\
       -(G_{q} + G_{r})
    \end{array}
    \right).
  \end{displaymath}
Therefore
  \begin{displaymath}
    \left(
    \begin{array}{c}
       E  \\
          \\
       F
    \end{array}
    \right) 
    = \dis\frac{\hat{\lambda}_{I}^{*}}
    { s^{*}(\psi) \{ (1 + G_{q} + G_{r})^{2} + (H_{q} + H_{r})^{2} \} } \times
  \end{displaymath}
  \begin{displaymath}
    \hspace{35mm} \left(
    \begin{array}{cc}
       1 + G_{q} + G_{r}  & H_{q} + H_{r}  \\
                    \\
       -(H_{q} + H_{r})    & 1 + G_{q} + G_{r}
    \end{array}
    \right)
    \left(
    \begin{array}{c}
       H_{q} + H_{r}  \\
              \\
       -(G_{q} + G_{r})
    \end{array}
    \right)
  \end{displaymath}
and we obtain
  \begin{displaymath}
    E = \dis\frac{\pd}{\pd s} {\rm Re} \ 
             \hat{\lambda}^{*}(s^{*}(\psi), \psi)
      = \dis\frac{\hat{\lambda}_{I}^{*} (H_{q} + H_{r})}
       { s^{*}(\psi) \{ (1 + G_{q} + G_{r})^{2} + (H_{q} + H_{r})^{2} \} }.
  \end{displaymath}
We have already proved $H_{q}>0$ and $H_{r}>0$ in the proof of Proposition 5.8 (i.e., (5.18) is positive). Hence we conclude $E>0$. This completes the proof of Theorem 2.9. 
\hfill $\Box$

\vspace{3mm}

In the rest of this subsection, we prove Proposition 2.11 which reveals the local behavior of real eigenvalues. We first trace the minimum point of the SLEP function $\widehat{G}_{ev}(\lambda; s, \psi)$ with respect to real $\lambda$ for fixed $\psi \in (0, \pi/2)$. In the following, for the notational simplicity, we denote $\hat{\lambda}^{*}$ by $\lambda$ and $G_{ev}(\lambda; s \cos \psi, s \sin \psi)$ by $G(\lambda, s)$:
  \begin{displaymath}
    \begin{array}{l}
      G(\lambda; s) := G_{ev}(\lambda; s \cos \psi, s \sin \psi)  \\
                      \\ \hspace{5mm}
        = \lambda - \hat{\zeta}_{0}^{*}
        + 4 x^{*} (\kappa^{*})^{2} \left[ 
           \alpha g_{+} (2 x^{*} \sqrt{1 + \lambda s \cos \psi} )
        + \dis\frac{\beta}{D^{2}} g_{+} 
           \left( \dis\frac{2 x^{*}}{D} \sqrt{1 + \lambda s \sin \psi} \right)
      \right ].
    \end{array}
  \end{displaymath}
Note that the domain of $G(\lambda; s)$ is 
  \begin{displaymath}
    \{ \lambda \in {\bf R} \ | \ \lambda > 
    - \min \{ 1/s \cos \psi, 1/s \sin \psi \} \},
  \end{displaymath}
which depends on $s$ and $\psi$. In order to show that there exists a local absolute minimum of $G(\lambda; s)$ for each $s>0$, the next lemma is useful.

\begin{lem}
{\rm (i)} 
  \begin{displaymath}
    G(0;s) = 4(\kappa^{*})^{2} \left( \alpha e^{-2 x^{*}} + \dis\frac{\beta}{D} e^{-2 x^{*}/D} \right) > 0, \ \ \ \dis\frac{\pd^{2} G}{\pd \lambda^{2}}(\lambda, s) > 0
  \end{displaymath}
for any $s>0$.

{\rm (ii)} For any $s>0$, there exists $\lambda=\underline{\lambda}(s)$ such that
  \begin{displaymath}
    \dis\frac{\pd G}{\pd \lambda}(\underline{\lambda}(s); s) = 0.
  \end{displaymath}

{\rm (iii)}
  \begin{displaymath}
    \underline{\lambda}(s) \left \{
    \begin{array}{l}
       < 0 \hspace{7mm} 0 < s < s_{c}  \\
       = 0 \hspace{7mm} s = s_{c}    \\
       > 0 \hspace{7mm} s_{c} < s.
    \end{array}
    \right.
  \end{displaymath}
where, $s_{c}$ is given by
  \begin{displaymath}
    s_{c} := -1 \Big / 4 (x^{*} \kappa^{*})^{2} \left[ 
           \alpha g'_{+} (2 x^{*}) \cos \psi
         + \dis\frac{\beta}{D^{3}} g'_{+} 
           \left( \dis\frac{2 x^{*}}{D} \right) \sin \psi \right ] > 0.
  \end{displaymath}

{\rm (iv)}
  \begin{displaymath}
    1 + s A \underline{\lambda}(s) = O(s^{2/3}), \ \ 
    \underline{\lambda}(s) = O \left( \dis\frac{1}{s} \right) \ \ 
    (s \rightarrow +0),
  \end{displaymath}
  \begin{displaymath}
    \underline{\lambda}(s) = O \left( \dis\frac{1}{s^{1/3}} \right) \ \ 
    (s \rightarrow \infty),
  \end{displaymath}
where $A=\cos \psi$ or $A=\sin \psi$.
\end{lem}

\begin{proof}
(i) This is a consequence of direct computation as shown below.
  \begin{displaymath}
    \begin{array}{l}
    \dis\frac{\pd^{2} G}{\pd \lambda^{2}}(\lambda;s) =
      4 \alpha (x^{*} \kappa^{*} \cos \psi )^{2} \left[ 
      - \dis\frac{1}{2(1 + \lambda s \cos \psi)^{3/2}} 
        \ g'_{+} (2 x^{*} \sqrt{1 + \lambda s \cos \psi} ) \right. \\
                    \\ \hspace{15mm}
      \left. + \dis\frac{x^{*}}{1 + \lambda s \cos \psi}
        \ g''_{+} (2 x^{*} \sqrt{1 + \lambda s \cos \psi} ) \right]  \\
                    \\ \hspace{15mm}
      + \dis\frac{4 \beta (x^{*} \kappa^{*} \sin \psi )^{2}}{D^{3}} \left[ 
      - \dis\frac{1}{2(1 + \lambda s \sin \psi)^{3/2}} 
        \ g'_{+} \left( \dis\frac{2 x^{*}}{D} \sqrt{1 + \lambda s \sin \psi} 
                   \right) \right.  \\
                    \\ \hspace{15mm}
      + \left. \dis\frac{x^{*}}{D(1 + \lambda s \sin \psi)}
        \ g''_{+} \left( \dis\frac{2 x^{*}}{D} \sqrt{1 + \lambda s \sin \psi} 
                   \right) \right]  \\
                    \\ \hspace{5mm}
      > 0.
    \end{array}
  \end{displaymath}

(ii) First note that the derivative of $G(\lambda;s)$ with respect to $\lambda$ is directly computed as 
  \begin{displaymath}
    \begin{array}{l}
      \dis\frac{\pd G}{\pd \lambda}(\lambda; s)
        = 1 + 4 (x^{*} \kappa^{*})^{2} \left[ 
           \dis\frac{\alpha \cos \psi}{\sqrt{1 + \lambda s \cos \psi}}
           \ g'_{+} (2 x^{*} \sqrt{1 + \lambda s \cos \psi} ) \right. \\
                           \\  \hspace{20mm} \left.
          + \dis\frac{\beta \sin \psi}{D^{3} \sqrt{1 + \lambda s \sin \psi}} 
          \ g'_{+} \left( \dis\frac{2 x^{*}}{D} \sqrt{1 + \lambda s \sin \psi} 
                   \right)
      \right ] s.
    \end{array}
  \end{displaymath}
We consider the asymptotic behavior of $G_{\lambda}(\lambda, s)$ as $\lambda \rightarrow$ $- \min \{1/s \cos \psi,$ $1/s \sin \psi\} + 0$ for fixed $s>0$. Without loss of generality, we focus on the case when $0 < \psi < \pi/4$ and $\lambda \rightarrow - 1/s \cos \psi + 0$.
Since
  \begin{displaymath}
    \dis\frac{g'_{+} (c \sqrt{1+\lambda s \cos \psi})}
        {\sqrt{1 + \lambda s \cos \psi}} 
        = O \left( \dis\frac{1}{(1+\lambda s \cos \psi)^{3/2}} \right)
  \end{displaymath}
for $c>0$ and $g'_{+}(y)<0$, we have
  \begin{displaymath}
    \lim_{\lambda \rightarrow -1/s \cos \psi + 0} G_{\lambda}(\lambda;s) 
        = - \infty.
  \end{displaymath}
For the case when $\lambda \rightarrow +\infty$, the fact
  \begin{displaymath}
    \dis\frac{g'_{+}(y)}{y} = O \left( \dis\frac{1}{y^{3}} \right) \ \ \ 
    (y \rightarrow \infty)
  \end{displaymath}
yields
  \begin{displaymath}
    \lim_{\lambda \rightarrow \infty} G_{\lambda}(\lambda;s) = 1 > 0
  \end{displaymath}
Noting that $G_{\lambda \lambda}(\lambda; s)>0$, we see that there exists a unique $\lambda=\underline{\lambda}(s)$ such that 
$G_{\lambda}(\underline{\lambda}(s);s)=0$.

\vspace{2mm}

(iii) Note that
  \begin{displaymath}
    G_{\lambda}(0; s) = 1 + 4 (x^{*} \kappa^{*})^{2} \left[ 
           \alpha g'_{+} (2 x^{*}) \cos \psi
         + \dis\frac{\beta}{D^{3}} g'_{+} 
           \left( \dis\frac{2 x^{*}}{D} \right) \sin \psi \right ] s.
  \end{displaymath}
Since the coefficient of the second term is negative, $G_{\lambda}(0; s) = 0$ has unique zero $s_{c}>0$. Then we can see that
  \[
    G_{\lambda}(0; s) \left \{
    \begin{array}{l}
       > 0 \hspace{7mm} 0 < s < s_{c}  \\
       = 0 \hspace{7mm} s = s_{c}    \\
       < 0 \hspace{7mm} s_{c} < s.
    \end{array}
    \right.
  \]
Combining the above facts and $G_{\lambda \lambda}(\lambda, s)>0$, we obtain conclusions (iii).

\vspace{2mm}

(iv) \ First we consider the asymptotic behavior when $s \rightarrow +0$. 
Since there exists $\lambda$ such that $G_{\lambda}(\lambda;s)=0$, the following asymptotic results hold 
  \begin{displaymath}
    g'_{+} (c_{1} \sqrt{1 + \lambda s \cos \psi}) 
         \dis\frac{s}{\sqrt{1 + \lambda s \cos \psi}} = O(1),
  \end{displaymath}
  \begin{displaymath}
    g'_{+} (c_{1} \sqrt{1 + \lambda s \sin \psi}) 
         \dis\frac{s}{\sqrt{1 + \lambda s \sin \psi}} = O(1)
  \end{displaymath}
for $c_{1}>0$ as $s \rightarrow +0$. Noting that
  \begin{displaymath}
    g'_{+}(y) = - \dis\frac{e^{y} + y + 1}{y^{2} e^{y}},
  \end{displaymath}
we see that $1 + s c_{2} \lambda$ ($c_{2}>0$) doesn't diverge but converges to $+0$ when $s \rightarrow +0$. That is,
  \begin{displaymath}
    \dis\frac{s}{(1 + s c_{2} \lambda) \sqrt{1 + s c_{2} \lambda}} 
       = O(1) \ \ \ (s \rightarrow +0).
  \end{displaymath}
In order to know the rate of decay, we set
  \begin{displaymath}
    1 + s c_{2} \lambda \approx a s^{\delta}
  \end{displaymath}
for $a$ and equating the power of $s$, we have $\delta=2/3$. Therefore
  \begin{displaymath}
    \lambda = O(s^{-1/3}) - \dis\frac{1}{s} = O \left( \dis\frac{1}{s} \right)
    \ \ \ (s \rightarrow +0).
  \end{displaymath}

When $s \rightarrow \infty$, we see that
  \begin{displaymath}
    \dis\frac{s}{(1+s c_{2} \lambda) \sqrt{1+s c_{2} \lambda}}. 
       = O(1)
  \end{displaymath}
In this case, $\lambda$ has to converge to $0$ so that we set
  \begin{displaymath}
    \lambda \approx a s^{\delta}
  \end{displaymath}
and equating the power of $s$, we have $\delta=-1/3$. Therefore we conclude that
  \begin{displaymath}
    \lambda \approx a s^{-1/3} = O \left( \dis\frac{1}{s^{1/3}} \right)
    \ \ \ (s \rightarrow \infty)
  \end{displaymath}
\end{proof}

Since $G_{\lambda \lambda}(\lambda;s)>0$, $G(\lambda;s)$ has a local absolute minimum at $\lambda=\underline{\lambda}(s)$. Let us define the minimum value as
  \begin{displaymath}
    m(s) := G(\underline{\lambda}(s);s)
  \end{displaymath}
and we study its dependency on $s$ precisely. First it holds that 

\begin{lem}
{\rm (i)}
  \begin{displaymath}
    \dis\frac{dm}{d s} \ \left \{
    \begin{array}{l}
       > 0 \hspace{7mm} 0 < s < s_{c}  \\
       = 0 \hspace{7mm} s = s_{c}    \\
       < 0 \hspace{7mm} s_{c} < s
    \end{array}
    \right.
  \end{displaymath}
{\rm (ii)}
  \begin{displaymath}
    \lim_{s \rightarrow +0} m(s) = -\infty, \ \ \ 
    m(s_{c}) = G(0; s_{c})>0, \ \ \ 
    \lim_{s \rightarrow \infty} m(s) = - \hat{\zeta}_{0}^{*} < 0
  \end{displaymath}
\end{lem}

\begin{proof}
(i) By the direct computation, we obtain
  \begin{displaymath}
    \begin{array}{l}
      \dis\frac{dm}{d s} 
         = \dis\frac{d}{d s} G(\underline{\lambda}(s); s)
         = G_{\lambda}(\underline{\lambda}(s); s)
             \dis\frac{d \underline{\lambda}}{d s}
             + G_{s}(\underline{\lambda}(s); s)
         = G_{s}(\underline{\lambda}(s); s),
    \end{array}
  \end{displaymath}
where
 \begin{displaymath}
   \begin{array}{l}
     G_{s}(\lambda; s) = 4 (x^{*} \kappa^{*})^{2} \left[ 
           \alpha g'_{+} (2 x^{*} \sqrt{1 + \lambda s \cos \psi} )
           \dis\frac{\cos \psi}{\sqrt{1 + \lambda s \cos \psi}}
       \right.  \\
                  \\ \hspace{25mm}
       \left. 
         + \dis\frac{\beta}{D^{3}} g'_{+} 
           \left( \dis\frac{2 x^{*}}{D} \sqrt{1 + \lambda s \sin \psi} \right)
           \dis\frac{\sin \psi}{\sqrt{1 + \lambda s \sin \psi}}
      \right ] \lambda.
    \end{array}
  \end{displaymath}
Noting that $g_{+}'<0$ and $0<\psi<\pi/4$, we can see that the coefficient of $\lambda$ is negative. 
This means that the sign of $dm/ds$ is opposite to that of $\underline{\lambda}(s)$. Combining the results of Lemma 5.10, we conclude (i).

\vspace{2mm}

(ii) \ $m(s)$ is computed as
  \begin{equation}
    \begin{array}{l}
      m(s) = G(\underline{\lambda}(s); s)  \\
                   \\ \hspace{8mm}
      = \underline{\lambda}(s) - \hat{\zeta}_{0}^{*}
      + 4 x^{*} (\kappa^{*})^{2} \left[ 
      \alpha g_{+} (2 x^{*} \sqrt{1 + \underline{\lambda}(s) s \cos \psi} ) 
      \right. \\
                   \\ \hspace{15mm}
      \left. + \dis\frac{\beta}{D^{2}} g_{+} 
           \left( \dis\frac{2 x^{*}}{D} \sqrt{1 
              + \underline{\lambda}(s) s \sin \psi} \right)
      \right ].
    \end{array}
  \end{equation}
From Lemma 5.10, $\underline{\lambda}(s)<0$ and $O(1/s)$ ($s \rightarrow +0$). On the other hand, in view of Lemma 5.10, the third term of (5.20) is negative and of $O(1/s^{1/3})$ ($s \rightarrow +0$). So we conclude that
  \begin{displaymath}
    \lim_{s \rightarrow +0} m(s) = 
    \lim_{s \rightarrow +0} G(\underline{\lambda}(s),s) = -\infty.
  \end{displaymath}

On the other hand, when $s \rightarrow \infty$, we have
  \begin{displaymath}
    \begin{array}{l}
    m(s) 
      = O \left( \dis\frac{1}{s^{1/3}} \right) - \hat{\zeta}_{0}^{*}
      + 4 x^{*} (\kappa^{*})^{2} \left[ 
         \alpha g_{+} \left(2 x^{*} \sqrt{1+O(s^{2/3})} \right) \right. \\
                     \\ \hspace{10mm}
      \left. + \dis\frac{\beta}{D^{2}} g_{+} 
           \left( \dis\frac{2 x^{*}}{D} \sqrt{1 + O(s^{2/3})} \right)
    \right ]
    \end{array}
  \end{displaymath}
since $\underline{\lambda}(s) = O(1/s^{1/3} )$ from Lemma 5.10. So we obtain
  \begin{displaymath}
    \lim_{s \rightarrow \infty} m(s) = 
    \lim_{s \rightarrow \infty} G(\underline{\lambda}(s),s) 
      = -\hat{\zeta}_{0}^{*}.
  \end{displaymath}
\end{proof}

Now we are ready to prove Proposition 2.11. 

\vspace{2mm}

{\it Proof of Proposition 2.11.} \ First note that $G(0;s) > 0$ and the graph of $G(\lambda;s)$ is convex. When $0 < s < s_{c}$, $m(s) < 0$ and this means that the graph of $G(\lambda;s)$ crosses the negative region of $\lambda$-axis at two points. As $s$ increases, $m(s)$ increases and becomes $m(s) = 0$ at some $s=\underline{s}(\psi)$ for the first time. After that, $m(s)>0$ for a while, until $m(s)$ becomes zero at some $s=\overline{s}(\psi) (> \underline{s}(\psi))$ again. When $s > \overline{s}(\psi)$, the graph of $G(\lambda;s)$ crosses the positive region of $\lambda$-axis at two points.

The asymptotic forms of $\hat{\lambda}^{*}$ near $s=\underline{s}(\psi)$ (or $\overline{s}(\psi)$) can be obtained by expanding $\widehat{G}_{ev}(\lambda;s,\psi)$ into double Taylor series at $(\hat{\lambda}^{*,-},\underline{s}(\psi))$ (or $(\hat{\lambda}^{*,+},\overline{s}(\psi))$) for fixed $\psi \in (0,\pi/2)$. 

\hfill $\Box$

\section{Concluding remarks and outlook}

In this paper, we have discussed the existence and stability properties of the standing pulse solutions to the three-component reaction diffusion system (1.1) by using the matched asymptotic expansion method (MAE) and the SLEP method. MAE combined with SLEP is one of the powerful analytical tools and plays a complementary role to the geometric singular perturbation method (GSP) for 1D problem. Here we make a list of advantages of it: firstly it allows us to describe a precise behavior of critical eigenvalues, as shown in section 5, especially for the complex ones, so that the parametric dependency of bifurcation points can be clarified in detail. In fact, we presented how a pair of complex critical eigenvalues emerge from the real axis, cross the imaginary axis, and back to the real axis again as a parameter varies. It should be noted that no additional conditions for the parameters $\hat{\tau}$, $\hat{\theta}$ and $D$ are necessary to study the trace of critical eigenvalues unlike \cite{vHDK}. This enables us to search for the singularities in much broader parameter space both analytically and numerically, which leads to finding a codimension two singularity of drift-Hopf type as in Theorem 2.9 and Propositions 2.11-2.12. Secondly MAE-SLEP approach shows a nice correspondence between the power of $\eps$ of a critical eigenvalue and how the associated eigenfunction behaves as $\eps \downarrow 0$. Such an eigenfunction in general does not remain as a usual function as $\eps \downarrow 0$ and an appropriate scaling is necessary to characterize it like Dirac's $\delta$-function. This scaling is directly linked to the order of critical eigenvalues and allows us to obtain the well-defined SLEP equation in the limit of $\eps \downarrow 0$ (see, for instance, (5.7)). Finally MAE-SLEP has a great potential to be extended to higher dimensional case as was shown in \cite{NS1} and \cite{NS2}. Recall that one of the necessities for a class of three-component systems is that it supports the coexistence of stable localized traveling patterns in higher dimensional space so that MAE-SLEP approach is one of the promising methods to study the existence and stabilities of those patterns. We close this section by presenting an outlook for future works.

\vspace{2mm}

\noindent
(i) {\it Unfolding the codimension two bifurcation points.} \ 
We have shown the existence of the codimension two points of drift-Hopf type by solving the SLEP equations as in Fig.3.  The next step is to unfold those singularities and show the existence of various types of moving objects including standing and traveling breathers. The collision dynamics among those emerging patterns is quite interesting, in particular for traveling breathers, because it depends not only on the velocity but the difference of phases. There are very few works in this direction \cite{TUN}, \cite{WIN2}, \cite{WIN1} so that this will open a new fertile ground in the class of three-component reaction diffusion systems. 


\vspace{2mm}

\noindent
(ii) {\it Localized moving patterns in higher dimensional space} \ 
Strong interactions among traveling spots present a variety of dynamics including annihilation, coalescence, and splitting when they collide each other \cite{NTU}, \cite{NTU3}. As a first step, existence and stability of traveling spots are necessary and there are already a couple of results in this direction \cite{vHS1}, \cite{vHS2}. To understand more detailed process of complex dynamics, one idea is to focus on the saddle solutions of high codimension, what is called the "scattors" introduced in \cite{NTU}, \cite{NTU3}. These saddle solutions arise during the large deformation of collision process and become a key to identify the transition path of whole collision process. Once we succeed to find them and examine their stabilities, then we can understand the onset of deformation encountered at a collision by looking at the profiles of unstable eigenfunctions of their saddles. As was mentioned above, one advantage of the MAE-SLEP method is its independency of the space dimension so that once the structure of internal transition layer is explicitly known, then it is basically possible to apply the SLEP method to investigate the stability. For the two-component reaction diffusion systems, there are some works \cite{NS1}, \cite{NS2} in this direction, however it is still a challenge to extend it to the three-component systems and to clarify the transition path of collision process via the dynamics around scattors.


\vspace{2mm}

\noindent
(iii) {\it Relation between the SLEP equation and the Evans function}\
The SLEP equation is closely related to the Evans function arising in GSP approach. In fact, for the stability analysis of the front solution in two component reaction diffusion systems, it was shown that the principal parts of both functions are equivalent up to a constant multiple (Ikeda, Nishiura and Suzuki \cite{INS}). This means that both methods work in a complementary way at least for 1D problem. It is a challenge to extend a geometrical method like GSP to a higher dimensional space and to get an insight complementary to analytical approach.

\vspace{4mm}

{\bf Acknowledgements}
This work was partially supported by JSPS KAKENHI Grant Number JP20K20341 (Y.N.).

\vspace{4mm}

\appendix
\section{Proof of Lemma 4.8} 

We only prove (4.17) since (4.18) is obtained similarly. For the symmetric mode, we consider the following problem:
  \begin{displaymath}
    \left \{
    \begin{array}{l}
       -D^{2} r_{xx} + r = R(x), \ \ x \in I, \\
                \\
       r_{x}(0) = 0, \ \ \dis\lim_{x \rightarrow \infty} r(x) = 0.
    \end{array}
    \right.
  \end{displaymath}
Let $r_{-}(x)$ and $r_{+}(x)$ be solutions of the equation $-D^{2} r_{xx} + r = 0$ with the boundary conditions
  \begin{displaymath}
    r_{-}(0)=1, \ \ \ \dis\frac{d}{dx}r_{-}(0)=0,
             \hspace{5mm}
    r_{+}(0)=1, \ \ \ \dis\lim_{x \rightarrow \infty} r_{+}(x)=0,
  \end{displaymath}
respectively. They are written as
  \begin{displaymath}
    r_{-}(x) = \cosh (x/D), \ \ \ r_{+}(x) = e^{-x/D}.
  \end{displaymath}
Then the Wronskian is $W(r_{-},r_{+})=-1/D$, and the green function $G(x,y)$ is given by
  \begin{displaymath}
    G(x,y)=\dis\frac{1}{D} \times
    \left \{
    \begin{array}{l}
       \cosh (x/D) e^{-y/D}, \hspace{5mm} 0 \leq x \leq y, \\
            \\
       \cosh (y/D) e^{-x/D}, \hspace{5mm} 0 \leq y \leq x.
    \end{array}
    \right.
  \end{displaymath}
So we have
  \begin{displaymath}
    \langle K_{r}^{e,*}(\delta_{*}), \ \delta_{*} \rangle 
       = G(x^{*}, x^{*}) = \dis\frac{\cosh (x^{*}/D)}{D e^{x^{*}/D}}.
  \end{displaymath}
$\langle K_{q}^{e,*}(\delta_{*}), \ \delta_{*} \rangle$ is obtained when we set $D=1$ at the above relation.  \hfill $\Box$

\section{Proof of Lemma 5.2}

We consider the case $\widehat{K}_{r}^{e,*}$, that is $\widehat{K}_{r}^{*}$ with Neumann boundary condition at $x=0$. The green function for $\hat{K}_{r}^{e,*}$ is given by
  \begin{displaymath}
    G(x,y)=\dis\frac{1}{\omega_{r} D^{2}} \times
    \left \{
    \begin{array}{l}
       \cosh (\omega_{r} x) e^{-\omega_{r} y} \hspace{5mm} 0 \leq x \leq y \\
            \\
       \cosh (\omega_{r} y) e^{-\omega_{r} x} \hspace{5mm} 0 \leq y \leq x,
    \end{array}
    \right.
  \end{displaymath}
where
  \begin{displaymath}
    \omega_{r} := \dis\frac{1}{D} \sqrt{1 + \hat{\theta} \hat{\lambda}^{*}}
  \end{displaymath}
and ${\rm Re} \ \hat{\lambda}^{*} > -1/\hat{\theta}$. Therefore we have
  \begin{displaymath}
    \langle \widehat{K}_{r}^{e,*}(\delta_{*}), \delta_{*} \rangle = G(x^{*}, x^{*})
     = \dis\frac{1}{\omega_{r} D^{2}} 
            \cosh (\omega_{r} x^{*}) e^{-\omega_{r} x^{*}}
     = \dis\frac{x^{*}}{D^{2}} \cdot \frac{1}{2 \omega_{r} x^{*}} 
          (1 +  e^{-2\omega_{r} x^{*}}).
  \end{displaymath}
Similarly we obtain
  \begin{displaymath}
    \langle \widehat{K}_{q}^{e,*}(\delta_{*}), \delta_{*} \rangle 
     = x^{*} \cdot \frac{1}{2 \omega_{q} x^{*}} 
          (1 +  e^{-2\omega_{q} x^{*}}), \ \ \ 
    \omega_{q} := \sqrt{1 + \hat{\tau} \hat{\lambda}^{*}},
  \end{displaymath}
where ${\rm Re} \ \hat{\lambda}^{*} >-1/\hat{\tau}$.

\vspace{3mm}

On the other hand, the green function for $\hat{K}_{r}^{o,*}$ is given by
  \begin{displaymath}
    G(x,y)=\dis\frac{1}{\omega_{r} D^{2}} \times
    \left \{
    \begin{array}{l}
       \sinh (\omega_{r} x) e^{-\omega_{r} y} \hspace{5mm} 0 \leq x \leq y \\
            \\
       \sinh (\omega_{r} y) e^{-\omega_{r} x} \hspace{5mm} 0 \leq y \leq x.
    \end{array}
    \right.
  \end{displaymath}
So we have
  \begin{displaymath}
    \langle \widehat{K}_{r}^{o,*}(\delta_{*}), \delta_{*} \rangle = G(x^{*}, x^{*})
     = \dis\frac{1}{\omega_{r} D^{2}} 
            \sinh (\omega_{r} x^{*}) e^{-\omega_{r} x^{*}}
     = \dis\frac{x^{*}}{D^{2}} \cdot \frac{1}{2 \omega_{r} x^{*}} 
          (1 -  e^{-2\omega_{r} x^{*}})
  \end{displaymath}
and
  \begin{displaymath}
    \langle \widehat{K}_{q}^{o,*}(\delta_{*}), \delta_{*} \rangle 
     = x^{*} \cdot \frac{1}{2 \omega_{q} x^{*}} 
          (1 -  e^{-2\omega_{q} x^{*}}).
  \end{displaymath}
\hfill $\Box$
\section{Algebraic multiplicity of the zero eigenvalue}

We discuss about the relation between the algebraic multiplicity of the zero eigenvalue at the drift bifurcation point and the degeneracy of zero solution of the SLEP equation. As was mentioned in Theorem 2.8, the drift bifurcation line $\Gamma_{d}$ is characterized by the common zero of the SLEP equation $G_{od}(\hat{\lambda}^{*}; \hat{\tau}, \hat{\theta}) = 0$ and its derivative  $(G_{od})_{\hat{\lambda}^{*}}(0; \hat{\tau}, \hat{\theta}) = 0$. 
Also the tangency at this point is equal to two from Lemma 5.4 (iii). We will show that the algebraic multiplicity of zero eigenvalue for the linearized problem (5.1) is exactly equal to the order of tangency of the SLEP equation at this point, namely two for our case.

 First recall that the condition that $(\hat{\tau}_{0},\hat{\theta}_{0})$ is a drift bifurcation point is equivalent to saying that there exist nonzero vectors $\varPhi$ and $\varPsi$ such that
  \begin{displaymath}
    {\mathcal L} \varPhi = {\bf 0}, \ \ \ 
    {\mathcal L} \varPsi = \widehat{T}_{0} \varPhi.\ \ \ 
  \end{displaymath}
As we mentioned in section 4, $\varPhi$ is a constant multiple of the derivative of the standing pulse solution. On the other hand, the existence of $\varPsi$ is equivalent to the solvability condition
  \begin{equation}
    \langle \widehat{T}_{0} \varPhi, \varPhi^{\#} \rangle = 0,
  \end{equation}
where $\varPhi^{\#}$ is a kernel function of the adjoint operator ${\mathcal L}^{\#}$ of ${\mathcal L}$. Therefore, we construct $\varPhi$ and $\varPhi^{\#}$ by the SLEP method and show (C.1). In fact, we can see that (C.1) is equivalent to
  \begin{displaymath}
    \dis\frac{\pd}{\pd \hat{\lambda}^{*}} 
        G_{od}(0; \hat{\tau}_{0}, \hat{\theta}_{0}) = 0.
  \end{displaymath}
For the later use, we set $\varPhi = (p, q, r)$ and $\varPhi^{\#}=(p^{\#}, q^{\#}, r^{\#})$. Then they satisfy
  \begin{equation}
      \left \{
      \begin{array}{l}
         L^{\eps} p - \eps \alpha q - \eps \beta r = 0 \\
                \\
          p + M q = 0  \\
                \\
          p + N r = 0
      \end{array}
    \right .
  \end{equation}
  \begin{equation}
    \left \{
    \begin{array}{l}
       L^{\eps} p^{\#} + q^{\#} + r^{\#} = 0 \\
              \\
       - \eps \alpha p^{\#} + M q^{\#} = 0  \\
              \\
       - \eps \beta p^{\#} + N r^{\#} = 0
    \end{array}
    \right .
  \end{equation}
We only prove that the principal part of the left hand side of (C.1) is equal to zero because the implicit function theorem guarantees that (C.1) holds for small $\eps>0$.

First, we solve (C.2) by the SLEP method. Note that the derivative of the stationary pulse solution is a candidate of the solution to (C.2). The first equation of (C.2) can be solved as
  \begin{displaymath}
    \begin{array}{rcl}
      p &=& (L^{\eps})^{-1} (\eps \alpha q + \eps \beta r)    \\
                \\
        &=& \dis \frac{\langle \alpha q + \beta r, 
                       \ \phi_{0}^{\eps}/\sqrt{\eps} \rangle}
              {\zeta_{0}^{\eps}/\eps^{2}} \phi_{0}^{\eps}/\sqrt{\eps}
            + \eps (L^{\eps})^{\dag}(\alpha q + \beta r).
    \end{array}
  \end{displaymath}
Substituting this into the second and the third equation of (C.2), we obtain
  \begin{equation}
    \left \{
      \begin{array}{l}
         - M q - \eps \alpha (L^{\eps})^{\dag} (q) = 
         \dis \frac{\langle \alpha q + \beta r, 
              \ \phi_{0}^{\eps}/\sqrt{\eps} \rangle}
              {\zeta_{0}^{\eps}/\eps^{2}} \phi_{0}^{\eps} /\sqrt{\eps}
          + \eps \beta (L^{\eps})^{\dag} (r), \\
                \\
         - N r - \eps \beta (L^{\eps})^{\dag} (r) = 
         \dis \frac{\langle \alpha q + \beta r, 
               \ \phi_{0}^{\eps}/\sqrt{\eps} \rangle}
              {\zeta_{0}^{\eps}/\eps^{2}} \phi_{0}^{\eps} /\sqrt{\eps}
          + \eps \alpha (L^{\eps})^{\dag} (q).
      \end{array}
    \right .
  \end{equation}
By the uniformity of the convergence of (C.4) as $\eps \downarrow 0$, we have only to consider the limiting systems. 

For the odd symmetric case, we take Dirichlet boundary condition at $x=0$. Now we define two operators $\widetilde{K}_{q}^{\eps}$ and $\widetilde{K}_{r}^{\eps}$ by
  \begin{displaymath}
    \widetilde{K}_{q}^{o,\eps} := 
    \left \{ - M - \eps \alpha (L^{\eps})^{\dag}(\cdot) \right \}^{-1}, \ \ 
    \widetilde{K}_{r}^{o,\eps} := 
    \left \{ - N - \eps \beta (L^{\eps})^{\dag}(\cdot)  \right \}^{-1}
  \end{displaymath}
with Dirichlet boundary conditions at $x=0$. Solving the limiting system of (C.4) with respect to $(q,r)$, we have
  \begin{equation}
    \left \{
    \begin{array}{rcl}
       q^{*} &=& 
      4 (\kappa^{*})^{2} \alpha \dis \frac{\langle q^{*}, \ \delta_{*} \rangle}
          {\hat{\zeta}_{0}^{*}} \widetilde{K}_{q}^{o,*}(\delta_{*}) + 
       4 (\kappa^{*})^{2} \beta \dis \frac{\langle r^{*}, \ \delta_{*} \rangle}
          {\hat{\zeta}_{0}^{*}} \widetilde{K}_{q}^{o,*}(\delta_{*}),    \\
               \\
       r^{*} &=& 
      4 (\kappa^{*})^{2} \alpha \dis \frac{\langle q^{*}, \ \delta_{*} \rangle}
          {\hat{\zeta}_{0}^{*}} \widetilde{K}_{r}^{o,*}(\delta_{*}) + 
       4 (\kappa^{*})^{2} \beta \dis \frac{\langle r^{*}, \ \delta_{*} \rangle}
          {\hat{\zeta}_{0}^{*}} \widetilde{K}_{r}^{o,*}(\delta_{*}),
    \end{array}
    \right.
  \end{equation}
where $\widetilde{K}_{q}^{o,*} = (-M)^{-1}$ and $\widetilde{K}_{r}^{o,*} = (-N)^{-1}$. Here we set 
  \begin{equation}
    q^{*} = A \widetilde{K}_{q}^{o,*}(\delta_{*}), \ \ \ r^{*} 
          = B \widetilde{K}_{r}^{o,*}(\delta_{*})
  \end{equation}
and substitute (C.6) into (C.5). Then we have
  \begin{displaymath}
    \left \{
    \begin{array}{l}
       A \widetilde{K}_{q}^{o,*}(\delta_{*}) =
       4 (\kappa^{*})^{2} A \alpha 
       \dis \frac{\langle \widetilde{K}_{q}^{o,*}(\delta_{*}), \ \delta_{*} \rangle}
          {\hat{\zeta}_{0}^{*}} \widetilde{K}_{q}^{*}(\delta_{*}) + 
       4 (\kappa^{*})^{2} B \beta 
       \dis \frac{\langle \widetilde{K}_{r}^{o,*}(\delta_{*}), \ \delta_{*} \rangle}
          {\hat{\zeta}_{0}^{*}} \widetilde{K}_{q}^{*}(\delta_{*})    \\
               \\
       B \widetilde{K}_{r}^{o,*}(\delta_{*}) =
       4 (\kappa^{*})^{2} A \alpha 
       \dis \frac{\langle \widetilde{K}_{q}^{o,*}(\delta_{*}), \ \delta_{*} \rangle}
          {\hat{\zeta}_{0}^{*}} \widetilde{K}_{r}^{*}(\delta_{*}) + 
       4 (\kappa^{*})^{2} B \beta 
       \dis \frac{\langle \widetilde{K}_{r}^{o,*}(\delta_{*}), \ \delta_{*} \rangle}
          {\hat{\zeta}_{0}^{*}} \widetilde{K}_{r}^{*}(\delta_{*})
    \end{array}
    \right.
  \end{displaymath}
or the equivalent system
  \begin{equation}
    \left(
    \begin{array}{cc}
      \hat{\zeta}_{0}^{*}
      - 4 (\kappa^{*})^{2} \alpha \langle \widetilde{K}_{q}^{o,*}(\delta_{*}), 
                  \ \delta_{*} \rangle &
      - 4 (\kappa^{*})^{2} \beta \langle \widetilde{K}_{r}^{o,*}(\delta_{*}), 
                  \ \delta_{*} \rangle \\
                \\
      - 4 (\kappa^{*})^{2} \alpha \langle \widetilde{K}_{q}^{o,*}(\delta_{*}), 
                  \ \delta_{*} \rangle &
       \hat{\zeta}_{0}^{*}
       - 4 (\kappa^{*})^{2} \beta \langle \widetilde{K}_{r}^{o,*}(\delta_{*}), 
                  \ \delta_{*} \rangle
    \end{array}
    \right)
    \left(
    \begin{array}{c}
       A \\
         \\
       B
    \end{array}
    \right) = 
    \left(
    \begin{array}{c}
       0 \\
         \\
       0
    \end{array}
    \right).
  \end{equation}
We can easily see that the determinant of the matrix in (C.7) is
  \begin{displaymath}
    \hat{\zeta}_{0}^{*} [ \hat{\zeta}_{0}^{*} - 4 (\kappa^{*})^{2} \{
        \alpha \langle \widetilde{K}_{q}^{o,*}(\delta_{*})
             \ \delta_{*} \rangle
      + \beta \langle \widetilde{K}_{r}^{o,*}(\delta_{*})
             \ \delta_{*} \rangle \} ]
  \end{displaymath}
Noting that $\widetilde{K}_{q}^{o,*} = K_{q}^{o,*}$, $\widetilde{K}_{r}^{o,*} = K_{r}^{o,*}$ and the formula of $\hat{\zeta}_{0}^{*}$ (see (4.19)), we find the determinant is equal to zero.

On the other hand, for the symmetric case, we obtain the similar system as (C.7) which has the Neumann boundary condition at $x=0$. Then the associated determinant
  \begin{displaymath}
    \hat{\zeta}_{0}^{*} [ \hat{\zeta}_{0}^{*} - 4 (\kappa^{*})^{2} \{
      \alpha \langle \widetilde{K}_{q}^{e,*}(\delta_{*}), \ \delta_{*} \rangle
    + \beta \langle \widetilde{K}_{r}^{e,*}(\delta_{*}), 
            \ \delta_{*} \rangle \} ]
  \end{displaymath}
is not zero, hence $A=B=0$. This means that $q^{*}$ and $r^{*}$ become trivial solutions. 
Therefore, we can conclude that the eigenfunction $\varPhi$ is an odd function with the dimension of the associated eigenspace is one. If we choose $A = B = 1$, we obtain
  \begin{displaymath}
    q^{*} = \widetilde{K}_{q}^{o,*}(\delta_{*}), \ \ \ 
    r^{*} = \widetilde{K}_{r}^{o,*}(\delta_{*}), \ \ \ 
    p^{*} = \dis \frac{4 (\kappa^{*})^{2} \langle \alpha q^{*} + \beta r^{*}, 
                \ \delta_{*} \rangle}
          {\hat{\zeta}_{0}^{*}} \delta_{*}
    = \delta_{*}.
  \end{displaymath}
This means that the asymptotic form $\varPhi^{*}$ of $\varPhi$ as $\eps \downarrow 0$ is given by
  \begin{displaymath}
    \varPhi^{*} = (p^{*}, q^{*}, r^{*})^{T}
     = (\delta_{*}, \ \widetilde{K}_{q}^{o,*}(\delta_{*}), \ 
                 \widetilde{K}_{r}^{o,*}(\delta_{*}))^{T}
  \end{displaymath}
Also the principal part of $\varPhi$ is represented as
  \begin{displaymath}
    \varPhi \approx \dis\frac{1}{2 \kappa^{*}} \left( \phi_{0}^{\eps}/\sqrt{\eps}, 
       \ \widetilde{K}_{q}^{o,*}(\phi_{0}^{\eps}/\sqrt{\eps}),
       \ \widetilde{K}_{r}^{o,*}(\phi_{0}^{\eps}/\sqrt{\eps}) \right)^{T}.
  \end{displaymath}
Here we used the fact that
  \begin{equation}
    \lim_{\eps \downarrow 0} \frac{1}{\sqrt{\eps}} \phi_{0}^{\eps} 
    = 2 \kappa^{*} \delta_{*} \ \ \ H^{-1}(I) \mbox{-sense}
  \end{equation}
(see Lemma 4.6).

\vspace{3mm}

Next, we solve (C.3). The first equation of (C.3) can be solved as
  \begin{displaymath}
    \begin{array}{rcl}
      p^{\#} &=& - (L^{\eps})^{-1} (q^{\#} + r^{\#})    \\
                \\
        &=& - \dis\frac{1}{\eps} \cdot 
              \dis \frac{\langle q^{\#} +r^{\#}, \ 
                 \phi_{0}^{\eps}/\sqrt{\eps} \rangle}
              {\zeta_{0}^{\eps}/\eps^{2}} \phi_{0}^{\eps}/\sqrt{\eps}
            - (L^{\eps})^{\dag}(q^{\#} + r^{\#}).
    \end{array}
  \end{displaymath}
Substituting this into the second and the third equation of (C.3), we obtain
  \begin{equation}
    \left \{
      \begin{array}{l}
         - M q^{\#} - \eps \alpha (L^{\eps})^{\dag} (q^{\#}) = 
         \alpha \dis \frac{\langle q^{\#} + r^{\#}, \ 
               \phi_{0}^{\eps}/\sqrt{\eps} \rangle}
              {\zeta_{0}^{\eps}/\eps^{2}} \phi_{0}^{\eps} /\sqrt{\eps}
          + \eps \alpha (L^{\eps})^{\dag} (r^{\#}), \\
                \\
         - N r^{\#} - \eps \alpha (L^{\eps})^{\dag} (r^{\#}) = 
         \beta \dis \frac{\langle q^{\#} + r^{\#}, \ 
               \phi_{0}^{\eps}/\sqrt{\eps} \rangle}
              {\zeta_{0}^{\eps}/\eps^{2}} \phi_{0}^{\eps} /\sqrt{\eps}
          + \eps \beta (L^{\eps})^{\dag} (q^{\#}).
      \end{array}
    \right .
  \end{equation}
For the odd symmetric case, solving the limiting system of (C.9) with respect to $(q^{\#}, r^{\#})$, we have
   \begin{equation}
    \left \{
    \begin{array}{l}
       q^{\# *} = 
       4 (\kappa^{*})^{2} \alpha \dis \frac{\langle q^{\# *}, \ \delta_{*} \rangle}
          {\hat{\zeta}_{0}^{*}} \widetilde{K}_{q}^{o,*}(\delta_{*}) + 
       4 (\kappa^{*})^{2} \alpha \dis \frac{\langle r^{\# *}, \ \delta_{*} \rangle}
          {\hat{\zeta}_{0}^{*}} \widetilde{K}_{q}^{o,*}(\delta_{*}),    \\
               \\
       r^{\# *} = 
       4 (\kappa^{*})^{2} \beta \dis \frac{\langle q^{\# *}, \ \delta_{*} \rangle}
          {\hat{\zeta}_{0}^{*}} \widetilde{K}_{r}^{o,*}(\delta_{*}) + 
       4 (\kappa^{*})^{2} \beta \dis \frac{\langle r^{\# *}, \ \delta_{*} \rangle}
          {\hat{\zeta}_{0}^{*}} \widetilde{K}_{r}^{o,*}(\delta_{*}).
    \end{array}
    \right.
  \end{equation}
We set $q^{\# *} = A \widetilde{K}_{q}^{o,*}(\delta_{*})$ and $r^{\# *} = B \widetilde{K}_{r}^{o,*}(\delta_{*})$, and substituting this into (C.10), we have
  \begin{displaymath}
    \left \{
    \begin{array}{l}
       A \widetilde{K}_{q}^{o,*}(\delta_{*}) 
       = \dis\frac{4 \alpha (\kappa^{*})^{2}}{\hat{\zeta}_{0}^{*}} [ 
          A \langle \widetilde{K}_{q}^{o,*}(\delta_{*}), \delta_{*} \rangle
        + B \langle \widetilde{K}_{r}^{o,*}(\delta_{*}), \delta_{*} \rangle ]
            \widetilde{K}_{q}^{o,*}(\delta_{*}),    \\
               \\
       B \widetilde{K}_{r}^{o,*}(\delta_{*}) 
       = \dis\frac{4 \beta (\kappa^{*})^{2}}{\hat{\zeta}_{0}^{*}} [ 
          A \langle \widetilde{K}_{q}^{o,*}(\delta_{*}), \delta_{*} \rangle
        + B \langle \widetilde{K}_{r}^{o,*}(\delta_{*}), \delta_{*} \rangle ]
            \widetilde{K}_{r}^{o,*}(\delta_{*}),
    \end{array}
    \right.
  \end{displaymath}
or the equivalent system
   \begin{equation}
    \left(
    \begin{array}{cc}
      \hat{\zeta}_{0}^{*}
      - 4 (\kappa^{*})^{2} \alpha \langle \widetilde{K}_{q}^{o,*}(\delta_{*}), 
                  \ \delta_{*} \rangle &
      - 4 (\kappa^{*})^{2} \alpha \langle \widetilde{K}_{r}^{o,*}(\delta_{*}), 
                  \ \delta_{*} \rangle \\
                \\
      - 4 (\kappa^{*})^{2} \beta \langle \widetilde{K}_{q}^{o,*}(\delta_{*}), 
                  \ \delta_{*} \rangle &
       \hat{\zeta}_{0}^{*}
      - 4 (\kappa^{*})^{2} \beta \langle \widetilde{K}_{r}^{o,*}(\delta_{*}), 
                  \ \delta_{*} \rangle
    \end{array}
    \right)
    \left(
    \begin{array}{c}
       A \\
         \\
       B
    \end{array}
    \right) = 
    \left(
    \begin{array}{c}
       0 \\
         \\
       0
    \end{array}
    \right)
  \end{equation}
We can easily see that the determinant of the matrix in (C.11) is
  \begin{displaymath}
    \hat{\zeta}_{0}^{*} \{ \hat{\zeta}_{0}^{*} - 4 (\kappa^{*})^{2} [
        \alpha \langle \widetilde{K}_{q}^{o,*}(\delta_{*}), 
               \ \delta_{*} \rangle
      + \beta \langle \widetilde{K}_{r}^{o,*}(\delta_{*}), 
               \ \delta_{*} \rangle ] \}
  \end{displaymath}
and it is zero (note the formula (4.19) of $\hat{\zeta}_{0}^{*}$ in Remark 4.9). Then we have $\beta A - \alpha B = 0$. On the other hand, for the symmetric case, the associated determinant is not zero and we obtain $A=B=0$.
Therefore $\varPhi^{\#}$ is an odd function. Now we choose $A = \alpha$ and $B = \beta$ respectively, and obtain
  \begin{displaymath}
    q^{\# *} = \alpha \widetilde{K}_{q}^{o,*}(\delta_{*}), \ \ \ 
    r^{\# *} = \beta \widetilde{K}_{r}^{o,*}(\delta_{*}).
  \end{displaymath}
Also we can see that the limit $\hat{p}^{\# *} := \dis\lim_{\eps \downarrow 0} \eps p^{\#}$ must exist. Then $\hat{p}^{\# *}$ satisfies
  \begin{displaymath}
    \hat{p}^{\# *} = -\dis \frac{4 (\kappa^{*})^{2}\langle q^{\# *} + r^{\# *},
          \ \delta_{*} \rangle} {\hat{\zeta}_{0}^{*}} \delta_{*}
         = - \delta_{*}.
  \end{displaymath}
Noting (C.8), we can see that the principal part of $\varPhi^{\# *}$ is given by
  \begin{displaymath}
    \varPhi^{\#} \approx 
      \dis\frac{1}{2 \kappa^{*}} \left (
             - \dis\frac{1}{\eps \sqrt{\eps}} \phi_{0}^{\eps}, \ 
        \alpha \widetilde{K}_{q}^{o,*}(\phi_{0}^{\eps}/\sqrt{\eps}), \ 
        \beta \widetilde{K}_{r}^{o,*}(\phi_{0}^{\eps}/\sqrt{\eps}) \right)^{T}.
  \end{displaymath}

Now we are ready to prove that $\langle \widehat{T}_{0} \varPhi, \varPhi^{\#} \rangle = 0$. In fact,
  \begin{equation}
    \begin{array}{l}
      \eps^{2} \langle \widehat{T}_{0} \varPhi, \varPhi^{\#} \rangle
      = - \dis\frac{1}{4 (\kappa^{*})^{2}}
        + \hat{\tau}_{0} \alpha 
        \langle (\widetilde{K}_{q}^{o,*})^{2}(\delta_{*}), \ \delta_{*} \rangle
        + \hat{\theta}_{0} \beta 
        \langle (\widetilde{K}_{r}^{o,*})^{2}(\delta_{*}), \ \delta_{*} \rangle   \\
                  \\ \hspace{5mm}
      = - \dis\frac{1}{4 (\kappa^{*})^{2}}
        - \left. \dis\frac{d}{d \hat{\lambda}^{*}} 
                 \right|_{\hat{\lambda}^{*}=0}
        [ \alpha \langle \widehat{K}_{q}^{o,*}(\delta_{*}), \delta_{*} \rangle
       + \beta \langle \widehat{K}_{r}^{o,*}(\delta_{*}), \delta_{*} \rangle ].
    \end{array}
  \end{equation}
Here we used the fact that $\widetilde{K}_{q}^{o,*}$ and $\widetilde{K}_{r}^{o,*}$ are self-adjoint operators and the following relations:
  \begin{equation}
    \left. \dis\frac{d}{d \hat{\lambda}^{*}} \right|_{\hat{\lambda}^{*}=0}
          \langle \widehat{K}_{q}^{o,*}(\delta_{*}), \delta_{*} \rangle
       = - \hat{\tau}_{0} 
         \langle (\widetilde{K}_{q}^{o,*})^{2}(\delta_{*}), \delta_{*} \rangle,
  \end{equation}
  \begin{equation}
    \left. \dis\frac{d}{d \hat{\lambda}^{*}} \right|_{\hat{\lambda}^{*}=0}
          \langle \widehat{K}_{r}^{o,*}(\delta_{*}), \delta_{*} \rangle
       = - \hat{\theta}_{0} 
         \langle (\widetilde{K}_{r}^{o,*})^{2}(\delta_{*}), \delta_{*} \rangle.
  \end{equation}
Because, let $q$ be a solution of
  \begin{displaymath}
    - \dis\frac{d^{2}}{dx^{2}} q + q + \hat{\tau}_{0} \hat{\lambda}^{*} q 
        = \delta_{*}
  \end{displaymath}
with Dirichlet boundary condition at $x=0$. Differentiating the above equation with respect to $\hat{\lambda}^{*}$, we have
  \begin{displaymath}
    - \dis\frac{d^{2}}{dx^{2}} 
      \left( \dis\frac{d q}{d \hat{\lambda}^{*}} \right)
       + \dis\frac{d q}{d \hat{\lambda}^{*}}
    + \hat{\tau}_{0} \hat{\lambda}^{*} \dis\frac{d q}{d \hat{\lambda}^{*}}
    + \hat{\tau}_{0} q = 0.
  \end{displaymath}
Since this implies $\widehat{T}_{q}( dq / d \hat{\lambda}^{*} ) = - \hat{\tau}_{0} q$, we obtain
  \begin{displaymath}
    \dis\frac{d}{d \hat{\lambda}^{*}} \widehat{K}_{q}^{o,*} (\delta_{*})
      = \dis\frac{d q}{d \hat{\lambda}^{*}} 
      = - \hat{\tau}_{0} \widehat{K}_{q}^{o,*} (q)
      = - \hat{\tau}_{0} (\widehat{K}_{q}^{o,*})^{2}(\delta_{*}).
  \end{displaymath}
Noting that $\widetilde{K}_{q}^{o,*} = \left. \widehat{K}_{q}^{o,*} \right|_{\hat{\lambda}^{*}=0}$, we obtain (C.13). (C.14) is proved similarly. On the other hand, differentiating $G_{od}(\hat{\lambda}^{*}; \hat{\tau}, \hat{\theta})$ with respect to $\hat{\lambda}^{*}$ at $(0, \hat{\tau}_{0}, \hat{\theta}_{0})$, we have
  \begin{displaymath}
   \dis\frac{\pd}{\pd \hat{\lambda}^{*}} G_{od}(0; \hat{\tau}_{0}, \hat{\theta}_{0})      = 1 + 4 (\kappa^{*})^{2} 
      \left. \dis\frac{d}{d \hat{\lambda}^{*}} \right|_{\hat{\lambda}^{*}=0}
        [ \alpha \langle \widehat{K}_{q}^{o,*}(\delta_{*}), \delta_{*} \rangle
       + \beta \langle \widehat{K}_{r}^{o,*}(\delta_{*}), \delta_{*} \rangle ].
  \end{displaymath}
This means that the principal part of $\langle \widehat{T}_{0} \varPhi, \varPhi^{\#} \rangle$ is equal to zero. 

Finally we will show that the algebraic multiplicity of the eigenvalue $\hat{\lambda}_{d}^{*}(\hat{\tau}_{0}, \hat{\theta}_{0})$ to the operator $\widehat{T}_{0}^{-1} {\mathcal L}$ is two. So we prove
  \begin{displaymath}
     \langle \widehat{T}_{0} \varPsi, \varPhi^{\#} \rangle \neq 0,
  \end{displaymath}
which means that the solvability condition for ${\mathcal L} \varXi = \widehat{T}_{0} \varPsi$ is not satisfied. Here we solve ${\mathcal L} \varPsi = \widehat{T}_{0} \varPhi$ by using SLEP method and compute the asymptotic form of $\varPsi$. Note that $\varPsi = (p,q,r)^T$ and $\varPhi = (\varphi_{p}, \varphi_{q}, \varphi_{r})^{T}$ satisfy the following system:
  \begin{equation}
    \left \{
      \begin{array}{l}
         L^{\eps} p - \eps \alpha q - \eps \beta r = \varphi_{p}, \\
                \\
          p + M q = \dis\frac{\hat{\tau}_{0}}{\eps^{2}} \varphi_{q},  \\
                \\
          p + N r = \dis\frac{\hat{\theta}_{0}}{\eps^{2}} \varphi_{r}.
      \end{array}
    \right .
  \end{equation}
The first equation of (C.15) can be solved as
  \begin{displaymath}
    \begin{array}{l}
      p = (L^{\eps})^{-1} (\eps \alpha q + \eps \beta r + \varphi_{p}) \\
                       \\ \hspace{2mm}
        = \dis \frac{\langle \alpha q + \beta r, 
                       \ \phi_{0}^{\eps}/\sqrt{\eps} \rangle}
              {\zeta_{0}^{\eps}/\eps^{2}} \phi_{0}^{\eps}/\sqrt{\eps}
            + \eps (L^{\eps})^{\dag}(\alpha q + \beta r)  \\
                       \\ \hspace{10mm}
           + \dis\frac{1}{\eps} \dis \frac{\langle \varphi_{p}, 
                       \ \phi_{0}^{\eps}/\sqrt{\eps} \rangle}
              {\zeta_{0}^{\eps}/\eps^{2}} \phi_{0}^{\eps}/\sqrt{\eps}
            + (L^{\eps})^{\dag}(\varphi_{p}).
    \end{array}
  \end{displaymath}
Substituting this into the second and third equations of (C.15), we obtain
  \begin{displaymath}
    \begin{array}{l}
         - M q - \eps \alpha (L^{\eps})^{\dag} (q) = 
         \dis \frac{\langle \alpha q + \beta r, 
              \ \phi_{0}^{\eps}/\sqrt{\eps} \rangle}
              {\zeta_{0}^{\eps}/\eps^{2}} \phi_{0}^{\eps} /\sqrt{\eps}
          + \eps \beta (L^{\eps})^{\dag} (r)  \\
               \\  \hspace{35mm}
          + \dis\frac{1}{\eps} \dis \frac{\langle \varphi_{p}, 
                       \ \phi_{0}^{\eps}/\sqrt{\eps} \rangle}
              {\zeta_{0}^{\eps}/\eps^{2}} \phi_{0}^{\eps}/\sqrt{\eps}
            + (L^{\eps})^{\dag}(\varphi_{p}) 
            - \dis\frac{\hat{\tau}_{0}}{\eps^{2}} \varphi_{q},
    \end{array}
  \end{displaymath}
  \begin{displaymath}
    \begin{array}{l}
         - N r - \eps \beta (L^{\eps})^{\dag} (r) = 
         \dis \frac{\langle \alpha q + \beta r, 
               \ \phi_{0}^{\eps}/\sqrt{\eps} \rangle}
              {\zeta_{0}^{\eps}/\eps^{2}} \phi_{0}^{\eps} /\sqrt{\eps}
          + \eps \alpha (L^{\eps})^{\dag} (q)  \\
                \\  \hspace{35mm}
          + \dis\frac{1}{\eps} \dis \frac{\langle \varphi_{p}, 
                       \ \phi_{0}^{\eps}/\sqrt{\eps} \rangle}
              {\zeta_{0}^{\eps}/\eps^{2}} \phi_{0}^{\eps}/\sqrt{\eps}
            + (L^{\eps})^{\dag}(\varphi_{p})
            - \dis\frac{\hat{\theta}_{0}}{\eps^{2}} \varphi_{r}.
    \end{array}
  \end{displaymath}
Here we introduce new variables as
  \begin{displaymath}
    \hat{q} := \eps^{2} q, \ \ \ , \hat{r} := \eps^{2} r,
  \end{displaymath}
we obtain a new system for $\hat{q}$ and $\hat{r}$:
 \begin{displaymath}
    \left \{
      \begin{array}{l}
         - M \hat{q} - \eps \alpha (L^{\eps})^{\dag} (\hat{q}) = 
         \dis \frac{\langle \alpha \hat{q} + \beta \hat{r}, 
              \ \phi_{0}^{\eps}/\sqrt{\eps} \rangle}
              {\zeta_{0}^{\eps}/\eps^{2}} \phi_{0}^{\eps} /\sqrt{\eps}
          + \eps \beta (L^{\eps})^{\dag} (\hat{r}) \\
                \\  \hspace{35mm}
          + \eps \dis \frac{\langle \varphi_{p}, 
                       \ \phi_{0}^{\eps}/\sqrt{\eps} \rangle}
              {\zeta_{0}^{\eps}/\eps^{2}} \phi_{0}^{\eps}/\sqrt{\eps}
            + \eps^{2} (L^{\eps})^{\dag}(\varphi_{p}) 
            - \hat{\tau}_{0} \varphi_{q}, \\
                \\
         - N \hat{r} - \eps \beta (L^{\eps})^{\dag} (\hat{r}) = 
         \dis \frac{\langle \alpha \hat{q} + \beta \hat{r}, 
               \ \phi_{0}^{\eps}/\sqrt{\eps} \rangle}
              {\zeta_{0}^{\eps}/\eps^{2}} \phi_{0}^{\eps} /\sqrt{\eps}
          + \eps \alpha (L^{\eps})^{\dag} (\hat{q}) \\
                \\  \hspace{35mm}
          + \eps \dis \frac{\langle \varphi_{p}, 
                       \ \phi_{0}^{\eps}/\sqrt{\eps} \rangle}
              {\zeta_{0}^{\eps}/\eps^{2}} \phi_{0}^{\eps}/\sqrt{\eps}
            + \eps^{2} (L^{\eps})^{\dag}(\varphi_{p})
            - \hat{\theta}_{0} \varphi_{r}.
      \end{array}
    \right .
  \end{displaymath}
Noting (C.8) and the fact
  \[
    \varphi_{p} \approx \dis\frac{1}{2 \kappa^{*}}
        \frac{\phi_{0}^{\eps}}{\sqrt{\eps}} \ \ \ \ \mbox{for small} \ \eps>0,
  \]
we obtain the limiting system as $\eps \downarrow 0$:
  \[
    \left \{
      \begin{array}{l}
         - M \hat{q} = 4 (\kappa^{*})^{2}
       \dis \frac{\langle \alpha \hat{q} + \beta \hat{r}, \ \delta_{*} \rangle}
              {\hat{\zeta}_{0}^{*}} \delta_{*}
            + \dis\frac{1}{\hat{\zeta}_{0}^{*}} \delta_{*}
            - \hat{\tau}_{0} \widetilde{K}_{q}^{o,*} (\delta_{*}), \\
                \\
         - N \hat{r} = 4 (\kappa^{*})^{2}
       \dis \frac{\langle \alpha \hat{q} + \beta \hat{r}, \ \delta_{*} \rangle}
              {\hat{\zeta}_{0}^{*}} \delta_{*}
            + \dis\frac{1}{\hat{\zeta}_{0}^{*}} \delta_{*}
            - \hat{\theta}_{0} \widetilde{K}_{r}^{o,*} (\delta_{*}).
      \end{array}
    \right .
  \]
Operating the operators $\widetilde{K}_{q}^{*}$ and $\widetilde{K}_{r}^{*}$ respectively, we have
  \[
    \left \{
    \begin{array}{l}
       \hat{q}^{*} =
    4 (\kappa^{*})^{2} \alpha \dis \frac{\langle \hat{q}^{*}, \ \delta_{*} \rangle}
          {\hat{\zeta}_{0}^{*}} \widetilde{K}_{q}^{o,*}(\delta_{*}) + 
    4 (\kappa^{*})^{2} \beta \dis \frac{\langle \hat{r}^{*}, \ \delta_{*} \rangle}
          {\hat{\zeta}_{0}^{*}} \widetilde{K}_{q}^{o,*}(\delta_{*})  \\
               \\ \hspace{20mm}
        + \dis\frac{1}{\hat{\zeta}_{0}^{*}} \widetilde{K}_{q}^{o,*}(\delta_{*})
        - \hat{\tau}_{0} (\widetilde{K}_{q}^{o,*})^{2} (\delta_{*}), \\
               \\
       \hat{r}^{*} =
    4 (\kappa^{*})^{2} \alpha \dis \frac{\langle \hat{q}^{*}, \ \delta_{*} \rangle}
          {\hat{\zeta}_{0}^{*}} \widetilde{K}_{r}^{o,*}(\delta_{*}) + 
    4 (\kappa^{*})^{2} \beta \dis \frac{\langle \hat{r}^{*}, \ \delta_{*} \rangle}
          {\hat{\zeta}_{0}^{*}} \widetilde{K}_{r}^{o,*}(\delta_{*}) \\
               \\ \hspace{20mm}
        + \dis\frac{1}{\hat{\zeta}_{0}^{*}} \widetilde{K}_{r}^{o,*}(\delta_{*})
        - \hat{\theta}_{0} (\widetilde{K}_{r}^{o,*})^{2} (\delta_{*}).
    \end{array}
    \right.    \\
  \]
We can see that the solutions $(q^{*}, r^{*})$ are given by
  \[
    \left(
    \begin{array}{c}
       \hat{q}^{*}  \\
              \\
       \hat{r}^{*}
    \end{array}
    \right) = c 
    \left(
    \begin{array}{c}
       \widetilde{K}_{q}^{o,*}(\delta_{*})  \\
              \\
       \widetilde{K}_{r}^{o,*}(\delta_{*})
    \end{array}
    \right) - 
    \left(
    \begin{array}{c}
       \hat{\tau}_{0}   (\widetilde{K}_{q}^{o,*})^{2} (\delta_{*})  \\
              \\
       \hat{\theta}_{0} (\widetilde{K}_{r}^{o,*})^{2} (\delta_{*})
    \end{array}
    \right)
  \]
for any $c$. So we chose the particular solutions $(\hat{q}^{*}, \hat{r}^{*})$ as
  \[
   \hat{q}^{*} = -\hat{\tau}_{0} (\widetilde{K}_{q}^{o,*})^{2} (\delta_{*}), \ \ \ 
   \hat{r}^{*} = -\hat{\theta}_{0} (\widetilde{K}_{r}^{o,*})^{2} (\delta_{*}).
  \]

On the other hand, $p$ is computed as follows:
  \begin{displaymath}
     p = \dis \frac{1}{\eps^{2} (\zeta_{0}^{\eps}/\eps^{2})}
         \left( \langle \alpha \hat{q} + \beta \hat{r}, 
                  \ \phi_{0}^{\eps}/\sqrt{\eps} \rangle
     + \dis\frac{1}{2 \kappa^{*}} \right) \phi_{0}^{\eps}/\sqrt{\eps} 
     + \dis\frac{1}{\eps} (L^{\eps})^{\dag}(\alpha \hat{q} + \beta \hat{r})
     + (L^{\eps})^{\dag}(\varphi_{p}).
  \end{displaymath}
Note that the limit of the coefficient of $O(\eps^{2})$-term is
  \begin{displaymath}
    \dis\frac{1}{\hat{\zeta}_{0}^{*}} \left(
    - 2 \kappa^{*} \hat{\tau}_{0} \alpha 
      \langle (\widetilde{K}_{q}^{o,*})^{2} (\delta_{*}), \ \delta_{*} \rangle 
    - 2 \kappa^{*} \hat{\theta}_{0} \beta 
      \langle (\widetilde{K}_{r}^{o,*})^{2} (\delta_{*}), \ \delta_{*} \rangle 
     + \dis\frac{1}{2 \kappa^{*}} \right),
  \end{displaymath}
which is equal to zero since that is equivalent to the drift bifurcation condition:
  \begin{displaymath}
    \dis\frac{d}{d \hat{\lambda}^{*}} G(0; \hat{\tau}_{0}, \hat{\theta}_{0}) = 0.
  \end{displaymath}
Therefore we can see that the principal part of $p$ is of $O(1/\eps)$. So we have
  \begin{displaymath}
    p \approx \dis\frac{1}{2 \eps}
        [ \hat{\tau}_{0} (\widetilde{K}_{q}^{o,*})^{2} (\delta_{*}) 
          + \hat{\theta}_{0} (\widetilde{K}_{r}^{o,*})^{2} (\delta_{*}) ].
  \end{displaymath}
Here we used the fact
  \[
    (L^{\eps} - \lambda)^{\dag}(h) \rightarrow 
            \dis\frac{h}{f_{u}^{*} - \lambda} \ \ \ (\eps \downarrow 0), \ \ \ 
     f_{u}^{*} \equiv -2.
  \]

Now we are ready to compute $\langle \widehat{T}_{0} \varPsi, \ \varPhi^{\#} \rangle$. Noting that
  \[
    \varPsi \approx \left(
       \dis\frac{1}{2 \eps} [ \hat{\tau}_{0} (\widetilde{K}_{q}^{o,*})^{2} (\delta_{*}) 
          + \hat{\theta}_{0} (\widetilde{K}_{r}^{o,*})^{2} (\delta_{*}) ], \ 
       - \dis\frac{1}{\eps^{2}} \hat{\tau}_{0} (\widetilde{K}_{q}^{o,*})^{2} (\delta_{*}), \ 
       - \dis\frac{1}{\eps^{2}} \hat{\theta}_{0} (\widetilde{K}_{r}^{o,*})^{2} (\delta_{*})
    \right ),
  \]
we obtain
  \begin{displaymath}
    \begin{array}{l}
      -2 \eps^{2} \langle \widehat{T}_{0} \varPsi, \varPhi^{\#} \rangle \approx
       \hat{\tau}_{0} \langle (\widetilde{K}_{q}^{o,*})^{2}(\delta_{*}), 
              \ \delta_{*} \rangle
       + \hat{\theta}_{0} \langle (\widetilde{K}_{r}^{o,*})^{2}(\delta_{*}), 
              \ \delta_{*} \rangle  \\
                  \\ \hspace{20mm}
       + 2 \alpha (\hat{\tau}_{0})^{2} 
           \langle (\widetilde{K}_{q}^{o,*})^{2}(\delta_{*}), 
              \ \widetilde{K}_{q}^{o,*}(\delta_{*}) \rangle
       + 2 \beta (\hat{\theta}_{0})^{2} 
           \langle (\widetilde{K}_{r}^{o,*})^{2}(\delta_{*}), 
              \ \widetilde{K}_{r}^{o,*}(\delta_{*}) \rangle  \\
                  \\ \hspace{10mm}
       = \hat{\tau}_{0} \langle (\widetilde{K}_{q}^{o,*})^{2}(\delta_{*}), 
              \ \delta_{*} \rangle
       + \hat{\theta}_{0} \langle (\widetilde{K}_{r}^{o,*})^{2}(\delta_{*}), 
              \ \delta_{*} \rangle
       + \dis\frac{\pd^{2}}{\pd \hat{\lambda}^{* 2}} 
       G(0; \hat{\tau}_{0}, \hat{\theta}_{0}) \\
                  \\ \hspace{10mm}
       > 0
    \end{array}
  \end{displaymath}
since $\widetilde{K}_{q}^{o,*}$ and $\widetilde{K}_{r}^{o,*}$ are self-adjoint operators. This completes the proof.
\hfill $\Box$

\begin{rmk}
It is shown here that the order of degeneracy of zero solution to the SLEP equation is equal to the algebraic multiplicity of the associated linearized eigenvalue problem for the drift bifurcation, which is consistent with the fact that the parabolic tangency of zero solution to the SLEP equation $G_{od}(\hat{\lambda}^{*}; \hat{\tau}, \hat{\theta}) = 0$.
The generalization of this correspondence to higher order multiplicity seems to be possible, however we defer the issue to one of the future works.
\end{rmk}

\section{Proof of Lemma 5.7}

(i) First, we prove $R(z; d)>0$ for $0 \leq z < \pi/4$. It is clear that $R(z; d)>0$ for $0 \leq Y + z < \pi/2$. For $Y + z \geq \pi/2$, it must be $Y > \pi/4$. Then, 
  \begin{displaymath}
    \begin{array}{rcl}
      \dis \frac{d}{\sqrt{\cos 2z}} R(z; d) &>& \cos z - e^{-X} \\
      &=& \sqrt{\dis\frac{Y^{2} + d^{2}}{2 Y^{2} + d^{2}} }
        - \dis\frac{1}{ \exp(\sqrt{Y^{2} + d^{2}}) }   \\
                    \\
      &>& \dis\frac{1}{\sqrt{2}} - \dis\frac{1}{ \exp(\sqrt{Y^{2} + d^{2}})} \\
                    \\
      &>&0
    \end{array}
  \end{displaymath}
since
  \begin{displaymath}
    \exp(\sqrt{Y^{2} + d^{2}}) > e^{Y} > e^{\pi/4}
  \end{displaymath}
and $e^{\pi} > 4$. 

Next we prove $I(z; d)<0$ for $0 < z < \pi/4$, we can see that $I(z; d)<0$ for $0 < Y + z < \pi$. When $Y + z \geq \pi$, it must be $Y > 3\pi/4$.
Then, 
  \begin{displaymath}
    \begin{array}{rcl}
      - \dis \frac{d}{\sqrt{\cos 2z}} I(z; d) &>& \sin z - e^{-X} \\
      &=& \dis\frac{Y}{ \sqrt{2 Y^{2} + d^{2}} }
        - \dis\frac{1}{ \exp(\sqrt{Y^{2} + d^{2}}) }   \\
                    \\
      &>& \dis\frac{3 \pi}{4\sqrt{2} \sqrt{Y^{2} + d^{2}}}
                - \dis\frac{1}{ e \sqrt{Y^{2} + d^{2}} } \\
                    \\
      &>&0
    \end{array}
  \end{displaymath}
since $e x \leq e^{x}$ for $x>0$ and $4 \sqrt{2} < 3 \pi e$. We easily check $I(0; d)=0$.

\vspace{3mm}

(ii) Differentiating $R(z;d)$ with respect to $z$, we have
  \begin{displaymath}
    \begin{array}{l}
      d \sqrt{\cos 2z} \ \dis\frac{dR}{dz}
        = - \sin 2z [\cos z + e^{-X} \cos (Y + z)] \\
                    \\ \hspace{24mm}
         - \cos 2z \ [ \sin z + e^{-X} 
              \{ X' \cos (Y + z) + (Y' + 1) \sin (Y + z) \} ]  \\
                        \\ \hspace{20mm}
      = - \sin 3z - e^{-X} \sin (Y + 3z) \\
                        \\ \hspace{24mm}
        - e^{-X} \cos 2z \ [ X' \cos (Y + z) + Y' \sin (Y + z) ]  \\
               \\ \hspace{20mm}
      = - \sin 3z - \dis\frac{1}{e^{X}} \sin (Y + 3z) 
        - \dis\frac{d}{e^{X} \sqrt{\cos 2z}} \sin (Y + 2z),
    \end{array}
  \end{displaymath}
where dash $'$ means the differentiation  with respect to $z$. Therefore we shall prove that
  \begin{equation}
    \sin 3z + \dis\frac{1}{e^{X}} \sin (Y + 3z) 
        + \dis\frac{d}{e^{X} \sqrt{\cos 2z}} \sin (Y + 2z)
  \end{equation}
is positive. Noting that $Y(z; d)$ is a monotonically increasing function with respect to $z$, we can see that each term of (D.1) is positive 
for $0 \leq z < \pi/2$ satisfying $Y(z; d) + 3z \leq \pi$.

Next we consider the case for $z$ such that $\pi < Y(z; d) + 3z$ and $Y(z; d) + 2z < \pi$. Then the third term of (D.1) is positive and $Y > \pi/4$ since $z < \pi/4$. Therefore, we shall show that the sum of the first and the second terms of (D.1) is positive when $\pi/4 \leq Y ( < \pi)$. Since
  \begin{displaymath}
    \begin{array}{l}
       \sin 3z = \sin z \cos 2z + \cos z \sin 2z  \\
               \\ \hspace{8mm}
       =\dis\frac{Y}{\sqrt{2 Y^{2} + d^{2}}} \dis\frac{d^{2}}{2 Y^{2}+d^{2}} +
       \sqrt{\dis\frac{Y^{2} + d^{2}}{2 Y^{2} + d^{2}}} 
              \dis\frac{2Y \sqrt{Y^{2} + d^{2}}}{2 Y^{2} + d^{2}}  \\
               \\ \hspace{8mm}
       = \dis\frac{2 Y^{3} + 3 d^{2} Y}{(2 Y^{2} + d^{2})^{3/2}},
    \end{array}
  \end{displaymath}
we obtain the following inequality:
  \begin{displaymath}
    \begin{array}{l}
      \sin 3z + \dis\frac{1}{e^{x}} \sin (Y + 3z) \geq 
        \dis\frac{Y (2 Y^{2} + 3 d^{2})}{(2 Y^{2} + d^{2})^{3/2}}
        - \dis\frac{1}{\exp(\sqrt{Y^{2}+d^{2}})}  \\
                  \\ \hspace{34mm}
      > \dis\frac{Y}{\sqrt{2} \sqrt{Y^{2} + d^{2}}}
        - \dis\frac{1}{\exp(\sqrt{Y^{2}+d^{2}})}  \\
                  \\ \hspace{34mm}
      > \dis\frac{\pi}{4 \sqrt{2} \sqrt{Y^{2} + d^{2}}}
        - \dis\frac{1}{e \sqrt{Y^{2}+d^{2}} }  \\
                  \\ \hspace{34mm}
      > 0.
    \end{array}
  \end{displaymath}
Here we used the fact that $e x \leq e^{x}$ for $x>0$ and $4 \sqrt{2} < \pi e$.

Finally we prove the case $Y(z; d) + 2z \geq \pi$. Then $Y \geq \pi/2$. Now it may be that the second and the third term are negative. Noting that
  \begin{displaymath}
    \begin{array}{l}
       {\rm (D.1)} = \dis\frac{Y(2 Y^{2} + 3 d^{2})}{(2 Y^{2} + d^{2})^{3/2}}
          - \dis\frac{1}{\exp(\sqrt{Y^{2}+d^{2}})}
          - \dis\frac{\sqrt{2 Y^{2} + d^{2}}}{\exp(\sqrt{Y^{2}+d^{2}})}  \\
                  \\ \hspace{5mm}
       = \dis\frac{Y(2 Y^{2} + 3 d^{2})}{(2 Y^{2} + d^{2})^{3/2}}
         - \dis\frac{1 + \sqrt{2 Y^{2} + d^{2}}}{\exp(\sqrt{Y^{2}+d^{2}})},
    \end{array}
  \end{displaymath}
we shall show 
  \begin{displaymath}
    \dis\frac{(2 y^{2} + d^{2})^{3/2} + (2 y^{2} + d^{2})^{2}}
        {y(2 y^{2} + 3 d^{2})} < \exp(\sqrt{y^{2}+d^{2}}).
  \end{displaymath}
for $Y>\pi/2$. Concerning the left hand side of the above inequality, we have
  \begin{displaymath}
    \begin{array}{l}
      \dis\frac{(2 Y^{2} + d^{2})^{3/2} + (2 Y^{2} + d^{2})^{2}}
        {Y(2 Y^{2} + 3 d^{2})} < 
      \dis\frac{(2 Y^{2} + 2 d^{2})^{3/2} + (2 Y^{2} + 2 d^{2})^{2}}
        {Y(2 Y^{2} + 2 d^{2})} \\
                \\ \hspace{43mm}
      = \dis\frac{\sqrt{2} (Y^{2} + d^{2})^{3/2} 
             + 2 (Y^{2} + d^{2})^{2}}
        {Y (Y^{2} + d^{2})}  \\
                \\ \hspace{43mm}
      < \dis\frac{2}{\pi} [\sqrt{2} \sqrt{Y^{2} + d^{2}} + 2 (Y^{2} + d^{2})].
    \end{array}
  \end{displaymath}
So we prove
  \begin{equation}
    \dis\frac{2}{\pi} [\sqrt{2} \sqrt{y^{2} + d^{2}} + 2 (y^{2} + d^{2})] 
        < \exp(\sqrt{y^{2}+d^{2}})
  \end{equation}
for $Y > \pi/2$. Let $b$ be $b = \sqrt{y^{2} + d^{2}}$, then the inequality we should prove is recast as
  \begin{displaymath}
    \dis\frac{2}{\pi} [\sqrt{2} b + 2 b^{2}] < e^{b}
  \end{displaymath}
Therefore we define $h_{2}(b)$ by
  \begin{displaymath}
    h_{2}(b) := \dis\frac{2}{\pi} \dis\frac{\sqrt{2} b + 2 b^{2}}{e^{b}}, 
  \end{displaymath}
and study the maximal value of $h_{2}(b)$ on 
$b \geq \sqrt{y^{2}+d^{2}} > \pi/2$. After some computation, we have
  \begin{displaymath}
    h_{2}'(b) = -\dis\frac{2}{\pi e^{b}} [2 b^{2} - (4-\sqrt{2}) b - \sqrt{2}].
  \end{displaymath}
This means that $h_{2}(b)$ has a maximal value at $b = \dis\frac{1}{2} (2 + \sqrt{2})$, and
  \begin{displaymath}
    h_{2} \left( \dis\frac{2 + \sqrt{2}}{2} \right) \approx 0.95183 < 1.
  \end{displaymath}
Thus we can conclude that (D.2) and $dR/dz \leq 0$ for $0 \leq z < \pi/4$. $\Box$

\bibliographystyle{unsrt} 
\bibliography{nishi_suzu}

\begin{thebibliography}{10}

\bibitem{L}
A.~W. Liehr.
\newblock {\em Dissipative Solitons in Reaction Diffusion Systems: Mechanisms,
  Dynamics, Interaction}.
\newblock Springer Series in Synergetics, Springer, 2013.

\bibitem{N2}
Y.~Nishiura.
\newblock {\em Far-from-Equilibrium Dynamics, Translations of Mathematical
  Monographs Iwanami Series in Modern Mathematics}.
\newblock AMS, 2002.

\bibitem{N3}
Y.~Nishiura.
\newblock Dynamics of particle patterns in dissipative systems -splitting 41
  destruction scattering-.
\newblock {\em SUGAKU EXPOSITIONS, American Mathematical Society},
  22(1):37--55, 2009.

\bibitem{NU2}
Y.~Nishiura and D.~Ueyama.
\newblock A skeleton structure of self-replicating dynamics.
\newblock {\em Physica D}, 130:73--104, 1999.

\bibitem{NU}
Y.~Nishiura and D.~Ueyama.
\newblock Spatio-temporal chaos for the gray-scott model.
\newblock {\em Physica D}, 150:137--162, 2001.

\bibitem{NTU2}
Y.~Nishiura, T.~Teramoto, and K.-I. Ueda.
\newblock Dynamic transitions through scatters in dissipative systems.
\newblock {\em Chaos}, 13(3):962--972, 2003.

\bibitem{NTU}
Y.~Nishiura, T.~Teramoto, and K.-I. Ueda.
\newblock Scattering and separators in dissipative systems.
\newblock {\em Phys. Rev. E}, 67:056210--1--056210--7, 2003.

\bibitem{NTU3}
Y.~Nishiura, T.~Teramoto, and K.-I. Ueda.
\newblock Scattering of traveling spots in dissipative systems.
\newblock {\em Chaos}, 15(4):047509, 2005.

\bibitem{VE}
V.~K. Vanag and I.~R. Epstein.
\newblock Localized patterns in reaction-diffusion systems.
\newblock {\em Chaos}, 17(3):037110, 2007.

\bibitem{BLSP}
M.~Bode, A.~W. Liehr, C.~P. Schenk, and H.-G. Purwins.
\newblock Interaction of dissipative solitons: particle-like behavior of
  localized structures in a three component reaction-diffusion system.
\newblock {\em Physica D}, 161:45--66, 2002.

\bibitem{OBSP}
M.~Or-Guil, M.~Bode, C.P. Schenk, and H.-G. Purwins.
\newblock Spot bifurcations in three-component reaction-diffusion systems: the
  onset of propagation.
\newblock {\em Phys. Rev. E}, 57:6432--6437, 1998.

\bibitem{PS}
H.-G. Purwins and L.~Stollenwerk.
\newblock Synergetic aspects of gas-discharge: lateral patterns in dc systems
  with a high ohmic barrier.
\newblock {\em Plasma Phys. Controlled Fusion}, 56(12):123001, 2014.

\bibitem{SOBP}
C.~P. Schenk, M.~Or-Guil, M.~Bode, and H.-G. Purwins.
\newblock Interacting pulses in three-component reaction-diffusion systems on
  two-dimensional domains.
\newblock {\em Phys. Rev. Lett.}, 78:3781--3784, 1997.

\bibitem{ME}
H.~Meinhardt.
\newblock {\em The Algorithmic Beauty of Sea Shells}.
\newblock Springer, 2009.

\bibitem{DVK}
A.~Doelman, P.~van Heijster, and T.~J. Kaper.
\newblock Pulse dynamics in a three-component system: existence analysis.
\newblock {\em J. Dyn. Diff. Eq.}, 21:73--115, 2009.

\bibitem{vHDK}
P.~van Heijster, A.~Doelman, and T.~J. Kaper.
\newblock Pulse dynamics in a three-component system: Stability and
  bifurcations.
\newblock {\em Physica D}, 237(24):3335--3368, 2008.

\bibitem{AGJ}
J.C. Alexander, R.A. Gardner, and C.K.R.T. Jones.
\newblock A topological invariant arising in the stability analysis of
  traveling waves.
\newblock {\em J. reine angew. Math.}, 410:167--212, 1990.

\bibitem{NF}
Y.~Nishiura and H.~Fujii.
\newblock Stability of singularly perturbed solutions to of reaction-diffusion
  equations.
\newblock {\em SIAM J. Math. Anal.}, 18:1726--1770, 1987.

\bibitem{NMIF}
Y.~Nishiura, M.~Mimura, H.~Ikeda, and H.~Fujii.
\newblock Singular limit analysis of stability of traveling wave solutions in
  bistable reaction-diffusion systems.
\newblock {\em SIAM J. Math. Anal.}, 21:85--122, 1990.

\bibitem{NS1}
Y.~Nishiura and H.~Suzuki.
\newblock Nonexistence of higher dimensional stable turing patterns in the
  singular limit.
\newblock {\em SIAM J. Math. Anal.}, 29:1087--1105, 1998.

\bibitem{NS2}
Y.~Nishiura and H.~Suzuki.
\newblock Higher dimensional slep equation and applications to morphological
  stability in polymer problems.
\newblock {\em SIAM J. Math. Anal.}, 36:916--966, 2004.

\bibitem{INS}
H.~Ikeda, Y.~Nishiura, and H.~Suzuki.
\newblock Stability of traveling waves and a relation between the evans
  function and the slep equation.
\newblock {\em J. reine angew. Math.}, 457:1--37, 1996.

\bibitem{BHIR}
M.~Chirilus-Bruckner, P.~van Heijster, H.~Ikeda, and J.D.M. Rademacher.
\newblock Unfolding symmetric bogdanov-takens bifurcations for front dynamics
  in a reaction-diffusion system.
\newblock {\em J. Nonlinear Sci.}, 29:2911--2953, 2019.

\bibitem{BDHR}
M.~Chirilus-Bruckner, A.~Doelman, P.~van Heijster, and J.D.M. Rademacher.
\newblock Butterfly catastrophe for fronts in a three-component
  reaction-diffusion system.
\newblock {\em J. Nonlinear Sci.}, 25:87--129, 2015.

\bibitem{vHCNT1}
P.~van Heijster, C.-N. Chen, Y.~Nishiura, and T.~Teramoto.
\newblock Localized patterns in a three-component fitzhugh-nagumo model
  revisited via an action functional.
\newblock {\em J. Dyn. Diff. Eq.}, 30:521--555, 2018.

\bibitem{vHCNT2}
P.~van Heijster, C.-N. Chen, Y.~Nishiura, and T.~Teramoto.
\newblock Pinned solutions in a heterogeneous three-component fitzhugh-nagumo
  model.
\newblock {\em J. Dyn. Diff. Eq.}, 31:153--203, 2019.

\bibitem{IN}
T.~Ikeda and Y.~Nishiura.
\newblock Pattern selection for two breathers.
\newblock {\em SIAM J. Appl. Math.}, 54:195--230, 1994.

\bibitem{NM}
Y.~Nishiura and M.~Mimura.
\newblock Layer oscillations in reaction-diffusion systems.
\newblock {\em SIAM J. Appl. Math.}, 49:481--514, 1989.

\bibitem{N}
Y.~Nishiura.
\newblock Coexistence of infinitely many stable solutions to reaction diffusion
  systems in the singular limit.
\newblock {\em Dynamics Reported (New Series)}, 3:25--103, 1994.

\bibitem{IMN}
H.~Ikeda, M.~Mimura, and Y.~Nishiura.
\newblock Global bifurcation phenomena of traveling wave solutions for some
  bistable reaction-diffusion systems.
\newblock {\em Nonlinear Anal.}, 13:507--526, 1989.

\bibitem{H}
D.~Henry.
\newblock {\em Geometric theory of semilinear parabolic equations}.
\newblock Springer-Verlag, 1981.

\bibitem{TUN}
T.~Teramoto, K.~Ueda, and Y.~Nishiura.
\newblock Phase-dependent output of scattering process for traveling breathers.
\newblock {\em Phys. Rev. E}, 69(4):056224--1--056224--8, 2004.

\bibitem{WIN2}
T.~Watanabe, M.~Iima, and Y.~Nishiura.
\newblock Spontaneous formation of travelling localized structures and their
  asymptotic behaviour in binary fluid convection.
\newblock {\em Journal of Fluid Mechanics}, 712:219--243, 12 2012.

\bibitem{WIN1}
T.~Watanabe, M.~Iima, and Y.~Nishiura.
\newblock A skeleton of collision dynamics - hierarchical network structure
  among even-symmetric steady pulses in binary fluid convection -.
\newblock {\em SIAM Journal on Applied Dynamical Systems}, 15(2):789--806,
  2016.

\bibitem{vHS1}
P.~van Heijster and B.~Sandstede.
\newblock Planar radial spots in a three-component fitzhugh-nagumo system.
\newblock {\em J. Nonlinear Science}, 21:705--745, 2011.

\bibitem{vHS2}
P.~van Heijster and B.~Sandstede.
\newblock Bifurcations to traveling planar spots in a three-component
  fitzhugh-nagumo system.
\newblock {\em Physica D}, 275:19--34, 2014.

\end{thebibliography}

\end{document}